\documentclass[matzei,final,numbook,envcountsame]{svjour}
\usepackage{amssymb,amsmath}
\begin{document}
\title{Motives and admissible representations 
of automorphism groups of fields}
\titlerunning{Motives and representations of automorphism groups of fields}
\author{M.Rovinsky\thanks{Supported in part by RFBR grant 02-01-22005}}
\institute{Independent University of Moscow, 
121002 Moscow, B.Vlasievsky Per. 11 
\email{marat@mccme.ru} and
Institute for Information Transmission Problems 
of Russian Academy of Sciences}
\date{Received: date / Revised version: date}
\maketitle
\begin{abstract}
Some of basic properties of the groups of automorphisms of 
algebraically closed fields and of their smooth representations 
are studied. In characteristic zero, Grothendieck motives modulo 
numerical equivalence are identified with a full subcategory 
in the category of graded smooth representations of certain 
automorphism groups of algebraically closed fields. 
\end{abstract}

\tableofcontents 

\section{Introduction}
Let $k$ be an algebraically closed field of characteristic zero. 
(This will be the case through out the paper, except Appendix 
\ref{subgr-p}, where all results of \S\ref{subgroups} are extended 
to positive characteristic.) Let $F$ be an algebraically closed 
extension of $k$ of transcendence degree $n$, $1\le n\le\infty$, 
and let $G=G_{F/k}$ be the group of automorphisms over $k$ of the 
field $F$. Let the set of subgroups $U_{k(x)}:={\rm Aut}(F/k(x))$ 
for all $x\in F$ be a base of neighbourhoods of the identity in $G$. 

This paper arose from an attempt to (i) compare properties of 
various ``geometric'' categories with properties of various 
categories of smooth (i.e., with open stabilizers) representations 
of $G$, and (ii) to find analogues for the group $G$ of familiar 
results of representation theory of $p$-adic groups. 

The group $G$ is very big, in particular, it contains the groups 
${\rm Aut}(L/k)$ as its sub-quotients for all sub-extensions 
$k\subset L\subset F$. All reduced irreducible algebraic groups 
of dimension $\le n$, the group ${\rm PGL}_{n+1}k$, some adelic 
groups are subgroups of groups of type ${\rm Aut}(L/k)$. 

One of the main results of the paper (Theorem \ref{exactness}) 
is the simplicity (in topological sense) of the group $G$ in the case 
$n=\infty$, and of the subgroup $G^{\circ}$ of $G$ generated by 
its compact subgroups in general. This is also true in positive 
characteristic, and implies (Corollary \ref{ex-fin-inf}) that in 
the case $n=\infty$ any non-trivial continuous representation of 
$G$ is faithful; and in the case $n<\infty$ any non-faithful 
continuous representation of $G$ factors through a discrete 
quotient $G/G^{\circ}$ of $G$. Another consequence is that $G^{\circ}$ 
(and $G$ if $n=\infty$) admits no smooth representations of finite degree. 

Unfortunately, I do not know much about the group $G/G^{\circ}$, 
and this is one of the reasons why I prefer to work 
in the ``stable'' case $n=\infty$. 

There is an evident link between representations of $G$ and some 
geometric objects. Namely, for a scheme $X$ over $k$ there is a 
natural smooth $G$-action on the group of cycles on $X_F:=X\times_kF$. 
Conversely, any smooth cyclic $G$-module is a quotient of the 
$G$-module of ``generic'' 0-cycles on $X_F$ for an appropriate 
irreducible variety $X$ of dimension $\le n$ over $k$. 
The Hecke algebras, playing an important r\^{o}le in representation 
theory of locally compact groups, become in our case algebras of 
non-degenerate correspondences on certain varieties over $k$, 
cf. \S\ref{Hecke-zykly}, p.\pageref{Hecke-zykly}.

In some cases one can identify the groups of morphisms between 
geometric objects with the groups of morphisms between 
corresponding $G$-modules (cf. Propositions \ref{geomod} 
and \ref{g-g-knought}, and Corollary \ref{Hom-G-mod}). 

(Homological) Grothendieck motives are pairs $(X,\pi)$ 
consisting of a smooth projective variety $X$ over $k$ with 
irreducible components $X_j$ and a projector $\pi=\pi^2\in
\bigoplus_jA^{\dim X_j}(X_j\times_kX_j)$ in the algebra 
of correspondences on $X$ modulo an adequate equivalence relation. 
The morphisms are defined by ${\rm Hom}((X',\pi'),(X,\pi))=
\bigoplus_{i,j}\pi_j\cdot A^{\dim X_j}(X_j\times_kX'_i)\cdot\pi'_i$. 
The category of Grothendieck motives carries an additive and a tensor 
structures: \begin{multline*}(X',\pi')\bigoplus(X,\pi):=
(X'\coprod X,\pi'\oplus\pi),\\ (X',\pi')\otimes(X,\pi)
:=(X'\times_kX,\pi'\times_k\pi).\end{multline*} 
A {\sl primitive $q$-motive} is a pair $(X,\pi)$ as above 
with $\dim X=q$ and $\pi\cdot A^q(X\times_kY\times{\mathbb P}^1)=0$ 
for any smooth projective variety $Y$ over $k$ with $\dim Y<q$. 
For instance, the category of the primitive $1$-motives modulo 
numerical equivalence is equivalent to the category of abelian 
varieties over $k$ with morphisms tensored with ${\mathbb Q}$. 
If the adequate equivalence relation is numerical equivalence, 
it follows from a result of Jannsen \cite{jan} that any Grothendieck 
motive is semi-simple and admits ``primitive'' decomposition 
$\bigoplus_{i,j}M_{ij}\otimes{\mathbb L}^{\otimes i}$, where 
$M_{ij}$ is a primitive $j$-motive and 
${\mathbb L}=({\mathbb P}^1,{\mathbb P}^1\times\{0\})$ 
(see Remark on p.\pageref{proof-of-prim-dec}). 

The major part of the results of 
Section \ref{Sgr} can be summarized as follows. 
\begin{theorem} \label{summ} \begin{enumerate} 
\item If $n=\infty$ there is a fully faithful functor 
${\mathbb B}^{\bullet}$: 
$$\left\{\begin{array}{c} \mbox{{\rm motives over $k$ modulo}}\\
\mbox{{\rm numerical equivalence}}\end{array}\right\}
\stackrel{{\mathbb B}^{\bullet}}{\longrightarrow}
\left\{\begin{array}{c} \mbox{{\rm graded semi-simple admissible}}\\
\mbox{$G$-{\rm modules of finite type}}\end{array}\right\}.$$
The grading corresponds to powers of the motive 
${\mathbb L}$ in the ``primitive'' decomposition above. 
\item \label{summ-2} For any $1\le n\le\infty$ and each $q\ge 0$ 
there is a functor ${\mathfrak B}^q$, fully faithful if $q\le n$: 
$$\left\{\begin{array}{c} \mbox{{\rm primitive $q$-motives over $k$}}\\
\mbox{{\rm modulo numerical equivalence}}\end{array}\right\}
\stackrel{{\mathfrak B}^q}{\longrightarrow}
\left\{\begin{array}{c} \mbox{{\rm semi-simple admissible}}\\ 
\mbox{$G$-{\rm modules of finite type}}\\ \mbox{{\rm and of level} 
$q$}\end{array}\right\}.$$ (One says that a semi-simple admissible 
$G$-module $W$ is of level $q$ if $N_qW=W$ and $N_{q-1}W=0$ for the 
filtration $N_{\bullet}$ defined in the beginning of \S\ref{func-I}.) 
\item \label{phora} If $n<\infty$ then the $G$-module 
${\mathfrak B}^n(M)$ carries a bilinear symmetric non-degenerate 
$G$-equivariant form with values 
in an oriented $G$-module ${\mathbb Q}(\chi)$ of degree 1, where 
$M=(X,\Delta_{k(X)})$ is the maximal primitive $n$-submotive 
of the motive $(X,\Delta_X)$ and $\dim X=n$. 

This form is definite if, for $(n-1)$-cycles on $2n$-dimensional 
complex varieties, numerical equivalence coincides with homological 
(e.g., for $n\le 2$), and therefore, ${\mathfrak B}^n$ factors 
through the subcategory of ``polarizable'' $G$-modules (i.e., 
carrying a positive form as above). \end{enumerate} \end{theorem}
This is a direct consequence of Corollary \ref{Hom-G-mod} and 
Propositions \ref{supply}, \ref{b-explic}, \ref{polar}. Roughly 
speaking, the functors ${\mathfrak B}^q$ and ${\mathbb B}^{\bullet}
=\oplus^{{\rm graded}}_j{\mathbb B}^{[j]}$ are defined as 
spaces of 0-cycles defined over $F$ modulo ``numerical 
equivalence over $k$''. Details are in \S\ref{functors-B}, 
p.\pageref{functors-B}, where it is shown that they 
are pro-representable. It follows from Proposition 
\ref{max-prim} that ${\mathfrak B}^q((X,\pi))$ depends only on the 
birational class of $X$. Moreover, the functor ${\mathfrak B}^1$ of 
Theorem \ref{summ}(\ref{summ-2}) is an equivalence of categories when 
$n=\infty$, cf. \S\ref{level-1}, and by Corollary \ref{for-fai}, the 
composition of the functor ${\mathfrak B}^1$ with the forgetful 
functor to the category of $G^{\circ}$-modules is also fully faithful. 

\begin{conjecture} \label{conj-equi} The functor 
${\mathbb B}^{\bullet}$ is an equivalence of categories. 
Equivalently, for any $q\ge 0$ the functor ${\mathfrak B}^q$ 
is an equivalence of categories if $n=\infty$. \end{conjecture} 

The section is concluded by showing (Corollary \ref{no-fi-di}) 
that any polarizable representation (in the sense of 
Theorem \ref{summ}(\ref{phora})) is infinite-dimensional. 
This is deduced from a vanishing result (Corollary \ref{hyperspe}) 
for representations of the Hecke algebra of the subgroup 
${\rm Gal}(F/L(x))$, induced by polarizable representations of $G$ 
(here $L$ is a subextension of $k$ in $F$ of finite type and $x$ an 
element of $F$ transcendental over $L$ with $F=\overline{L(x)}$) 
corresponding to the triviality of the primitive $n$-submotives 
of $(Y\times{\mathbb P}^1,\pi)$, where $\dim Y<n$. 

Corollary \ref{hyperspe} suggests also that the category of the 
primitive $n$-motives is not too far from the category of the 
polarizable $G$-modules (at least, if $n\le 2$). However, as the 
twists of the polarizables by order-two characters of $G$ are again 
polarizable, but not ``motivic'' (cf. Corollary \ref{0-1}), one 
should impose some additional conditions. One of such conditions 
could be the ``stability'' in the sense that for any algebraically 
closed extension $F'$ of $F$ any ``motivic'' representation is 
isomorphic to $W^{G_{F'/F}}$ for some smooth $G_{F'/k}$-module $W$; 
another one could be the ``arithmeticity'' in the sense that any 
``motivic'' representation admits an extension $W$ by a module of 
lower level (in the sense of filtration $N_{\bullet}$) with $W$ 
isomorphic to the restriction to $G\subseteq G_{F/k'}$ of a smooth 
$G_{F/k'}$-module for a subfield $k'\subseteq k$ of finite type 
over ${\mathbb Q}$... 

By analogy with the Langlands correspondences, one can 
call the $G$-modules in the image of ${\mathfrak B}^n$ {\sl cuspidal}. 
For the groups ${\rm GL}$ over a local non-archimedian field 
there are several equivalent definitions of quasicuspidal 
representations. One of them: all matrix coefficients 
(the functions $\langle\sigma w,\widetilde{w}\rangle$, cf. 
p.\pageref{geomor}) are compactly supported modulo the center. 
However, it is shown in \S\ref{geomor} that there are no 
such representations for any subgroup of $G$ containing $G^{\circ}$. 
This is a consequence of the irreducibility of the smooth 
$G$-modules $F/k$ and $F^{\times}/k^{\times}$, considered 
as modules over the subgroup $G^{\circ}$ of $G$, and their 
faithfulness as modules over the algebra of compactly supported 
measures on $G$ shown in \S\ref{muladd}. 

In \S\ref{ext-examp}, in the case $n=\infty$, various analogues 
of Hilbert Theorem 90 are verified. In particular, it is shown 
in Proposition \ref{torsors}, that any $G$-torsor under 
${\mathcal A}(F)$ is trivial for any algebraic group 
${\mathcal A}$ over $k$. There exist, however, 
interesting examples of torsors in the case $n<\infty$. 

If $n=\infty$ and ${\mathcal A}$ is an irreducible commutative 
algebraic group over $k$, we show in Corollary \ref{Ext-1} that 
${\rm Ext}^1_{{\mathcal S}m_G}({\mathcal A}(F)/{\mathcal A}(k),
{\mathbb Q})={\rm Hom}({\mathcal A}(k),{\mathbb Q})$, where 
${\mathcal S}m_G$ is the category of smooth $G$-modules. If 
${\mathcal A}$ is an abelian variety then ${\mathcal A}(F)/
{\mathcal A}(k)={\mathfrak B}^1({\mathcal A}^{\vee})$ (here 
${\mathcal A}^{\vee}:={\rm Pic}^{\circ}{\mathcal A}$ is the dual 
abelian variety), so this should correspond to the identity 
${\rm Ext}^1_{\mathcal{MM}_k}({\mathbb Q}(0),H_1({\mathcal A}))
={\mathcal A}(k)_{{\mathbb Q}}$ in the category of mixed motives 
over $k$. If ${\mathcal A}={\mathbb G}_m$ then the identity 
${\rm Ext}^1_{\mathcal{MM}_k}({\mathbb Q}(0),{\mathbb Q}(1))=
k^{\times}\otimes{\mathbb Q}$ suggests that the non-admissible 
$G$-module $F^{\times}/k^{\times}$ admits a motivic 
interpretation analogous to ${\mathbb Q}(1)$. 

The purpose of \S\ref{I-G} is to introduce an abelian category 
${\mathcal I}_G$ of ``homotopy invariant'' representations having 
some properties of the Chow groups. If $n=\infty$ then it 
contains all admissible $G$-modules (Proposition \ref{adm-i}), 
it is a Serre subcategory in ${\mathcal S}m_G$ (Proposition 
\ref{thick}), and it is closed under the inner ${\mathcal H}om$ 
functor on ${\mathcal S}m_G$ (Proposition \ref{inner-hom}). 
There are no smooth projective representations of $G$, 
if $n=\infty$ (cf. Remark on p.\pageref{no-proj}). However, it 
is shown in Corollary \ref{projectivity} that ${\mathcal I}_G$ 
has enough projective objects. Namely, the inclusion functor 
${\mathcal I}_G\longrightarrow{\mathcal S}m_G$ admits the left adjoint 
${\mathcal S}m_G\stackrel{{\mathcal I}}{\longrightarrow}{\mathcal I}_G$, 
and to any subextension $L$ of finite type one associates the projective 
object $C_L:={\mathcal I}{\mathbb Q}[G/G_{F/L}]$ of ${\mathcal I}_G$. 
For any smooth proper model $X$ of $L/k$ there is a natural surjection 
$C_L\longrightarrow CH_0(X\times_kF)_{{\mathbb Q}}$. One can expect 
(Conjecture \ref{C=CH}) that this is an isomorphism if $n=\infty$. 
If ${\rm tr.deg}(L/k)=1$ this is Corollary \ref{cl-pic}. 

At the end of \S\ref{level-1}, p.\pageref{filtra-2}, a functorial 
decreasing filtration ${\mathcal F}^{\bullet}$ on objects of 
${\mathcal I}_G$ is introduced. It is likely that in the case of 
$G$-modules of type $CH_0(X_F)_{{\mathbb Q}}$ for smooth proper 
$X$ over $k$ it is the motivic filtration. This agrees with 
Corollary \ref{filt-2}: $C_{k(X)}\cong{\mathbb Q}\oplus
{\rm Alb}X(F)_{{\mathbb Q}}\oplus{\mathcal F}^2C_{k(X)}.$ 

If $n=\infty$ then Conjecture \ref{C=CH}, the semi-simplicity 
conjecture and Bloch--Beilinson filtration conjecture would imply 
(i) Conjecture \ref{conj-equi}, (ii) that any irreducible object of 
${\mathcal I}_G$ is admissible, and (iii) that the $G$-modules $gr^N_jW$ 
are semi-simple for any object $W$ of ${\mathcal I}_G$ (the latter is 
Conjecture \ref{ss-conj}). Indeed, for some collection of subfields 
$L\subset F$ of finite type and of transcendence degree $j$ over $k$ 
there is a surjective morphism $\oplus_L{\mathbb Q}[G/G_{F/L}]
\stackrel{\xi}{\longrightarrow}gr^N_jW$, which factors through 
$\oplus_Lgr^N_jC_L$, cf. Proposition \ref{def-I}. If for a smooth 
proper model $Y_{[L]}$ of $L/k$ one has 
$C_L=CH_0(Y_{[L]}\times_kF)_{{\mathbb Q}}$ 
then $gr^N_jC_L=CH^j(L\otimes_kF)_{{\mathbb Q}}$, 
so $\xi$ factors through $\oplus_LCH^j(L\otimes_kF)_{{\mathbb Q}}$. 
One can deduce from the semi-simplicity conjecture and the 
filtration conjecture (cf. \cite{n-mot} \S1.4, or \cite{mpi} 
Prop.1.1.1)\footnote{where it is shown that under above assumptions 
the localization surjection $CH^{\ast}(Y_{[L]}\times_kY_{[L']})
\longrightarrow CH^{\ast}(L\otimes_kL')_{{\mathbb Q}}$ kills the 
numerically trivial cycles.} that $CH^j(L\otimes_kF)_{{\mathbb Q}}$ 
coincides with ${\mathfrak B}^j(M)$, where $M$ is the maximal 
primitive $j$-submotive of the motive $(Y_{[L]},\Delta_{Y_{[L]}})$. 
Finally, by semi-simplicity, there are projectors $\pi_L$ 
and an isomorphism $$\oplus_L{\mathfrak B}^j((Y_{[L]},\pi_L))
\stackrel{\sim}{\longrightarrow}gr^N_jW.$$ This shows 
(iii), and taking irreducible $W$ (which coincides with 
$gr^N_jW$ for some $j$) we get also (i) and (ii). 

It is also conjectured\footnote{and deduced from Conjecture 
\ref{ss-conj}} that the level filtration $N_{\bullet}$ is strictly 
compatible with the morphisms in ${\mathcal I}_G$ (cf. Corollary 
\ref{str-comp}), so that, in particular, extensions of $G$-modules 
in ${\mathcal I}_G$ of lower level by irreducible $G$-modules in 
${\mathcal I}_G$ of higher level are (canonically) split. Obviously, 
this is motivated by Hodge theory, and one would like to find a 
category bigger than ${\mathcal I}_G$ and modify the filtration 
to keep this property. 

However, if we want to consider $G$-modules like 
$F^{\times}/k^{\times}$, the notion of weight should be more subtle. 
Usually, for a pair $W_1,W_2$ of irreducible objects, $W_1$ is of 
higher weight if ${\rm Ext}^1(W_1,W_2)\neq 0$, so 
this would give weight$({\mathbb Q})<$weight$(F^{\times}/k^{\times})<$%
weight$({\mathcal A}(F)/{\mathcal A}(k))$, cf. \S\ref{ext-examp}, for 
any abelian variety ${\mathcal A}$ over $k$, which is not good if 
${\mathcal A}(F)/{\mathcal A}(k)$ corresponds to $H_1({\mathcal A})$. 

If $n=\infty$, the category ${\mathcal I}_G$ carries a tensor 
structure compatible with the inner ${\mathcal H}om$, but its 
associativity depends on Conjecture \ref{C=CH}. 

If $n<\infty$ then the category of smooth $G$-modules has sufficiently 
many projective objects. Namely, any smooth $G$-module is a quotient 
of a direct sum of ${\mathbb Q}[G/U_j]$ for some open compact subgroups 
$U_j$ of $G$. However, the $G$-modules ${\mathbb Q}[G/U]$ seem to be very 
complicated. The last section contains two examples of pairs of essentially 
different open compact subgroups $U_1$ and $U_2$ of $G$ with the same 
irreducible subquotients of ${\mathbb Q}[G/U_1]$ and ${\mathbb Q}[G/U_2]$. 
As in both examples the primitive motives of maximal level of models of 
$F^{U_1}$ and $F^{U_2}$ are trivial, one could expect that collections 
of irreducible subquotients of ${\mathbb Q}[G/U]$ are of motivic nature. 

In Appendix \ref{Hecke-center} one shows that the centers of 
the Hecke algebras of the pairs $(G,U)$ and $(G^{\circ},U)$ 
(see \S\ref{NCT} for the definition) consist of scalars for 
any compact subgroup $U$ in $G$. Compared to the analogous 
question for $p$-adic groups, this is a negative result. 

\subsection{Notations, conventions and terminology} \label{NCT} 
For a field $F$ and a collection of its subrings 
$F_0,(F_{\alpha})_{\alpha\in I}$ we denote by 
$G_{\{F,(F_{\alpha})_{\alpha\in I}\}/F_0}$ the group of automorphisms 
of the field $F$ over $F_0$ preserving all $F_{\alpha}$, and set 
$G_{F/F_0}:=G_{\{F\}/F_0}$. If $K$ is a subfield of $F$ then 
$\overline{K}$ denotes its algebraic closure in $F$, 
${\rm tr.deg}(F/K)$ the transcendence degree of the extension $F/K$ 
(possibly infinite, but countable), and $U_K$ denotes the group 
$G_{F/K}$. Throughout the paper $k$ is an algebraically closed 
field, $F$ its algebraically closed extension with 
${\rm tr.deg}(F/k)=n\ge 1$ and $G=G_{F/k}$. Everywhere, except 
the appendix, $k$ is of characteristic zero. 

For a totally disconnected topological group $H$ we denote by 
$H^{\circ}$ its subgroup generated by the compact subgroups. 
Obviously, $H^{\circ}$ is a normal subgroup in $H$, which is 
open at least if $H$ is locally compact. 

In what follows, ${\mathbb Q}$ is the field of rational numbers, 
and a {\sl module} is always a ${\mathbb Q}$-vector space. 
For an abelian group $A$ set $A_{{\mathbb Q}}=A\otimes{\mathbb Q}$. 

A representation of $H$ in a vector space $W$ over a field is 
called {\sl smooth}, if stabilizers of all vectors in $W$ are 
open. A smooth representation $W$ is called {\sl admissible} if fixed 
vectors of each open subgroup form a finite-dimensional subspace in $W$. 
A representation of $H$ in $W$ is called {\sl continuous}, if 
stabilizers of all vectors in $W$ are closed. Any cyclic $G$-module 
is a quotient of the continuous $G$-module ${\mathbb Q}[G]$, so any 
$G$-module is a quotient of a continuous $G$-module. 

${\mathbb Q}(\chi)$ is the quotient of the free abelian group generated 
by the set of compact open subgroups in $G^{\circ}$ by the relations 
$[U]=[U:U']\cdot[U']$ for all $U'\subset U$. If $n<\infty$ 
it is a one-dimensional ${\mathbb Q}$-vector space oriented 
by $[U]>0$ for any $U$. The group $G$ acts on it by conjugation. 
$\chi:G\longrightarrow{\mathbb Q}^{\times}_+$ is the modulus. 

${\bf D}_E(H):=\!\lim\limits_{\longleftarrow_U}\!E[H/U]$ and 
$\widehat{H}:=\lim\limits_{\longleftarrow_U}H/U$, where, 
for a field $E$ of characteristic zero, the inverse 
systems are formed with respect to the projections 
$E[H/V]\stackrel{r_{VU}}{\longrightarrow}E[H/U]$ and 
$H/V\stackrel{r_{VU}}{\longrightarrow}H/U$ induced by inclusions 
$V\subset U$ of open subgroups in $H$. $\widehat{H}$ is a semigroup.  
For any $\nu\in{\bf D}_E(H)$, any $\sigma\in H$ and an open subgroup 
$U$ we set $$\nu(\sigma U):=\mbox{coefficient of $[\sigma U]$ of the 
image of $\nu$ in $E[H/U]$}.$$ The {\sl support} of $\nu$ is the 
minimal closed subset $S$ in $\widehat{H}$ such that $\nu(\sigma U)=0$ 
if $\sigma U\bigcap S=\emptyset$. Define a pairing ${\bf D}_E(H)\times 
W\longrightarrow W$ for each smooth $E$-representation $W$ of $H$ by 
$(\nu,w)\longmapsto\sum_{\sigma\in H/V}\nu(\sigma V)\cdot\sigma w$, 
where $V$ is an arbitrary open subgroup in the stabilizer of $w$. When 
$W=E[H/U]$ this pairing is compatible with the projections $r_{VU}$, 
so we get a pairing ${\bf D}_E(H)\times\lim\limits_{\longleftarrow_U}
\!E[H/U]\longrightarrow\lim\limits_{\longleftarrow_U}\!E[H/U]=
{\bf D}_E(H)$, and thus an associative multiplication ${\bf D}_E(H)
\times{\bf D}_E(H)\stackrel{\ast}{\longrightarrow}{\bf D}_E(H)$ 
extending the convolution of compactly supported measures. Set 
${\bf D}_E={\bf D}_E(G)$ and ${\bf D}^{\circ}_E={\bf D}_E(G^{\circ})$. 

\label{def-Hecke}
If $U$ is a compact subgroup in $H$ the {\sl Hecke algebra} of the 
pair $(H,U)$ is the subalgebra ${\mathcal H}_E(H,U):=
h_U\ast{\bf D}_E(H)\ast h_U$ in ${\bf D}_E(H)$ of $U$-bi-invariant 
measures. Here $h_U$ is the Haar measure on $U$ defined by the system 
$(h_U)_V=[U:U\bigcap V]^{-1}\sum_{\sigma\in U/U\bigcap V}[\sigma V]
\in{\mathbb Q}[H/V]$ for all open subgroups $V\subset H$. $h_U$ is the 
identity in ${\mathcal H}_E(H,U)$ and $h_Uh_{U'}=h_U$ for a closed 
subgroup $U'\subseteq U$. Set ${\mathcal H}_E(U)=
{\mathcal H}_E(G,U)$, ${\mathcal H}(U)={\mathcal H}_{{\mathbb Q}}(G,U)$ 
and ${\mathcal H}^{\circ}_E(U)={\mathcal H}_E(G^{\circ},U)$. 

For any variety $X$ over $k$ and any field extension $E/k$ we set 
$X_E:=X\times_kE$, and denote by $\widetilde{X}$ one of its 
desingularizations. For an extension $L/k$ of finite type, 
$Y_{[L]}=Y_{U_L}$ denotes a smooth proper model of $L/k$. 

${\mathbb P}^M_K$ denotes the $M$-dimensional projective space 
over a field $K$. 

For a commutative group scheme ${\mathcal A}$ over $k$ 
we set $W_{{\mathcal A}}={\mathcal A}(F)/{\mathcal A}(k)$. 

${\mathcal S}m_H(E)$ is the category of smooth $E$-representations 
of $H$. ${\mathcal I}_G(E)$ is the full subcategory in ${\mathcal S}m_G(E)$ 
consisting of those representations $W$ of $G$ for which 
$W^{G_{F/L}}=W^{G_{F/L'}}$ for any extension $L$ of $k$ 
in $F$ and any purely transcendental extension $L'$ of $L$ in $F$. 
When discussing ${\mathcal I}_G(E)$, the principal case will be 
$n=\infty$. We set ${\mathcal S}m_H={\mathcal S}m_H({\mathbb Q})$ 
and ${\mathcal I}_G={\mathcal I}_G({\mathbb Q})$. 
The {\sl level} filtration $N_{\bullet}$ is defined in 
the beginning of \S\ref{func-I}. 

\section{Preliminaries on closed subgroups in $G$} \label{subgroups} 
The topology on $G$ described in Introduction has been studied 
in \cite{jac}, p.151, Exercise 5, \cite{p-s-sha}, \cite{shimura} 
Ch.6, \S6.3, and \cite{ihara} Ch.2, Part 1, Section 1. It is shown 
that $G$ is Hausdorff, locally compact if $n<\infty$, and totally 
disconnected; the subgroups $G_{\{F,(F_{\alpha})_{\alpha\in I}\}/k}$ 
are closed in $G$, there is an injective morphism of unitary semigroups 
$$\{\mbox{{\rm subfields in} $F$ {\rm over} $k$}\}
\longrightarrow\{\mbox{{\rm closed subgroups in} $G$}\}$$ 
given by $K\longmapsto{\rm Aut}(F/K)$, its image is stable under 
passages to sub-/sup-groups with compact quotients, and it induces 
bijections 
\begin{itemize} \item $\{\mbox{{\rm subfields} 
$K\subset F$ {\rm over} $k$ {\rm with} $F=\overline{K}$}\}
\leftrightarrow\{\mbox{{\rm compact subgroups of} $G$}\}$; 
\item $\left\{\begin{array}{c} \mbox{{\rm subfields} $K$ {\rm of} $F$ 
{\rm of finite type}} \\ \mbox{{\rm over} $k$ {\rm with} 
$F=\overline{K}$} \end{array} \right\}\leftrightarrow
\left\{\begin{array}{c} \mbox{{\rm compact open}} \\ 
\mbox{{\rm subgroups of} $G$}\end{array} \right\}$. 

The inverse correspondences are given by $G\supset H\longmapsto F^H$
(the subfield in $F$ fixed by $H$). \end{itemize} 
The first is \cite{p-s-sha}, \S3, Lemma 1, or \cite{shimura} 
Prop.6.11; the second is immediate from loc.cit., or \cite{shimura} 
Prop.6.12. 

\begin{lemma} \label{2.0} If $L\subseteq F$ containing $k$ is the 
intersection of a collection $\{L_{\alpha}\}_{\alpha}$ of its 
algebraic extensions then the subgroup in $G$ generated by all 
$U_{L_{\alpha}}$ is dense in $U_L$. \end{lemma} 
{\it Proof.} Clearly, $U_L$ contains the subgroup in $G$ generated 
by all $U_{L_{\alpha}}$, and $G_{F/\overline{L}}$ is a normal subgroup 
in each of  $U_{L_{\alpha}}$ and in $U_L$. Set 
$\overline{U}_{L_{\alpha}}=U_{L_{\alpha}}/G_{F/\overline{L}}$, and 
similarly $\overline{U}_L=U_L/G_{F/\overline{L}}$. These are compact 
subgroups in $G_{\overline{L}/k}$. Then, under the above Galois 
correspondence, the closure of the subgroup in $G_{\overline{L}/k}$ 
generated by all $\overline{U}_{L_{\alpha}}$ corresponds to the subfield 
$\overline{L}^{\langle\overline{U}_{L_{\alpha}}~|~\alpha\rangle}=
\bigcap_{\alpha}L_{\alpha}=L$, i.e., to the same subfield as 
$\overline{U}_L$, which means that 
$\langle\overline{U}_{L_{\alpha}}~|~\alpha\rangle$ is dense in 
$\overline{U}_L$, and therefore, 
$\langle U_{L_{\alpha}}~|~\alpha\rangle$ is dense in $U_L$. \qed 

\begin{lemma} \label{norm-int-phi-1} Let $L\subset F$ be such an 
extension of $k$ that any subextension of transcendence degree $\le 2$ 
is of finite type over $k$. Then the common normalizer in $G$ of all 
normal closed subgroups of index $\le 3$ in $U_L$ coincides with $U_L$. 
\end{lemma} 
{\it Proof.} Any element $\tau$ in the common normalizer in $G$ of all 
closed subgroups of index $\le 2$ in $U_L$ satisfies $\tau(L(f^{1/2}))
=L(f^{1/2})$ for all $f\in L^{\times}$. If $\tau\not\in U_L$ then 
there is an element $x\in L^{\times}$ such that $\tau x/x\neq 1$. 
Then $\tau x/x=y^2$ for some $y\in F^{\times}-\{\pm 1\}$. Set 
$f=x+\lambda$ for a variable $\lambda\in k$. By Kummer theory, $\tau 
f/f\in L^{\times 2}$, and therefore, $L$ contains $L_0:=k\left(
y\frac{(x+\lambda y^{-2})^{1/2}}{(x+\lambda)^{1/2}}~|~
\lambda\in k\right)\subset\overline{k(x,y)}$. 

As ${\rm tr.deg}(\overline{k(x,y)}/k)\le 2$, by our assumption 
on $L$, the subfield $L_0$ of $L$ should be finitely 
generated over $k$. But this is possible only if $y^2=1$, 
i.e., if $\tau\in U_L$. (To see this, one can choose 
a smooth model of the extension $L_0(x)/k(x,y)$ 
over $k$ and look at its branch locus.) \qed 

\begin{corollary} \label{tri-cent} For any $\xi\in G-\{1\}$ 
and any open subgroup $U\subseteq G$ there exists an element 
$\sigma\in U$ such that $[\sigma,\xi]\in U-\{1\}$. \end{corollary} 
{\it Proof.} Let $L$ be such a subfield of $F$ finitely 
generated over $k$ that $U_L\subseteq U$. 

Set $L'=L\xi(L)$. Then $\xi^{-1}(\sigma\xi\sigma^{-1})|_L=id$ 
for any $\sigma\in U_{L'}$. By Lemma \ref{norm-int-phi-1}, the 
centralizer in $G$ of $U_{L'}$ is trivial. (Let $\{L_{\alpha}\}$ 
be the set of all finitely generated extensions of $L'$. Then 
$\bigcup_{\alpha}L_{\alpha}=F$. Any element $\tau\in G$ 
centralizing $U_{L'}$ normalizes all subgroups in all 
$U_{L_{\alpha}}$, and thus, by Lemma \ref{norm-int-phi-1}, 
$\tau\in\bigcap_{\alpha}U_{L_{\alpha}}=\{1\}$.) In particular, as 
$\xi$ does not centralize $U_{L'}$, there is $\sigma\in U_{L'}$ 
with $\xi^{-1}\sigma\xi\sigma^{-1}\in U_L-\{1\}$. \qed 

\begin{lemma} \label{uutuun} Let $H\neq\{1\}$ be a normal closed 
subgroup in $G^{\circ}$ such that $H\bigcap U\neq\{1\}$ for a compact 
subgroup $U$. Then $H$ contains $U_{k'(x)}$ for any algebraically 
closed extension $k'$ of $k$ in $F$ such that $k'=
\overline{k'\bigcap F^U}$ and ${\rm tr.deg}(F/k')=1$, and for 
some $x\in F-k'$. \end{lemma} 
{\it Proof.} Let $\sigma\in H\bigcap U-\{1\}$ and $k'$ be the 
algebraic closure in $F$ of any subfield in $F^{\langle\sigma\rangle}$ 
with ${\rm tr.deg}(F/k')=1$. As the extension 
$F/F^{\langle\sigma\rangle}$ is abelian there is an element 
$x\in F-k'$ and an integer $N\ge 2$ such that $\sigma x\neq x$ and 
$\sigma x^N=x^N$. Then one has $\sigma(k'(x))=k'(x)$. 

Let $L$ be a finite Galois extension of $k'(x)$. Its smooth proper 
model over $k'$ is unramified outside a finite set $S$ of points on a 
smooth proper model of $k'(x)$ over $k'$. Then there is an 
element $\overline{\alpha}_L\in{\rm Aut}(k'(x)/k')$ such that 
the set $\overline{\alpha}_L^{-1}\sigma\overline{\alpha}_L(S)$ 
does not intersect $S$, and therefore, for an extension 
$\alpha_L\in G^{\circ}_{F/k'}$ of $\overline{\alpha}_L$ 
to $F$, a smooth proper model over $k'$ of the field 
$L\bigcap\alpha_L^{-1}\sigma\alpha_L(L)$ is unramified over the model 
of $k'(x)$, so $L\bigcap\alpha_L^{-1}\sigma\alpha_L(L)=k'(x)$. 

Let $\beta\in G^{\circ}_{F/k'}$ be given on $k'(x)$ by 
$\overline{\alpha}_L$, and somehow extended to the field 
$\alpha_L^{-1}\sigma\alpha_L(L)$. Then $\beta^{-1}\circ\sigma^{-1}
\circ\beta\circ\alpha_L^{-1}\circ\sigma\circ\alpha_L$ is the identity 
on $k'(x)$, and therefore, induces an automorphism of $L$. Since $L$ 
and $\alpha_L^{-1}\sigma\alpha_L(L)$ are Galois extensions of 
$L\bigcap\alpha_L^{-1}\sigma\alpha_L(L)=k'(x)$, for any given 
automorphism $\tau$ of $L$ over $k'(x)$ there is an extension of 
$\beta$ to $F$ such that for its restriction to $L$ the composition 
$\beta^{-1}\circ\sigma^{-1}\circ\beta\circ\alpha_L^{-1}\circ\sigma
\circ\alpha_L$ coincides with $\tau$. This means that the natural 
projection $H\bigcap U_{k'(x)}\longrightarrow{\rm Gal}(L/k'(x))$ is 
surjective for any  Galois extension $L$ of $k'(x)$, i.e., that 
$H\bigcap U_{k'(x)}$ is dense in $U_{k'(x)}$. As $H\bigcap U_{k'(x)}$ 
is closed, we have $H\supseteq U_{k'(x)}$. \qed 

\begin{lemma} \label{opnor} If $n=1$ then there are no 
proper normal open subgroups in $G^{\circ}$. \end{lemma}
{\it Proof.} Let $H$ be a normal open subgroup in $G^{\circ}$. 
Then for some subfield $L\subset F$ finitely generated over $k$ 
one has $U_L\subseteq H$. For any purely transcendental extension 
$L'\subset F$ of $k$ with $\overline{L'}=F$ one also has 
$U_{LL'}\subseteq H$, as well as $H\supseteq\langle\sigma 
U_{LL'}\sigma^{-1}~|~\sigma\in N_{G^{\circ}}U_{L'}\rangle$. 

A smooth proper model over $k$ of the extension $LL'/L'$ is ramified 
only over a divisor on the model ${\mathbb P}^1_k$ of $L'$ over $k$, 
but the group $N_{G^{\circ}}U_{L'}/U_{L'}\cong{\rm PGL}_2k$ 
does not preserve this divisor, so the intersection 
$\bigcap_{\sigma\in N_{G^{\circ}}U_{L'}}\sigma(LL')$ 
is unramified over $L'$, i.e., the field 
$\bigcap_{\sigma\in N_{G^{\circ}}U_{L'}}\sigma(LL')$ coincides 
with $L'$. By Lemma \ref{2.0}, this shows that $H\supseteq U_{L'}$ 
for any purely transcendental extension $L'\subset F$ of $k$ 
with $\overline{L'}=F$, and therefore, $H$ contains all compact 
subgroups in $G$, so finally, $H\supseteq G^{\circ}$. \qed 

\begin{lemma} \label{reduction} Let $H$ be a closed subgroup of 
$G^{\circ}$. Assume that for any algebraically closed extension $k'$ 
of $k$ in $F$ with ${\rm tr.deg}(F/k')=1$ and any $x\in F-k'$ the 
subgroup $H$ contains $U_{k'(x)}$. Then $H=G^{\circ}$. \end{lemma} 
{\it Proof.} Any element $\sigma$ of a compact subgroup in $G^{\circ}$ 
can be presented as the limit of a compatible collection $(\sigma_L)$ 
of embeddings $\sigma_L$ of finitely generated extensions $L$ of $k$ 
in $F$ into $F$. Replacing $L$ with the compositum of the images of 
$L$ under powers of $\sigma$, we may suppose that $L$ is 
$\sigma$-invariant. 

For a finitely generated $k$-algebra $R$ with the fraction field 
$L^{\langle\sigma\rangle}$, let $X$ be the affine variety over $k$ 
with the coordinate ring $R[t,\sigma t,\sigma^2t,\dots]\subset L$ 
for some $t\in L$ such that $L=L^{\langle\sigma\rangle}(t)$. 
Let $Y:={\bf Spec}(R)\subseteq{\mathbb A}^N_k$. 

By induction on $\dim X$ we show that there is a fibration of a 
Zariski open $\langle\sigma\rangle$-invariant subset of $X$ by 
$\langle\sigma\rangle$-invariant irreducible curves. Namely, if 
$\dim X=1$ there is nothing to prove. If $\dim X>1$ then, by Th{\'e}or{\`e}me 
6.3 4) of \cite{jou},\footnote{For any irreducible scheme $X$ of 
finite type over $k$ and a morphism $X\stackrel{f}{\longrightarrow}
{\mathbb A}_k$ to an affine space with $\dim\overline{f(X)}\ge 2$ the 
preimages of almost all hyperplanes in ${\mathbb A}_k$ are 
irreducible.} pull-backs of sufficiently general hyperplanes in 
${\mathbb A}^N_k$ under the composition $X\longrightarrow Y
\hookrightarrow{\mathbb A}_k$ are irreducible, and thus, there is a 
fibration of a Zariski open $\langle\sigma\rangle$-invariant subset 
$U$ of $X$ by $\langle\sigma\rangle$-invariant irreducible divisors 
with the base $S$. By induction assumption, there is a fibration of a 
Zariski open $\langle\sigma\rangle$-invariant subset of 
$U\times_S\overline{k(S)}$ by $\langle\sigma\rangle$-invariant 
irreducible curves over $\overline{k(S)}$, i.e., there is a fibration 
of a Zariski open $\langle\sigma\rangle$-invariant subset $U'$ of 
$X$ by $\langle\sigma\rangle$-invariant irreducible curves over $k$. 

Then the function field $k''$ of the base is a subfield in 
$L^{\langle\sigma\rangle}$ algebraically closed in $L$ and containing 
$k$. Let $x\in L^{\langle\sigma\rangle}-k''$, and let $k'$ 
be a maximal algebraically closed subfield in 
$F$ containing $k''$ but not $x$. Then $\sigma|_L$ coincides with
restriction to $L$ of an element of $U_{k'(x)}$. 
This shows that $\sigma$ belongs to the closure of the union of 
$U_{k'(x)}$ over all $k'$ and all $x\in F-k'$. \qed 

\begin{lemma} \label{choice} If $\xi$ is a non-trivial element of 
$G$ and $2m\le n$ then there exist elements $w_1,\dots,w_m\in F$ 
such that $w_1,\dots,w_m,\xi w_1,\dots,\xi w_m$ are algebraically 
independent over $k$. \end{lemma} 
{\it Proof.} We proceed by induction on $m$, the case $m=0$ being 
trivial. We wish to find $w_m\in F$ such that $w$ and $\xi w_m$ 
are algebraically independent over $k'$ generated over $k$ by 
$w_1,\dots,w_{m-1},\xi w_1,\dots,\xi w_{m-1}$. Suppose that there 
is no such $w_m$. Then for any $u\in F-\overline{k'}$ and 
any $v\in F-\overline{k'(\xi u)}$ one has the following vanishings 
in $\Omega^2_{F/k'}$: $du\wedge d\xi u=dv\wedge d\xi v=0$, 
$d(u+v)\wedge d\xi(u+v)=0$, and $d(u+v^2)\wedge d\xi(u+v^2)=0$. 
Applying the first two to the third, we get 
$2(v-\xi v)dv\wedge d\xi u=0$, which means that 
$\xi v=v$ for any $v\in F-\overline{k'(\xi u)}$, i.e., $\xi=1$. 
This contradiction shows that there exists desired $w_m\in F$. \qed 

\begin{lemma} \label{inf-cl} Let $L$ be a subfield of $F$ with 
${\rm tr.deg}(F/L)=\infty$. Then $G^{\circ}_{F/L}$ is dense in 
$G_{F/L}$. \end{lemma} 
{\it Proof.} Fix a transcendence basis $x_1,x_2,x_3,\dots$ of $F$ over 
$L$. We wish to show that for each $\sigma\in G_{F/L}$, any integer 
$m\ge 1$ and any $y_1,\dots,y_m\in F$ there is a triplet 
$(\tau_1,\tau_2,\tau_3)$ of elements of some compact subgroups in 
$G_{F/L}$ such that $\tau_3\tau_2\tau_1\sigma y_s=y_s$ for all $1\le s\le 
m$. As $y_1,\dots,y_m$ are algebraic over $k_0:=L(x_1,\dots,x_M)$ 
for some integer $M\ge m$, it is enough to show that for each 
$\sigma\in G_{F/L}$ and any integer $M\ge 1$ there is a pair 
$(\tau_1,\tau_2)$ of elements of some compact subgroups in $G_{F/L}$ 
such that $\tau_2\tau_1\sigma x_s=x_s$ for all $1\le s\le M$. 

Let $k_1:=\overline{k_0\sigma(k_0)}$, $z_j\in F-k_j$ and 
$k_{j+1}:=\overline{k_j(z_j)}$ for all $1\le j\le M$. Let 
$\tau_1\sigma x_j=z_j$ and $\tau_1z_j=\sigma x_j$ for all 
$1\le j\le M$ and $\tau_1$ is somehow extended to an element of a 
compact subgroup in $G_{F/L}$. Let $\tau_2x_j=z_j$ and $\tau_2z_j=x_j$ 
for all $1\le j\le M$ and $\tau_2$ is somehow extended to an element 
of a compact subgroup in $G_{F/L}$. Then $(\tau_2\tau_1)\sigma x_j=x_j$. \qed 

\begin{theorem} \label{exactness} If $n<\infty$ then any non-trivial 
subgroup in $G$ normalized by $G^{\circ}$ is dense in $G^{\circ}$. 
If $n=\infty$ then any non-trivial 
normal subgroup in $G$ is dense. \end{theorem} 
{\it Proof.} Let $H$ be a non-trivial closed subgroup in $G$ 
normalized by $G^{\circ}$. By Corollary \ref{tri-cent}, 
$H$ intersects non-trivially any open compact subgroup in 
$G^{\circ}$ if $n<\infty$. Then, by Lemma \ref{uutuun}, the 
group $H\bigcap G^{\circ}_{F/k'}$ is an open normal subgroup 
in $G^{\circ}_{F/k'}$ for any algebraically closed extension 
$k'$ of $k$ in $F$ with ${\rm tr.deg}(F/k')=1$, and thus, by 
Lemma \ref{opnor}, $H\bigcap G^{\circ}_{F/k'}=G^{\circ}_{F/k'}$. 
Then Lemma \ref{reduction} implies that $H\supseteq G^{\circ}$. 

In the case $n=\infty$ and $H$ is normal, we fix a pair of 
transcendence basis $x_1,x_2,x_3,\dots$ and $y_1,y_2,y_3,\dots$ of $F$ 
over $k$ and show that for each $m$ there is an element $\sigma\in H$ 
such that $\sigma x_j=y_j$ for all $1\le j\le m$. Fix a non-trivial 
element $\xi\in H$. Choose some elements $z_1,\dots,z_m$ algebraically 
independent over $k(x_1,\dots,x_m,y_1,\dots,y_m)$. By Lemma 
\ref{choice}, there exist $w_1,\dots,w_m\in F$ such that 
$w_1,\dots,w_m,\xi w_1,\dots,\xi w_m$ are algebraically independent 
over $k$. Then there exist elements $\tau_1,\tau_2\in G$ such that 
$\tau_1x_j=w_j$, $\tau_1z_j=\xi w_j$ and $\tau_2y_j=w_j$, 
$\tau_2z_j=\xi w_j$ for all $1\le j\le m$. Then 
$\sigma:=\tau^{-1}_2\circ\xi^{-1}\circ\tau_2\circ\tau^{-1}_1
\circ\xi\circ\tau_1$ sends $x_j$ to $y_j$ for all $1\le j\le m$. 

This implies that there is a compact subgroup $U$ intersecting $H$ 
non-trivially, so by Lemma \ref{uutuun}, $H$ contains $U_{k'(x)}$ 
for some algebraically closed extension $k'$ of $k$ in $F$ with 
${\rm tr.deg}(F/k')=1$, and for some $x\in F-k'$. By Lemma 
\ref{opnor}, this implies that $H$ contains $G^{\circ}_{F/k'}$. 

The algebraically closed extensions $k'$ of $k$ in $F$ with 
${\rm tr.deg}(F/k')=1$ form a single $G$-orbit, so one can 
apply Lemma \ref{reduction}, which gives $H\supseteq G^{\circ}$. 
Finally, it follows from Lemma \ref{inf-cl} that $H=G$. \qed 

\begin{corollary} \label{liss-fin-inf} Any smooth 
representation of $G^{\circ}$ of finite degree is trivial. 

({\it Proof.} {\rm This is clear, since there are no proper open 
normal subgroups in} $G^{\circ}$.) \qed \end{corollary} 

\begin{corollary} \label{ex-fin-inf} For any subgroup $H$ of $G$ 
containing $G^{\circ}$ and any continuous homomorphism $\pi$ from 
$H$ either $\pi$ is injective, or the restriction of $\pi$ to 
$G^{\circ}$ is trivial. \end{corollary} 
{\it Proof.} If $\pi$ is not injective then, by Corollary 
\ref{tri-cent}, its kernel has a non-trivial intersection with 
$G^{\circ}$. Then, by Theorem \ref{exactness}, the kernel of 
$\pi$ contains $G^{\circ}$. \qed 

\begin{lemma} \label{lll} Let $d\ge 0$ be an integer, $F_1$ 
and $F_2$ be algebraically closed subfields of $F$ such that 
$F_1\bigcap F_2=k$ and ${\rm tr.deg}(F/F_2)\ge d$. Then for 
any subfield $L$ in $F$ with ${\rm tr.deg}(L/k)=d$ there is 
$\xi\in H:=\langle G_{F/F_1},G_{F/F_2}\rangle$ such that 
${\rm tr.deg}(\xi(L)F_2/F_2)=d$. \end{lemma} 
{\it Proof.} We proceed by induction on $d$, the case $d=0$ being 
trivial. If $d>0$ fix a subfield $L_1\subset L$ with 
${\rm tr.deg}(L_1/k)=d-1$ and some $t\in L$ transcendental over $L_1$. 

Replacing $F_1$ with the algebraic closure in $F$ of the subfield 
generated over $F_1$ by a transcendence basis of $F$ over $F_1F_2$ 
(thus, making $H$ smaller), we may assume that $F$ is algebraic over 
$F_1F_2$. In particular, there exists a subfield $K\subset F_1$ over 
$k$ with ${\rm tr.deg}(K/k)=d$ and ${\rm tr.deg}(KF_2/F_2)=d$. 

By the induction assumption, there exists an element $\tau\in H$ such 
that ${\rm tr.deg}(\tau(L_1)F_2/F_2)=d-1$, i.e., we may suppose that 
${\rm tr.deg}(L_1F_2/F_2)=d-1$. Moreover, we may suppose that 
$L_1=k(t_1,\dots,t_{d-1})$ is purely transcendental over $k$ 
and $L=L_1(t)$. 

The subgroup $G_{F/F_2}$ acts transitively on the set of purely 
transcendental extensions of $F_2$ of a given transcendence degree, 
so for any collection $x_1,\dots,x_{d-1}$ of elements of $F_1$ 
algebraically independent over $F_2$ there is $\sigma\in G_{F/F_2}$ 
such that $\sigma t_j=x_j$. 

If $t\not\in\overline{L_1F_2}$ then induction is completed, 
so we assume that $t$ is algebraic over $L_1F_2$, i.e., there is 
an irreducible polynomial $P\in F_2[X_1,\dots,X_d]-k[X_1,\dots,X_d]$ 
with $P(t_1,\dots,t_d)=0$, where $t_d:=t$. 

Consider the irreducible hypersurface $$W=\{(y_1,\dots,y_d)~|~
P(y_1,\dots,y_d)=0\}\hookrightarrow{\mathbb A}^d_{F_2}\longrightarrow
{\mathbb A}^d_k$$ and the projection $W\stackrel{\pi}{\longrightarrow}
{\mathbb A}^{d-1}_{F_2}$ to the first $d-1$ coordinates. Suppose that 
for any $\sigma$ as above one has $\sigma t\in F_1$. Then for any 
generic point $(x_1,\dots,x_{d-1})\in{\mathbb A}^{d-1}_F$ as above the 
points of the fiber of the projection $\pi$ over $(x_1,\dots,x_{d-1})$ 
are defined over $F_1$. This means that $W$ is defined over $F_1$, and 
therefore, over $F_1\bigcap F_2=k$, contradicting ${\rm tr.deg}(L/k)=d$. 

Therefore, there is $\sigma\in G_{F/F_2}$ such that 
$\sigma(L_1)\subset F_1$ and $\sigma t\not\in F_1$, so in the rest 
of the proof we assume that $L_1\subset F_1$ and $t\not\in F_1$. 

The $G_{F/F_1}$-orbit of $t$ coincides with $F-F_1$. As the 
intersection $$\left(F-F_1\right)\bigcap\left(F-\overline{L_1F_2}
\right)=F-\left(F_1\bigcup\overline{L_1F_2}\right)$$ is non-empty, 
there is an element $\xi$ in $G_{F/F_1}$ such that 
$\xi t\in F-\overline{L_1F_2}$, 
so finally, ${\rm tr.deg}(\xi(L)F_2/F_2)=d$. \qed 

\begin{corollary} \label{ccc} In notations of Lemma \ref{lll}, 
for any $\sigma\in G$ there is $\tau\in H$ such that 
$\sigma|_L=\tau|_L$. \end{corollary} 
{\it Proof.} According to Lemma \ref{lll}, there exist elements 
$\xi,\xi'\in H$ such that $${\rm tr.deg}(\xi'(L)F_2/F_2)=
{\rm tr.deg}(\xi\sigma(L)F_2/F_2)={\rm tr.deg}(L/k).$$ Evidently, 
there is $\lambda\in G_{F/F_2}\subset H$ inducing an 
isomorphism $\xi'(L)\stackrel{\sim}{\longrightarrow}\xi\sigma(L)$ 
such that $\sigma|_L=\xi^{-1}\circ\lambda\circ\xi'|_L$. \qed 

\begin{proposition} \label{2.14} 
Let $L_1$ and $L_2$ be subextensions of $k$ in $F$ such that 
$\overline{L_1}\bigcap\overline{L_2}$ is algebraic over 
$L_1\bigcap L_2$ and ${\rm tr.deg}(F/L_2)=\infty$. 
Then the subgroup in $G$ generated by $G_{F/L_1}$ and $G_{F/L_2}$ 
is dense in $G_{F/L_1\bigcap L_2}$. \end{proposition} 
{\it Proof.} The inclusion $\langle G_{F/L_1},G_{F/L_2}\rangle
\subseteq G_{F/L_1\bigcap L_2}$ is evident. In Corollary \ref{ccc} 
we may replace $k$ with $\overline{L_1}\bigcap\overline{L_2}$ 
to get that the subgroup in $G$ generated by 
$G_{F/\overline{L_1}}$ and $G_{F/\overline{L_2}}$ is dense in 
$G_{F/\overline{L_1}\bigcap\overline{L_2}}$. 

It remains to show that the compact group 
$\overline{\langle G_{F/L_1},G_{F/L_2}\rangle}/
G_{F/\overline{L_1}\bigcap\overline{L_2}}$ coincides with 
${\rm Gal}(\overline{L_1}\bigcap\overline{L_2}/L_1\bigcap L_2)$. But 
this is the Galois theory: $$\left(\overline{L_1}\bigcap\overline{L_2}
\right)^{\langle G_{F/L_1},G_{F/L_2}\rangle}=L_1\bigcap L_2.\qquad\qed$$ 

\begin{lemma} \label{G-0-dec} If $n<\infty$ and 
$F\neq\overline{k'}$, for a subfield $k'$ over $k$ in $F$, 
then $G=G_{F/k'}\cdot G^{\circ}$. \end{lemma} 
{\it Proof.} We may suppose that $k'$ is algebraically closed and 
maximal in $F-\{x\}$ for some $x\in F-k$. Fix a transcendence basis 
$x_1,\dots,x_{n-1}$ of $k'$ over $k$. We have to show that for each 
$\sigma\in G$ there is such an element $\tau\in G^{\circ}$ that 
$\sigma\tau x_s=x_s$ for all $1\le s<n$. 

Set $x_n=x$. Let $0\le j<n$ be the maximal integer with the property 
that there is an element $\tau\in G^{\circ}$ such that 
$\sigma\tau\in G_{F/k''}$, where $k''=\overline{k(x_1,\dots,x_j)}$. 
For such $\tau$ set $\sigma\tau=\sigma'$. Suppose that $0\le j<n-1$. 
Then, by Lemma \ref{choice}, there exists $w\in F-k''$ such that $w$ 
and $\sigma'w$ are algebraically independent over $k''$. There are 
integers $s,s'>j$ such that $w$ is transcendental over $k''(x_s)$ and 
$\sigma'w$ is transcendental over $k''(x_{s'})$. 

Then there are elements $\tau_1,\tau_2,\tau_3,\tau_4$ of some compact 
subgroups of $G_{F/k''}$ such that $\tau_1x_{j+1}=x_s$, $\tau_2x_s=w$, 
$\tau_3(\sigma'w)=x_{s'}$, $\tau_4x_{s'}=x_{j+1}$. Then the 
automorphism $\tau_4\tau_3\sigma'\tau_2\tau_1=\sigma\tau_0$, where 
$\tau_0\in G^{\circ}$, acts trivially on $k''(x_{j+1})$, contradicting 
our assumptions. This means that $j=n-1$. \qed 

\begin{lemma} \label{pur-tra} Let $L$ be a subfield of $F$ and 
$S$ be a set of elements of $F$ algebraically independent over 
$L$. Then the group $H$ generated by $G_{F/L(S-\{x\})}$ 
for all $x\in S$ is dense in $G_{F/L}$. \end{lemma} 
{\it Proof.} If $\overline{H}$ contains the subgroups $G_{F/L(S')}$ 
for all subsets $S'\subset S$ with finite complement then 
it contains the closure of $\bigcup_{S'}G_{F/L(S')}$ 
coincident with $G_{F/L(\bigcap_{S'}S')}=G_{F/L}$. 
So it suffices to treat the case $S=\{x,y\}$. 

For any $z\in F-\overline{L}$ there is an element $\sigma\in 
G_{F/L(x)}G_{F/L(y)}$ such that $\sigma x=z$, so $H$ contains 
$\sigma G_{F/L(x)}\sigma^{-1}=G_{F/L(z)}$, and thus, $H$ is a normal 
subgroup in $U_L$. By Theorem \ref{exactness}, $\overline{H}$ 
contains $G_{F/\overline{L}}$ if ${\rm tr.deg}(F/L)=\infty$, and 
$\overline{H}$ contains $G^{\circ}_{F/\overline{L}}$ otherwise. 
By Lemma \ref{G-0-dec}, in the latter case the projection 
$G_{F/\overline{L}(x)}\longrightarrow G_{F/\overline{L}}/
G^{\circ}_{F/\overline{L}}$ is surjective, and thus, $\overline{
\langle G_{F/\overline{L}(x)},G_{F/\overline{L}(y)}\rangle}$ 
surjects onto $G_{F/\overline{L}}/G^{\circ}_{F/\overline{L}}$, 
so $\overline{\langle G_{F/\overline{L}(x)},
G_{F/\overline{L}(y)}\rangle}=G_{F/\overline{L}}$. 

As the subgroup $G_{F/L(x)}$ of $H$ surjects onto 
$G_{F/L}/G_{F/\overline{L}}\cong{\rm Aut}(\overline{L}/L)$, and $H$ 
is an extension of a subgroup in $G_{F/L}/G_{F/\overline{L}}$ by 
$G_{F/\overline{L}}\subseteq H$, we get $\overline{H}=G_{F/L}$. \qed 

\section{Some geometric representations} \label{Sgr}
Now we are going to construct a supply of semi-simple admissible 
representations of $G$. Recall (\cite{bz}), that for each smooth 
$E$-representation $W$ of $G$ and each compact subgroup $U$ of $G$ 
the Hecke algebra ${\mathcal H}_E(U)$ acts on the space $W^U$, since 
$W^U=h_U(W)$. 

Note that the definition in \S\ref{def-Hecke} on 
p.\pageref{def-Hecke} is equivalent to the usual 
definition of $G$ the Hecke algebra when $U$ is open. 
\begin{proposition} \label{2.10} 
\begin{enumerate} \item \label{irr-loc-glob} 
Let $1\le n\le\infty$. Then a smooth $E$-representation $W$ of $G$ 
is irreducible if and only if for each compact subgroup $U$ with 
$F^U$ purely transcendental over an extension of $k$ of finite 
type the ${\mathcal H}_E(U)$-module $W^U$ is irreducible. 
\item \label{loceq} Let $W_j$ for $j=1,2$ be smooth irreducible 
$E$-representations of $G$ and $W_1^U\neq 0$ for some 
compact subgroup $U$. Then $W_1$ is equivalent to $W_2$ 
if and only if $W_1^U$ is equivalent to $W_2^U$. 
\item \label{existc} For each open compact $U\subset G$ 
and each irreducible $E$-representation $\tau$ 
of the algebra ${\mathcal H}_E(U)$ there is a smooth irreducible 
representation $W$ of $G$ with $\tau\cong W^U$. 
\end{enumerate} \end{proposition}
{\it Proof.} If $n<\infty$ then this is Proposition 2.10 of 
\cite{bz}. In the case $n=\infty$ the proof is modified as follows. 
\begin{enumerate} \item Suppose that $W$ is irreducible and $W^U\neq 0$. 
If $V$ is a non-zero ${\mathcal H}_E(U)$-submodule in $W^U$ then 
the natural morphism ${\bf D}_Eh_U\otimes_{{\mathcal H}_E(U)}V
\longrightarrow W$ is non-zero, and thus, surjective. Then 
$V=h_U{\bf D}_Eh_U\otimes_{{\mathcal H}_E(U)}V$ surjects onto 
$W^U$, i.e., $V=W^U$. 

Conversely, let $W_1\subset W$ be a non-zero proper subrepresentation. 
Fix an element $e\in W-W_1$. There is a compact subgroup $U\subset
{\rm Stab}_e$  with $F^U$ purely transcendental over an extension of 
$k$ of finite type. Then $e\in W^U_1\neq W^U$, giving contradiction. 
\item ``Only if'' part is evident, so suppose that the 
${\mathcal H}_E(U)$-modules $W^U_1$ and $W^U_2$ are isomorphic. 
Then one has ${\bf D}_Eh_U\otimes_{{\mathcal H}_E(U)}W^U_1\cong
{\bf D}_Eh_U\otimes_{{\mathcal H}_E(U)}W^U_2$. Set $W_U:=
\{w\in W~|~h_U{\bf D}_Ew=0\}$ and $_UW:=W/W_U$.
 
Then $_U({\bf D}_Eh_U\otimes_{{\mathcal H}_E(U)}W^U_j)$ is 
irreducible, since otherwise the inclusion 
$0\neq W_3\subset~_U({\bf D}_Eh_U\otimes_{{\mathcal H}_E(U)}W^U_j)$ 
would imply $W^U_3=W^U_j$, and therefore, that 
$W_3=~_U({\bf D}_Eh_U\otimes_{{\mathcal H}_E(U)}W^U_j)$. 
This gives that $$W_1=~_U({\bf D}_Eh_U\otimes_{{\mathcal H}_E(U)}
W^U_1)\cong~_U({\bf D}_Eh_U\otimes_{{\mathcal H}_E(U)}W^U_1)=W_2.\quad 
\qed$$ \end{enumerate} 

\begin{lemma} \label{loc-glo-ss} Let $W$ be a smooth $E$-representation 
of $G$. Suppose that for each compact subgroup $U\subset G$ with $F^U$ 
purely transcendental over an extension of $k$ of finite 
type the ${\mathcal H}_E(U)$-submodule $W^U$ is semi-simple. 
Then $W$ is semi-simple. \end{lemma} 
{\it Proof.} Let $W_0$ be the sum of all irreducible $E$-subrepresentations 
of $G$ in $W$. Suppose $W\neq W_0$ and $e\in W-W_0$. Let $V$ be the 
subrepresentation of $G$ in $W$ generated by $e$, and $V_0=V\bigcap W_0$. 
Let $V_2$ be a maximal subrepresentation of $G$ in $V-\{e\}$ containing 
$V_0$. Then the representation $V_1=V/V_2$ of $G$ is irreducible. Our goal 
is to embed $V_1$ into $W$. Let $U$ be a compact subgroup in $G$ such that 
$F^U$ is purely transcendental over an extension of $k$ of finite type and 
$V_1^U\neq 0$. By Proposition \ref{2.10}(\ref{irr-loc-glob}), the 
${\mathcal H}_E(U)$-module $V_1^U$ is irreducible. As the 
${\mathcal H}_E(U)$-module $V^U$ is semi-simple and $V^U/V^U_2$ is 
isomorphic to $V_1^U$, there is an ${\mathcal H}_E(U)$-submodule $N_U$ 
in $V^U$ isomorphic to $V_1^U$ and intersecting $V_2$ trivially. 

For each open subgroup $U'$ of $U$ the ${\mathcal H}_E(U')$-submodule 
$N'_{U'}$ of $U'$-invariants of the subrepresentation of $G$ in $V$ 
generated by $N_U$ splits into a direct sum $\bigoplus_jA_j^{m_j}$ of 
irreducible ${\mathcal H}_E(U')$-modules, where $m_j\le\infty$ and 
$A_i\not\cong A_j$ if $i\neq j$. 

As $(N'_{U'})^U=h_U\ast{\bf D}_EN_U=h_U\ast{\bf D}_E\ast h_UN_U=N_U$, 
there is exactly one index $j$ such that $A_j^U\neq 0$. For this index 
one has $m_j=1$. By  Proposition \ref{2.10}(\ref{existc}), there are 
smooth irreducible representations $W_1$ and $W_2$ of $G$ such that 
$W_1^{U'}\cong A_j$ and $W_2^{U'}\cong V_1^{U'}$. 
Since $W_1^U\cong W_2^U$, by Proposition \ref{2.10}(\ref{loceq}), 
one has $W_1\cong W_2$, so $A_j\cong V_1^{U'}$. Set $N_{U'}=A_j$.

Then $N:=\bigcup_{U'\subset U}N_{U'}$ is a $E$-subrepresentation of $G$ 
in $V$ intersecting $V_0$ trivially. (The action of $\sigma\in G$ on 
$N_{U_L}$ for some finite extension $L$ of $F^U$ is given by the 
composition of the embedding $N_{U_L}\subseteq N_{U_{L\sigma(L)}}$ 
with $h_{L\sigma(L)}\ast\sigma\ast h_{L\sigma(L)}\in
{\rm End}(N_{U_{L\sigma(L)}})$. Here $h_L:=h_{U_L}$ is the Haar 
measure on $U_L$, and we have used $\sigma:W^{U_L}\longrightarrow 
W^{U_{\sigma(L)}}\subseteq W^{U_{L\sigma(L)}}$.) By Proposition 
\ref{2.10}(\ref{irr-loc-glob}), the $E$-representation $N$ of $G$ 
is irreducible. This contradicts the definition of $W_0$. \qed 

\vspace{4mm}

\label{Hecke-zykly}
For any irreducible variety $Y$ over $k$ with $k(Y)=F^U$ for a 
compact open subgroup $U$ in $G$ one can identify the Hecke algebra 
${\mathcal H}(U)$ with the ${\mathbb Q}$-algebra of 
non-degenerate correspondences on $Y$ (i.e., of formal linear combinations 
of $n$-subvarieties in $Y\times_kY$ dominant over both factors $Y$). 
This follows from the facts that the set of double classes 
$U\backslash G/U$ can be identified with a basis of ${\mathcal H}(U)$ as 
a ${\mathbb Q}$-space via $[\sigma]\longmapsto h_U\ast\sigma\ast h_U$; 
that irreducible $n$-subvarieties in $Y\times_kY$ dominant over 
both factors $Y$ are in a natural bijection with the set of 
maximal ideals in $F^U\otimes_kF^U$, and the following lemma. 
\begin{lemma} \label{max-spectrum} Let $L,L'\subseteq F$ be field 
subextension of $k$ with ${\rm tr.deg}(L/k)=q<\infty$. Then the set of 
double classes $U_{L'}\backslash G/U_L$ is canonically identified with 
the set of all points in ${\bf Spec}(L\otimes_kL')$ of codimension 
$\ge q-{\rm tr.deg}(F/L')$ (so $U_{L'}\backslash G/U_L
={\bf Max}(L\otimes_kL')$, if $F=\overline{L'}$). Here 
$G/U_L=\{\mbox{embeddings of $L$ into $F$ over $k$}\}$. \end{lemma} 
{\it Proof.} To any embedding $\sigma:L\hookrightarrow F$ over 
$k$ one associates the ideal in $L\otimes_kF$ generated by 
elements $x\otimes 1-1\otimes\sigma x$ for all $x\in F$. It is 
maximal, since it is the kernel of the surjection $L\otimes_kF
\stackrel{\sigma\cdot id}{\longrightarrow}F$. Conversely, any 
maximal ideal ${\mathfrak m}$ in $L\otimes_kF$ determines a homomorphism 
$L\otimes_kF\stackrel{\varphi}{\longrightarrow}(L\otimes_kF)/
{\mathfrak m}=\Xi$ with $\Xi$ a field. Since its restriction to the 
subfield $k\otimes_kF$ is an embedding, one can regard $F$ as 
a subfield of $\Xi$. Let $t_1,\dots,t_m$ be a transcendence basis of 
$L$ over $k$. As $L$ is algebraic over $L_0:=k(t_1,\dots,t_m)$, $\Xi$ 
is integral over $\varphi(L_0\otimes_kF)$, so the latter is a field. 
One has ${\bf Spec}(L_0\otimes_kF)\subset{\bf Spec}(F[t_1,\dots,t_m])=
{\mathbb A}^m_F$. On any subvariety of ${\mathbb A}^m_F$ outside the 
union of all divisors defined over $k$ there is an $F$-point which 
also lies outside the union of all divisors defined over $k$, 
and therefore, ${\bf Max}(L_0\otimes_kF)\subset
{\bf Max}(F[t_1,\dots,t_m])$. This means that 
$\varphi(L_0\otimes_kF)=F$, and therefore, $\Xi=F$. 
The restriction of $\varphi$ to $L\otimes_kk$ gives 
an embedding $\sigma:L\hookrightarrow F$ over $k$. \qed 

\vspace{4mm}

Let $U$ and $V$ be open compact subgroups of $G$, and $F^U=L$, $F^V=L'$. 
Then, in notation of \S\ref{NCT}, p.\pageref{def-Hecke}, one has 
the isomorphism $h_V{\bf D}_{{\mathbb Q}}h_U\stackrel{\sim}{\longrightarrow}
{\rm Hom}_G({\bf D}_{{\mathbb Q}}h_V,{\bf D}_{{\mathbb Q}}h_U)$ 
given by right multiplication 
(the inverse sends a $G$-homomorphism to its value on $h_V$), 
which is evidently compatible with multiplication in 
${\bf D}_{{\mathbb Q}}$ and composing $G$-homomorphisms 
$${\rm Hom}_G({\bf D}_{{\mathbb Q}}h_W,{\bf D}_{{\mathbb Q}}h_V)
\times{\rm Hom}_G({\bf D}_{{\mathbb Q}}h_V,{\bf D}_{{\mathbb Q}}h_U)
\longrightarrow
{\rm Hom}_G({\bf D}_{{\mathbb Q}}h_W,{\bf D}_{{\mathbb Q}}h_U).$$ 

The canonical projection ${\bf D}_{{\mathbb Q}}\longrightarrow
{\mathbb Q}[G/V]$ identifies the $G$-module ${\bf D}_{{\mathbb Q}}h_V$ 
with ${\mathbb Q}[\{L'\stackrel{/k}{\hookrightarrow}F\}]$. 
The latter can be also regarded as the $G$-module of generic 0-cycles 
${\mathbb Q}[{\bf Max}(L'\otimes_kF)]$. Similarly, ${\bf D}_{{\mathbb Q}}h_U
={\mathbb Q}[\{L\stackrel{/k}{\hookrightarrow}F\}]=
{\mathbb Q}[{\bf Max}(L\otimes_kF)]$. 

The correspondences ${\mathbb Q}[{\bf Max}(L\otimes_kL')]={\mathbb Q}
[V\backslash\{L\stackrel{/k}{\hookrightarrow}F\}]$ act on the space 
${\mathbb Q}[{\bf Max}(L'\otimes_kF)]$ as follows. A correspondence 
$[\tau]\in V\backslash\{L\stackrel{/k}{\hookrightarrow}F\}$ sends 
a cycle $L'\stackrel{\sigma}{\hookrightarrow}F$ to
$\sum\limits_{\xi\in V/U_{\tau(L)L'}}[\sigma\xi\tau]\in
{\mathbb Q}[\{L\stackrel{/k}{\hookrightarrow}F\}]$. 

This gives an isomorphism $${\mathbb Q}[{\bf Max}(L\otimes_kL')]
\stackrel{\sim}{\longrightarrow}{\rm Hom}_G({\mathbb Q}
[\{L'\stackrel{/k}{\hookrightarrow}F\}],{\mathbb Q}
[\{L\stackrel{/k}{\hookrightarrow}F\}])$$ compatible with composing 
of correspondences and of $G$-homomorphisms. Its inverse is given 
by the composition \begin{multline*}{\rm Hom}_G({\mathbb Q}
[\{L'\stackrel{/k}{\hookrightarrow}F\}],{\mathbb Q}
[\{L\stackrel{/k}{\hookrightarrow}F\}])
\stackrel{(id_{L'})}{\longrightarrow}
{\mathbb Q}[\{L\stackrel{/k}{\hookrightarrow}F\}] \\
\stackrel{\alpha}{\longrightarrow}{\mathbb Q}[V\backslash
\{L\stackrel{/k}{\hookrightarrow}F\}]=
{\mathbb Q}[{\bf Max}(L\otimes_kL')],\end{multline*}
where $(id_{L'})$ is the value on the element 
$L'\stackrel{id}{\hookrightarrow}F$ and 
$\alpha([\tau])=\frac{1}{[L'\tau(L):L']}[\tau]$. 

\vspace{4mm}

\label{shrink}
Let $A^q(Y)$ be the quotient of the ${\mathbb Q}$-space $Z^q(Y)$ of 
cycles on a smooth proper variety $Y$ over $k$ of codimension $q$ by 
the ${\mathbb Q}$-subspace $Z^q_{\sim}(Y)$ of cycles $\sim$-equivalent 
to zero for an adequate equivalence relation $\sim$. According to 
Hironaka, each smooth variety $X$ admits an open embedding $i$ 
into a smooth proper variety $\overline{X}$ over $k$. 
Then $A^q(-)$ can be extended to arbitrary smooth variety $X$ as 
the cokernel of the map $Z^q_{\sim}(\overline{X})
\stackrel{i^{\ast}}{\longrightarrow}Z^q(X)$ induced by restriction 
of cycles. This is independent of the choice of variety 
$\overline{X}$.\footnote{since for any pair of smooth 
compactifications $(\overline{X},\overline{X}')$ of $X$ there is their 
common refinement $\overline{X}\stackrel{\beta}{\longleftarrow}
\overline{X}''\stackrel{\beta'}{\longrightarrow}\overline{X}'$, 
$i^{\ast}$ factors through $Z^q_{\sim}(\overline{X})
\stackrel{\beta^{\ast}}{\longrightarrow}Z^q_{\sim}(\overline{X}'')
\stackrel{(i'')^{\ast}}{-\!\!\!-\!\!\!\longrightarrow}Z^q(X)$ and 
$i^{\ast}Z^q_{\sim}(\overline{X})
=(i'')^{\ast}Z^q_{\sim}(\overline{X}'')$.} 

In the standard way one extends the contravariant functors $A^q(~)$ 
and $Z^q(~)$ to contravariant functors on the category of smooth 
pro-varieties over $k$. Namely, if for a set of indices $I$, an 
inverse system $(X_j)_{j\in I}$ of smooth varieties over $k$ is 
formed with respect to flat morphisms and $X$ is the limit, then 
$Z^q(X)=\lim\limits_{_{j\in I}\longrightarrow}Z^q(X_j)$, where the 
direct system is formed with respect to the pull-backs, and similarly 
for $A^q(~)$. This is independent of the choice of the projective 
system defining $X$. 

In particular, as for any commutative $k$-algebra $R$ the scheme 
${\bf Spec}(R)$ is an inverse limit of a system of $k$-varieties, 
$A^q(R):=A^q({\bf Spec}(R))$ is defined. Any automorphism 
$\alpha$ of the $k$-algebra $R$ induces a morphism of a system 
$(X_j)_{j\in I}$ defining ${\bf Spec}(R)$ to a system 
$(\alpha^{\ast}(X_j))_{j\in I}$ canonically equivalent to 
$(X_j)_{j\in I}$, and therefore, induces an automorphism of 
$A^q(Y_R)$ for any $k$-scheme $Y$. This gives a contravariant 
functor from a category of varieties over $k$ to the category 
of ${\rm Aut}(R/k)$-modules. 
Set $B^q(X)=A^q(X)$ for $\sim=$numerical equivalence. 

In what follows $X$ will be of type $Y_F$ for a $k$-subscheme 
$Y$ in a variety over $k$. It should be stressed that in this 
case $B^q(X)=B^q(Y_F)$ means {\sl not} some sort of numerical 
equivalence over $F$, but a limit of cycle groups modulo 
numerical equivalence of varieties over $k$. 

\label{shrink-e}

\vspace{4mm}

The homomorphism of algebras ${\mathcal H}(U)\longrightarrow 
A^{\dim Y}(Y\times_kY)$ is surjective for any smooth projective 
$Y$, as one can see from the following ``moving lemma'', 
applied in the case $X_1=X_2=Y$ and $Z=X_1\times_kX_2$. 
(Its present form is suggested by the referee.)

\begin{lemma} \label{moving} Let $Z,X_1,\dots,X_r$ be irreducible 
projective varieties over $k$, and let $Z\stackrel{p_j}{\longrightarrow}X_j$ 
be surjective maps. Let $\alpha\subset Z$ be an irreducible subvariety 
of dimension $q\ge\max\limits_{1\le j\le r}\dim X_j$. Then $\alpha$ 
is rationally equivalent to a linear combination of some irreducible 
subvarieties in $Z$ surjective (under maps $p_j$) over all $X_j$'s. 
\end{lemma}
{\it Proof.} Choose a closed irreducible subvariety $W$ in $Z$ 
containing $\alpha$ as a divisor such that $W$ is surjective (under 
maps $p_j$) over all $X_j$. We can replace $Z$ by $W$ and assume 
that $\alpha$ has codimension 1 in $Z$. Then we can replace $Z$ by 
the blowup of $Z$ along $\alpha$. Since $\alpha$ is the direct 
image of its pullback, we can assume $\alpha$ is a Cartier divisor on 
$Z$. Since our varieties are taken to be projective, any Cartier 
divisor is rationally equivalent to a difference of two very ample 
divisors, which we take to dominate the $X_j$'s. \qed 

\begin{corollary} \label{cor-mov} The natural map 
$Z^q(k(X)\otimes_kk(Y))\longrightarrow A^q(X\times_kY)$ is surjective 
if $q\le\dim X\le\dim Y$; and $Z^q(k(Y)\otimes_kF)\longrightarrow 
A^q(Y_F)$ is surjective if $q\le\dim Y\le n$. \qed \end{corollary}

\begin{proposition} \label{geomod} 
Let $Y$ be a smooth irreducible proper variety over $k$ and 
$\dim Y\le n$. Let $X$ be a smooth variety 
over $k$. For each $q\ge 0$ there are canonical isomorphisms 
\begin{multline*}A^q(X_{k(Y)})\stackrel{\sim}{\longrightarrow}{\rm Hom}_G
(A^{\dim Y}(Y_F),A^q(X_F))\\ \stackrel{\sim}{\longrightarrow}
{\rm Hom}_G(Z^{\dim Y}(k(Y)\otimes_kF),A^q(X_F)).\end{multline*}

For each pair of reduced irreducible group schemes ${\mathcal A}$ 
and ${\mathcal B}$ over $k$ there is a natural bijection 
$${\rm Hom}({\mathcal A},{\mathcal B}):={\rm Hom}_{{\rm group~schemes}/k}
({\mathcal A},{\mathcal B})\stackrel{\sim}{\longrightarrow}{\rm Hom}_G
({\mathcal A}(F),{\mathcal B}(F)).$$ \end{proposition}
{\it Proof} of the first part uses Corollary \ref{cor-mov}, Lemma 
\ref{max-spectrum} and elementary intersection theory as follows. 
Each embedding $\sigma:k(Y)\stackrel{/k}{\hookrightarrow}F$ 
induces an identification $Z^{\dim Y}(k(Y)\otimes_kF)
\stackrel{\sim}{\longrightarrow}{\mathbb Q}[G/U_{\sigma(k(Y))}]$, 
and thus, for each $G$-module $M$ one has an isomorphism 
$M^{U_{\sigma(k(Y))}}\stackrel{\sim}{\longrightarrow}
{\rm Hom}_G(Z^{\dim Y}(k(Y)\otimes_kF),M)$ given by 
$m\longmapsto[\tau\sigma\longmapsto\tau m]$ for any $\tau\in G$. 

Let $Z^{\dim Y}(k(Y)\otimes_kF)\stackrel{\varphi}{\longrightarrow}
A^q(X_F)$ be a $G$-homomorphism. Fix an embedding 
$k(Y)\stackrel{\sigma}{\hookrightarrow}F$. For any 
$\xi\in G$ one has $\varphi(\xi\sigma)=\xi\varphi(\sigma)$, 
in particular, if $\xi|_{\sigma(k(Y))}=id$ then $\varphi(\sigma)
\in A^q(X_F)^{U_{\sigma(k(Y))}}=A^q(X_{\sigma(k(Y))})$.\footnote{For a 
subextension $L\subseteq F$ of $k$ let $L'$ be a purely transcendental
subextension of $L$ over which $F$ is algebraic. By Galois descent 
property, $A^q(X_F)^{U_{L'}}=A^q(X_{L'})$. By the homotopy invariance, 
$A^q(X_{L'})=A^q(X_L)$, so $A^q(X_L)\subseteq A^q(X_F)^{U_L}\subseteq 
A^q(X_F)^{U_{L'}}=A^q(X_{L'})=A^q(X_L)$.} From this and the fact that 
the pairing $A^q(X\times_kY)\otimes A^{\dim Y}(Y_F)\longrightarrow 
A^q(X_F)$ factors through $A^q(X_{k(Y)})\otimes A^{\dim Y}(Y_F)
\longrightarrow A^q(X_F)$,\footnote{since any correspondence 
supported on $X\times_kD$ for a divisor $D$ on $Y$ sends 
$A^{\dim Y}(Y_F)$ to zero} one deduces that sending a cycle 
$\alpha\in A^q(X_{k(Y)})$ to the action of its arbitrary lifting 
to a correspondence on $X\times_kY$ determines a homomorphism 
$$A^q(X_{k(Y)})\longrightarrow{\rm Hom}_G(A^{\dim Y}(Y_F),A^q(X_F)),$$ 
which is surjective and canonical. On the other hand, it has the 
inverse given by $\varphi\longmapsto
\sigma^{-1}_{\ast}\varphi(\sigma)$, where $\sigma_{\ast}:
A^q(X_{k(Y)})\stackrel{\sim}{\longrightarrow}A^q(X_{\sigma(k(Y))})$. 

The map ${\rm Hom}({\mathcal A},{\mathcal B})\longrightarrow
{\rm Hom}_G({\mathcal A}(F),{\mathcal B}(F))$ is clearly 
injective. Let ${\mathcal A}(F)\stackrel{\varphi}{\longrightarrow}
{\mathcal B}(F))$ be a $G$-homomorphism. Fix an irreducible curve 
$C\subset{\mathcal A}$ over $k$ generating ${\mathcal A}$ as an 
algebraic group, i.e., with the dominant multiplication map 
$C^N\longrightarrow{\mathcal A}$ for any $N\ge\dim{\mathcal A}$. 
Fix a generic point $x\in C(F)$. Since $\varphi(\tau x)=\tau\varphi(x)$, 
the element $\varphi(x)$ is fixed by any element $\tau$ fixing $x$, 
so any coordinate of $\varphi(x)$ is a rational function over $k$ 
in coordinates of $x$, and therefore, this gives rise to a rational 
$k$-map $h:C--\!\rightarrow{\mathcal B}$. Consider the rational $k$-map 
$C^N--\!\rightarrow{\mathcal B}$ given by $(x_1,\dots,x_N)\longmapsto 
h(x_1)\cdots h(x_N)$. On the set of $F$-points out of the union
of ``vertical'' divisors, (i.e., on 
${\bf Max}(\underbrace{k(C)\otimes_k\cdots\otimes_kk(C)}
\limits_{N~{\rm copies}}\otimes_kF)$) this map coincides with one 
given by $(x_1,\dots,x_N)\longmapsto\varphi(x_1)\cdots\varphi(x_N)=
\varphi(x_1\cdots x_N)$. 

As the multiplication map ${\bf Max}
(\underbrace{k(C)\otimes_k\cdots\otimes_kk(C)}\limits_{N~{\rm copies}}
\otimes_kF)\longrightarrow{\mathcal A}(F)$ is surjective for any integer 
$N\ge 2\dim{\mathcal A}$, the rational map $C^N--\!\rightarrow{\mathcal B}$ 
factors through $C^N--\!\rightarrow{\mathcal A}
\stackrel{\widehat{h}}{\longrightarrow}{\mathcal B}$. 
Since $\varphi$ is a homomorphism, $\widehat{h}$ should 
also be a homomorphism, and in particular, regular. \qed 

\begin{corollary} \label{Hom-G-mod} For any field $L'$ of finite 
type and of transcendence degree $m\le n$ over $k$, any field 
$L$ of finite type over $k$ and any integer $q\ge m$ 
there is a canonical isomorphism 
\begin{equation} \label{prim-mor} \phantom{\qquad\qed}
A^q(L\otimes_kL')\stackrel{\sim}{\longrightarrow}
{\rm Hom}_G(A^m(L'\otimes_kF),A^q(L\otimes_kF)),\end{equation}
where the both groups are zero if $q>m$. \end{corollary}
{\it Proof.} Let $Y$ be an irreducible smooth 
projective variety over $k$ with the function field $k(Y)=L'$. 
The group $A^m(L'\otimes_kF)$ is the quotient of $A^m(Y_F)$ by 
the sum of the images of $A^{m-1}(\widetilde{D}_F)$ for all 
divisors $D$ on $Y$ and their desingularizations $\widetilde{D}$. 

Then the target of (\ref{prim-mor}) is 
a subgroup in ${\rm Hom}_G(A^m(Y_F),A^q(L\otimes_kF))$. 

By Proposition \ref{geomod}, 
${\rm Hom}_G(A^m(Y_F),A^q(L\otimes_kF))=A^q(L\otimes_kL')$, and 
$${\rm Hom}_G(A^{m-1}(\widetilde{D}_F),A^q(L\otimes_kF))=
A^q(L\otimes_kk(D)),$$ which is zero, since $\dim D<q$. 
The vanishing of 
${\rm Hom}_G(A^{m-1}(\widetilde{D}_F),A^q(L\otimes_kF))$ 
for all divisors $D$ on $Y$ implies the coincidence 
of both sides in (\ref{prim-mor}). \qed

\begin{proposition} \label{supply} The $G$-module\footnote{Recall 
that $B^q(X\times_kF)$ is a limit of certain quotients of 
${\mathbb Q}$-spaces of classes of numerical equivalence of cycles 
on smooth proper varieties over $k$, but not over $F$, cf. 
pp.\pageref{shrink}--\pageref{shrink-e} before Lemma \ref{moving}.} 
${\bf B}^q_X=B^q(X\times_kF)$ is admissible for any 
smooth proper $k$-variety $X$ and any $q\ge 0$. 
If $n<\infty$ then ${\bf B}^q_X$ is semi-simple. 

If $q=0$, or $q=1$, or $q=\dim X\le n$ then ${\bf B}^q_X$ 
is semi-simple and of finite length. \end{proposition} 
{\it Proof.} By a standard argument, we may suppose that $k$ is 
embedded into the field of complex numbers ${\mathbb C}$, and thus, 
for any smooth proper $k$-variety $Y$ with $k(Y)\subset F$ the 
space $B^q(X\times_kF)^{U_{k(Y)}}$ is a quotient of the 
finite-dimensional space $Z^q(X\times_kY)/\sim_{{\rm hom}}
\subseteq H^{2q}((X\times_kY)({\mathbb C}),{\mathbb Q}(q))$, 
so the representation $B^q(X\times_kF)$ is admissible. 

For each smooth irreducible variety $Y$ over $k$ with 
$k(Y)\subset F$ and $\overline{k(Y)}=F$ the kernel 
of $A^q(X\times_kY)\longrightarrow A^q(X\times_kk(Y))$ is a 
$A^n(Y\times_kY)$-submodule in $A^q(X\times_kY)$, since for its 
arbitrary element $\alpha$ and for any element 
$\beta\in A^n(Y\times_kY)$ one has $\alpha\circ\beta={\rm pr}_{13\ast}
({\rm pr}_{12}^{\ast}\alpha\cdot{\rm pr}_{23}^{\ast}\beta)$, so the 
projection to $Y$ of the support of $\alpha\circ\beta$ is contained 
in ${\rm pr}_2((D\times_kY)\bigcap{\rm supp}(\beta))$ for a divisor 
$D$ on $Y$, which is of dimension $n-1$, so cannot dominate $Y$. 
This implies that $A^q(X\times_kk(Y))$ has a natural structure 
of a $A^n(Y\times_kY)$-module. 

By \cite{jan}, the algebra $B^n(Y\times_kY)$ is semi-simple, 
so the $B^n(Y\times_kY)$-module $B^q(X\times_kk(Y))$ is also 
semi-simple. By the moving lemma \ref{moving}, the ring homomorphism 
${\mathcal H}_{U_{k(Y)}}\longrightarrow B^n(Y\times_kY)$, induced by 
the identification of the Hecke algebra ${\mathcal H}_{U_{k(Y)}}$ 
with the algebra of non-degenerate correspondences on $Y$, 
(see p.\pageref{Hecke-zykly}) is surjective. This gives a 
(semi-simple) ${\mathcal H}_{U_{k(Y)}}$-module structure on any 
$B^n(Y\times_kY)$-module. Then, by Lemma \ref{loc-glo-ss}, 
the $G$-module $B^q(X\times_kF)$ is semi-simple. 

Now suppose that $n=\infty$. By the same result of Jannsen \cite{jan}, 
the category of motives modulo numerical equivalence is semi-simple. 
Let $(X,\Delta_X)=\bigoplus_j(X,\pi_j)$ be a decomposition into 
a (finite) direct sum of irreducible submotives. Then $B^q(X_F)=
\bigoplus_j\pi_jB^q(X_F)$. If $W=\pi_jB^q(X_F)$ is reducible 
there is a non-zero proper $G$-submodule $W_0$ in $W$. Fix elements 
$e_0\in W_0-{0}$ and $e_1\in W-W_0$. Then the common stabilizer 
of $e_0$ and $e_1$ is an open subgroup in $G$, so it contains a 
subgroup $U_L$ for a subfield $L$ of $F$ finitely generated over $k$. 
Let $F'$ be an algebraically closed extension of $L$ with 
$\dim X\le{\rm tr.deg}(F'/k)<\infty$. Then $W_0^{G_{F/F'}}$ 
is a non-zero proper $G_{F'/k}$-module in $W^{G_{F/F'}}=\pi_j\left(
B^q(X_F)\right)^{G_{F/F'}}$, so the length of the $G_{F'/k}$-module 
$\pi_j\left(B^q(X_F)\right)^{G_{F/F'}}=\pi_jB^q(X_{F'})$ is $\ge 2$. 
By Proposition \ref{geomod}, there is a canonical surjection 
${\rm End}_{{\rm motive}/k}\left((X,\pi_j)\right)\longrightarrow
{\rm End}_{G_{F'/k}}(\pi_jB^q(X_{F'}))$ if $q\in\{0,1,\dim X\}$. 
By the irreducibility of $(X,\pi_j)$ we have a division algebra 
on the left hand side, but the algebra on the right hand side 
has divisors of zero since the $G_{F'/k}$-module $\pi_jB^q(X_{F'})$ 
is semi-simple, but not irreducible, giving contradiction. 

Any cyclic semi-simple $G$-module, ${\bf B}^{\dim X}_X$ in 
particular (Corollary \ref{cor-mov}), is of finite length. 

It follows from Lefschetz theorem on $(1,1)$-classes that for 
any smooth proper $k$-variety $Y$ one has $B^1(X\times_kY)=B^1(X)
\oplus{\rm Hom}({\rm Alb}X,{\rm Alb}Y)\oplus B^1(Y)$. By Lefschetz 
hyperplane section theorem, inclusion $C\hookrightarrow X$ of any 
smooth 1-dimensional plane section $C$ of $X$ induces a surjection 
${\rm Alb}C\longrightarrow{\rm Alb}X$. This implies that 
${\bf B}^1_X$ is embedded into $B^1(X)
\oplus{\bf B}^1_C$, so it is also of finite length. \qed 

\begin{corollary} \label{H-G-m} One has ${\rm Hom}_G(
B^q(L'\otimes_kF),B^p(L\otimes_kF))=0$ for any pair of 
fields $L,L'$ finitely generated over $k$ with ${\rm tr.deg}(L/k)=p$, 
${\rm tr.deg}(L'/k)=q$ and $p\neq q$. \end{corollary}
{\it Proof.} If $p>n$, or $q>n$, then at least one of $B^q
(L'\otimes_kF)$ and $B^p(L\otimes_kF)$ is zero, so we may assume 
that $\max(p,q)\le n$. By Proposition \ref{supply}, the $G$-modules 
$B^q(L'\otimes_kF)$ and $B^p(L\otimes_kF)$ are semi-simple, so 
${\rm Hom}_G(B^q(L'\otimes_kF),B^p(L\otimes_kF))$ is isomorphic to 
${\rm Hom}_G(B^p(L\otimes_kF),B^q(L'\otimes_kF))$, so we may assume 
that $p>q$. Then, by Corollary \ref{Hom-G-mod}, one has ${\rm Hom}_G
(B^q(L'\otimes_kF),B^p(L\otimes_kF))=B^p(L\otimes_kL')=0$. \qed 

\begin{corollary} \label{contr-0-cycl} Let $X$ and $Y$ be smooth 
irreducible proper varieties over $k$. Then the ${\mathbb Q}$-vector spaces 
$B^{\dim X}(X_{k(Y)})$ and $B^{\dim Y}(Y_{k(X)})$ are naturally dual. 
If $n<\infty$ and $\dim X\le n$ this duality induces a non-degenerate 
$G$-equivariant pairing 
$$B^{\dim X}(X_F)\otimes\lim\limits_{_U\longrightarrow}
B^n((Y_U)_{k(X)})\longrightarrow{\mathbb Q}(\chi),$$ where $U$ runs 
over the set of open compact subgroups in $G$, $Y_U$ is a smooth 
proper model of $F^U$ over $k$ (thus, $\dim Y_U=n$) and the direct 
system is formed with respect to the pull-backs on the cycles. 
\end{corollary}
{\it Proof.} Let $\dim Y\ge\dim X$. 
Set $n=\dim Y$. By Proposition \ref{geomod}, 
$$B^{\dim X}(X_{k(Y)})={\rm Hom}_G(B^n(Y_F),B^{\dim X}(X_F))$$ and 
$B^n(Y_{k(X)})={\rm Hom}_G(B^{\dim X}(X_F),B^n(Y_F))$. 
By Proposition \ref{supply}, the $G$-modules $B^n(Y_F)$ and 
$B^{\dim X}(X_F)$ are semi-simple and of finite length. For any 
$\alpha\in B^{\dim X}(X_{k(Y)})$ and $\beta\in B^n(Y_{k(X)})$ set 
$\langle\alpha\cdot\beta\rangle={\rm tr}(\alpha\circ\beta)$ 
($={\rm tr}(\beta\circ\alpha)$). Here $\alpha$ and $\beta$ are 
considered as $G$-homomorphisms. If $\alpha\neq 0$ there is an 
element $\gamma\in B^n(Y_{k(X)})$ such that $\alpha\circ\gamma$ 
is a non-zero projector in ${\rm End}_GB^{\dim X}(X_F)$, so the 
form $\langle\phantom{-}\cdot\phantom{-}\rangle$ is non-degenerate. 

The form $B^{\dim X}(X_F)\otimes\lim\limits_{_U\longrightarrow}
B^n((Y_U)_{k(X)})\longrightarrow{\mathbb Q}(\chi)$ is defined by 
$\alpha\otimes\beta\longmapsto\langle\alpha\cdot\beta\rangle
\cdot[U]$, for any $\alpha\in B^{\dim X}(X_F)^U$ and $\beta\in 
B^n((Y_U)_{k(X)})$. By the projection formula, it is well-defined. \qed 

\vspace{4mm}

The above examples of $G$-modules are obtained from some 
(pro-)varieties over $k$ by extending the base field to $F$. 
More generally, one can construct a $G$-module starting from 
some birationally invariant functor ${\mathcal F}$ on a category 
of varieties over $k$, or on a category of field extensions 
of $k$ (as in Corollary \ref{contr-0-cycl}). 

Starting with the functor ${\rm Div}_{{\rm alg}}$ of algebraically 
trivial divisors on the category of smooth proper varieties over 
$k$,\footnote{${\rm Div}_{{\rm alg}}(Y_U)$ is independent of the 
choice of $Y_U$, since any birational morphism induces a homomorphism 
of the groups of algebraically trivial divisors, which is an 
isomorphism of the subgroups of linearly trivial divisors (= 
multiplicative groups of the function fields modulo $k^{\times}$) and 
induces an isomorphism of the quotients modulo linear equivalence (= 
${\rm Pic}^{\circ}$-groups).} or with the functor 
${\rm Pic}^{\circ}_{{\mathbb Q}}$, we get another examples of 
$G$-modules of this type:  ${\rm Div}^{\circ}_{{\mathbb Q}}=
\lim\limits_{_U\longrightarrow}{\rm Div}_{{\rm alg}}(Y_U)_{{\mathbb Q}}$, 
and ${\rm Pic}^{\circ}_{{\mathbb Q}}=\lim\limits_{_U\longrightarrow}
{\rm Pic}^{\circ}(Y_U)_{{\mathbb Q}}$, where $U$ runs over the set of 
open subgroups of type $U_L$ and $Y_U$ is a smooth projective model 
of $F^U=L$ over $k$. 

If ${\mathcal A}$ is a commutative group scheme over $k$, 
we set $W_{{\mathcal A}}={\mathcal A}(F)/{\mathcal A}(k)$. 
\begin{proposition} \label{pic-knot} ${\rm Pic}^{\circ}_{{\mathbb Q}}=
\bigoplus_{{\mathcal A}}{\mathcal A}(k)\otimes_{{\rm End}({\mathcal A})}
W_{{\mathcal A}^{\vee}}$, where ${\mathcal A}$ runs over the isogeny classes 
of simple abelian varieties over $k$, and 
${\mathcal A}^{\vee}:={\rm Pic}^{\circ}{\mathcal A}$ is the dual abelian 
variety. \end{proposition}
{\it Proof.} For any open compact subgroup $U$ there is 
a canonical decomposition $$\bigoplus_{{\mathcal A}}{\mathcal A}(k)
\otimes_{{\rm End}({\mathcal A})}{\rm Hom}({\mathcal A},
{\rm Pic}^{\circ}Y_U)_{{\mathbb Q}}\stackrel{\sim}{\longrightarrow}
{\rm Pic}^{\circ}(Y_U)_{{\mathbb Q}}=
\left({\rm Pic}^{\circ}_{{\mathbb Q}}\right)^U$$ given by 
$a\otimes\varphi\longmapsto\varphi(a)$ for any $a\in{\mathcal A}(k)$ and 
any $\varphi\in{\rm Hom}({\mathcal A},{\rm Pic}^{\circ}Y_U)_{{\mathbb Q}}$, 
where ${\mathcal A}$ runs over the isogeny classes of simple abelian 
varieties over $k$. (Clearly, the image of $a\otimes t\varphi$, i.e., 
$t\varphi(a):=\varphi(ta)$ coincides with the image of 
$ta\otimes\varphi$ for any $t\in{\rm End}({\mathcal A})$, so the map 
is well-defined.) Passing to the direct limit with respect to $U$, 
we get ${\rm Pic}^{\circ}_{{\mathbb Q}}=\bigoplus_{{\mathcal A}}
\left({\rm Pic}^{\circ}_{{\mathbb Q}}\right)_{{\mathcal A}}$, where 
$\left({\rm Pic}^{\circ}_{{\mathbb Q}}\right)_{{\mathcal A}}:={\mathcal A}(k)
\otimes_{{\rm End}({\mathcal A})}\lim\limits_{_U\longrightarrow}
{\rm Hom}({\mathcal A},{\rm Pic}^{\circ}Y_U)_{{\mathbb Q}}$ is the 
${\mathcal A}$-isotypic component. 

Using the identifications \begin{multline*}
{\rm Hom}({\mathcal A},{\rm Pic}^{\circ}
Y_U)_{{\mathbb Q}}={\rm Hom}({\rm Alb}Y_U,{\mathcal A}^{\vee}
)_{{\mathbb Q}}\\ =\left({\rm Mor}(Y_U,{\mathcal A}^{\vee})/
{\mathcal A}^{\vee}(k)\right)_{{\mathbb Q}}=\left(
W_{{\mathcal A}^{\vee}}\right)^U,\end{multline*} we get 
$\lim\limits_{_U\longrightarrow}{\rm Hom}({\mathcal A},
{\rm Pic}^{\circ}Y_U)_{{\mathbb Q}}=W_{{\mathcal A}^{\vee}}$, so 
$$\left({\rm Pic}^{\circ}_{{\mathbb Q}}\right)_{{\mathcal A}}={\mathcal A}(k)
\otimes_{{\rm End}({\mathcal A})}W_{{\mathcal A}^{\vee}}.\qquad\qed$$ 

\begin{corollary} \label{semspl} For any smooth irreducible variety 
$X$ of dimension $\le n+1$ over $k$ and $q\in\{0,1,2,\dim X\}$ 
there is a unique $G$-submodule in ${\bf B}^q_X=B^q(X\times_kF)$ 
isomorphic to $B^q(k(X)\otimes_kF)$. \end{corollary}
{\it Proof.} For any $q\ge 0$ and any adequate 
relation $\sim$ one has the short exact sequence 
$$\bigoplus_{D\in X^1}A^{q-1}(\widetilde{D}_F)
\longrightarrow A^q(X_F)\longrightarrow A^q(k(X)\otimes_kF)
\longrightarrow 0.$$ Let $F'$ be an algebraically closed extension 
of $k$ in $F$ with $\dim X-1\le{\rm tr.deg}(F'/k)<\infty$. As 
$\left({\bf B}^q_X\right)^{G_{F/F'}}$ is semi-simple, there is a 
subrepresentation of $G_{F'/k}$ in $\left({\bf B}^q_X\right)^{G_{F/F'}}$ 
isomorphic to $B^q(k(X)\otimes_kF')$. 

By Proposition \ref{geomod}, there is an embedding 
$${\rm Hom}_{G_{F'/k}}(A^{q-1}(\widetilde{D}_{F'}),
A^q(k(X)\otimes_kF'))\hookrightarrow A^q(k(X)\otimes_kk(D))=0$$ 
for $q=\dim X$; ${\rm Hom}_{G_{F'/k}}(A^{q-1}(\widetilde{D}_{F'}),
A^q(k(X)\otimes_kF'))=0$ for $q\in\{0,1\}$; and 
$A^1(\widetilde{D}_{F'})$ is isomorphic to a subquotient of 
$A^1(C_{F'})\oplus{\mathbb Q}^N$ for a smooth proper curve $C$, 
so $${\rm Hom}_{G_{F'/k}}(A^1(\widetilde{D}_{F'}),A^2(k(X)\otimes_kF'))
\subseteq A^2(k(X)\otimes_kk(C))=0.$$ This implies that for any $q$ 
in the range of the statement, any $G_{F'/k}$-equivariant homomorphism 
$$\left({\bf B}^q_X\right)^{G_{F/F'}}\longrightarrow 
B^q(k(X)\otimes_kF')$$ 
factors through an endomorphism of $B^q(k(X)\otimes_kF')$, and 
therefore, by semi-simplicity, that there is a unique 
$G_{F'/k}$-submodule in $\left({\bf B}^q_X\right)^{G_{F/F'}}$ 
isomorphic to $B^q(k(X)\otimes_kF')$. Then the union over all 
$F'$ of such $G_{F'/k}$-submodules is the unique $G$-submodule 
in ${\bf B}^q_X$ isomorphic to $B^q(k(X)\otimes_kF)$. \qed 

\vspace{4mm}

For each open compact subgroup 
$U\subset G$ and a smooth irreducible variety $Y$ over $k$ with 
$k(Y)=F^U$ we define a semi-simple $G$-module (of finite length) 
${\bf B}^q_{Z,Y}$ as the minimal one such that the 
${\mathcal H}(U)$-module $({\bf B}^q_{Z,Y})^U$ is isomorphic 
to $B^q(Z\times_kY)$. By Proposition \ref{2.10}, it exists and it is unique. 
\begin{lemma} \label{redvan} Let $X$, $Y$ and $Z$ be smooth irreducible 
$k$-varieties, $\dim X=\dim Y=n\ge\dim Z$, and $p,q\ge 0$ integers. 
Let $${\mathcal B}^{p,q}_{X,Y,Z}={\rm Hom}_{{\mathcal H}(U)}(B^q(Z\times_kY),
B^p(k(X)\otimes_kk(Y))).$$ Then ${\mathcal B}^{p,q}_{X,Y,Z}=0$, 
if either $q=\dim Z<p$, or $q=n$ and $\dim Z<p$, or $q>n$ and 
$p+q>\dim Z+n$, or $q<p$ and $q\in\{0,1\}$. \end{lemma} 
{\it Proof.} \begin{itemize} \item Let $\dim Z=q<p$. As 
$B^q(Z\times_kY)=({\bf B}^q_{Z,Y})^U$, it follows from the moving 
lemma \ref{moving} that $W_1:={\bf B}^q_{Z,Y}$ is a quotient of the 
module $Z^q(k(Z)\otimes_kF)$ (since $W_1^U=W^U$, where $W$ 
is the quotient of $Z^q(k(Z)\otimes_kF)$ by its $G$-submodule 
generated by the kernel of $Z^q(k(Z)\otimes_kk(Y))
\longrightarrow B^q(Z\times_kY)$). 

As ${\rm Hom}_G(Z^q(k(Z)\otimes_kF),B^p(k(X)\otimes_kF))=0$, 
this gives $${\rm Hom}_G(W_1,B^p(k(X)\otimes_kF))=0.$$ 

We need to show that ${\rm Hom}_{{\mathcal H}(U)}
(W^U,B^p(k(X)\otimes_kk(Y)))=0$. 

By Proposition \ref{2.10}, for any pair $(W_1,W_2)$ of semi-simple 
$G$-modules the natural homomorphism ${\rm Hom}_G(W_1,W_2)
\longrightarrow{\rm Hom}_{{\mathcal H}(U)}(W_1^U,W_2^U)$ is 
surjective. Then, as its source is zero when $W_2=B^p(k(X)\otimes_kF)$, 
we get the vanishing of the space 
${\rm Hom}_{{\mathcal H}(U)}(({\bf B}^q_{Z,Y})^U,
B^p(k(X)\otimes_kk(Y)))$. 
\item Let $c=n-\dim Z$, so the variety $Z\times{\mathbb P}^c$ is 
$n$-dimensional. As there is an embedding of 
${\mathcal H}(U)$-modules $B^n(Z\times_kY)\hookrightarrow 
B^n(Z\times_k
{\mathbb P}^c\times_kY)$, it is enough to show the vanishing of 
${\mathcal B}^{p,n}_{X,Y,Z\times{\mathbb P}^c}$ for $p>\dim Z$. 
By semi-simplicity, the latter is a subgroup in ${\rm Hom}_G
({\bf B}^n_{Z\times{\mathbb P}^c,Y},B^p(k(X)\otimes_kF))\subseteq 
{\rm Hom}_G(Z^n(k(Z\times{\mathbb P}^c)\otimes_kF),
B^p(k(X)\otimes_kF))$. By Proposition \ref{geomod}, the latter 
coincides with $B^p(k(X)\otimes_kk(Z\times{\mathbb P}^c))$, which 
is dominated by $B^p({\mathbb A}^c_{k(X)\otimes_kk(Z)})
\stackrel{\sim}{\longleftarrow}B^p(k(X)\otimes_kk(Z))=0$, 
and thus, ${\mathcal B}^{p,n}_{X,Y,Z}=0$, if $p>\dim Z$. 
\item As the ${\mathcal H}(U)$-module $\bigoplus_{D\in Z^{q-n}}
B^n(\widetilde{D}\times_kY)$ surjects onto the 
${\mathcal H}(U)$-module $B^q(Z\times_kY)$ when $q>n$, the 
vanishing of ${\mathcal B}^{p,q}_{X,Y,Z}$ follows from 
${\mathcal B}^{p,n}_{X,Y,\widetilde{D}}=0$ 
for each subvariety $D$ of codimension $q-n$ on $Z$. 
\item If $q=1$ then $B^1(Z\times_kY)$ is a subquotient of 
$B^1(C\times_kY)$ for a smooth curve $C$, so we are reduced to the 
case $q=\dim Z<p$. The case $q=0$ is trivial. \qed \end{itemize} 

\begin{proposition} \label{split} Let $X$ and $Y$ be smooth 
irreducible varieties over $k$, and either $q\in\{0,1,2\}$, 
or $q=\dim X=\dim Y$. Then there is a unique submodule 
in $B^q(X\times_kY)$ over $\left(B^{\dim X}
(X\times_kX)\otimes B^{\dim Y}(Y\times_kY)^{{\rm op}}\right)$ 
isomorphic to its quotient $B^q(k(X)\otimes_kk(Y))$. \end{proposition} 
{\it Proof.} The existence of such submodule follows from the 
semi-simplicity  of $B^q(X\times_kY)$ (\cite{jan}). By Lemma \ref{redvan}, 
$${\rm Hom}_{B^{\dim Y}(Y\times_kY)}(B^{q-1}(\widetilde{D}\times_kY),
B^q(k(X)\otimes_kk(Y)))=0.$$ As the kernel of the projection 
$B^q(X\times_kY)\longrightarrow B^q(k(X)\otimes_kk(Y))$ is generated 
by the images of $B^{q-1}(X\times_k\widetilde{E})$ and 
$B^{q-1}(\widetilde{D}\times_kY)$ for all divisors $D$ on $X$ and all 
divisors $E$ on $Y$, this implies the uniqueness. \qed 

\subsection{The projector $\Delta_{k(X)}$}
For any pair of varieties $X,Y$ let $~^t\phantom{\alpha}$ 
be the transposition of cycles, induced by $X\times Y
\stackrel{\sim}{\longrightarrow}Y\times X$. Denote by 
$\Delta_{k(X)}=~^t\Delta_{k(X)}$ the identity (diagonal) element 
in $B^n(k(X)\otimes_kk(X))$ considered as an element of 
$B^n(X\times_kX)$. 
\begin{lemma} \label{proj} For any irreducible smooth proper 
$k$-variety $X$ of dimension $n$ the element $\Delta_{k(X)}$ is 
a central projector in the algebra $B^n(X\times_kX)$. The left 
(equivalently, right) ideal generated by $\Delta_{k(X)}$ coincides 
with (the image of) $B^n(k(X)\otimes_kk(X))$. \end{lemma} 
{\it Proof.} Denote by $\varphi$ the projection $B^n(X\times_kX)
\longrightarrow B^n(k(X)\otimes_kk(X))$ and by $\psi$ its unique
section $B^n(k(X)\otimes_kk(X))\longrightarrow B^n(X\times_kX)$. 
The kernel of $\varphi$ coincides with the sum of the kernels of 
$B^n(X\times_kX)\longrightarrow B^n(X\times_kk(X))$ and of $B^n(X\times_kX)
\longrightarrow B^n(k(X)\times_kX)$, where the projections 
are induced by the ring homomorphisms $$B^n(X\times_kX)
\longrightarrow{\rm End}_GB^n(X\times_kF)=B^n(X\times_kk(X))$$ 
and $B^n(X\times_kX)\longrightarrow{\rm End}_GB^n(F\times_kX)
=B^n(k(X)\times_kX)$, so $\ker\varphi$ is a two-sided ideal. 
Then the image of $\varphi$ is a $B^n(X\times_kX)$-bi-module, 
and thus, the image of $\psi$ is a two-sided ideal in 
$B^n(X\times_kX)$. 

As $\alpha\circ\Delta_X=\Delta_X\circ\alpha$ for any $\alpha\in 
B^n(X\times_kX)$, and $\varphi$ and $\psi$ are morphisms of 
$B^n(X\times_kX)$-bi-modules, one has $\alpha\psi(\varphi(
\Delta_X))=\psi(\varphi(\alpha\circ\Delta_X))=\psi(\varphi(
\Delta_X\circ\alpha))=\psi(\varphi(\Delta_X))\alpha$, so $\alpha
\Delta_{k(X)}=\Delta_{k(X)}\alpha$. The $B^n(X\times_kX)$-action on 
$B^n(k(X)\otimes_kk(X))$ factors through $B^n(k(X)\otimes_kk(X))$, 
so $\Delta_{k(X)}^2=\Delta_X\Delta_{k(X)}=\Delta_{k(X)}$. \qed 

\begin{lemma} \label{max-proj} $\Delta_{k(X)}B^{\dim X}(X_L)=
B^{\dim X}(k(X)\otimes_kL)$ for any irreducible smooth proper 
$k$-variety $X$ and any field extension $L$ of $k$. \end{lemma} 
{\it Proof.} Set $d=\dim X$. Using homotopy invariance and the Galois 
descent property of $B^{\ast}$, we may replace $L$ by an algebraically 
closed extension $F$ with ${\rm tr.deg}(F/k)=n\ge d$ and then 
$\Delta_{k(X)}B^d(X_L)=(\Delta_{k(X)}B^d(X_F))^{G_{F/L}}$ and 
$B^d(k(X)\otimes_kL)=B^d(k(X)\otimes_kF)^{G_{F/L}}$. 

By Proposition \ref{geomod} and Corollary \ref{Hom-G-mod}, 
the canonical maps \begin{multline*}B^d(k(X)\otimes_kk(X))
\longrightarrow{\rm End}_GB^d(k(X)\otimes_kF)\\ \longrightarrow
{\rm Hom}_G(B^d(X_F),B^d(k(X)\otimes_kF))\end{multline*} 
are isomorphisms. As $\Delta_{k(X)}$ is the identity element 
in ${\rm End}_GB^d(k(X)\otimes_kF)$, this means that 
$\Delta_{k(X)}B^d(X_F)=B^d(k(X)\otimes_kF)$. \qed 

\begin{proposition} \label{max-prim} $(X,\Delta_{k(X)})$ is the 
maximal primitive $n$-submotive of the motive $(X,\Delta_X)$ for 
any smooth irreducible proper $n$-dimensional $k$-variety $X$. 
The motive $(X,\Delta_{k(X)})$ is a birational invariant of $X$. 
\end{proposition} 
{\it Proof.} To show that $(X,\Delta_{k(X)})$ 
is a primitive $n$-motive, we have to check that the ${\mathbb Q}$-vector 
space $W:=\Delta_{k(X)}B^n(X\times_kY\times{\mathbb P}^1)$ is zero 
for any variety $Y$ of dimension $<n$. Replacing $Y$ with 
$Y\times{\mathbb P}^{n-\dim Y-1}$ we may suppose that $\dim Y=n-1$. 

$W$ is a left $B^n(k(X)\otimes_kk(X))$-module, since by Lemma 
\ref{proj}, $\Delta_{k(X)}$ is a central projector in 
$B^n(X\times_kX)$ and $B^n(X\times_kX)\Delta_{k(X)}=
B^n(k(X)\otimes_kk(X))$. 

In notations of Lemma \ref{redvan} one has \begin{multline*}{\rm Hom}
_{B^n(X\times_kX)}(B^n(X\times_kY\times{\mathbb P}^1),
B^n(k(X)\otimes_kk(X)))=:{\mathcal B}^{n,n}_{X,X,Y\times{\mathbb P}^1}\\ 
={\mathcal B}^{n,n}_{X,X,Y}\bigoplus{\mathcal B}^{n,n-1}_{X,X,Y}. 
\end{multline*} 
By the first case of Lemma \ref{redvan}, one has 
${\mathcal B}^{n,n-1}_{X,X,Y}=0$; by the second case of Lemma 
\ref{redvan}, one has ${\mathcal B}^{n,n}_{X,X,Y}=0$, so 
${\mathcal B}^{n,n}_{X,X,Y\times{\mathbb P}^1}=0$. Then the semi-simplicity 
implies that $${\rm Hom}_{B^n(X\times_kX)}(B^n(k(X)\otimes_kk(X)),
B^n(X\times_kY\times{\mathbb P}^1))=0.$$ As $W$ is a quotient of 
a direct sum of several copies of $B^n(k(X)\otimes_kk(X))$, 
but there are no non-zero $B^n(X\times_kX)$-module quotients of 
$B^n(k(X)\otimes_kk(X))$ in $B^n(X\times_kY\times{\mathbb P}^1)$, 
this means that $W=0$. 

For the maximality of $(X,\Delta_{k(X)})$ among primitive 
$n$-submotives of the motive $(X,\Delta_X)$, we have to show that 
${\rm Hom}(M,X)={\rm Hom}(M,(X,\Delta_{k(X)}))$ for any 
primitive $n$-motive $M=(Z,\pi)$. ${\rm Hom}(X,M)=\pi 
B^n(Z\times_kX)$. For any divisor $D$ on $X$ one has 
$\pi B^{n-1}(Z\times_k\widetilde{D})=0$, so from the exact sequence 
$$\bigoplus_{D\in X^1}\pi B^{n-1}(Z\times_k\widetilde{D})
\longrightarrow\pi B^n(Z\times_kX)\longrightarrow
\pi B^n(Z\times_kk(X))\longrightarrow 0$$ 
we get $\pi B^n(Z\times_kX)=\pi B^n(Z\times_kk(X))$. By Proposition 
\ref{geomod}, $$B^n(Z\times_kk(X))={\rm Hom}_G
(B^n(X_F),B^n(Z_F)),$$ so ${\rm Hom}(X,M)=\pi{\rm Hom}_G
(B^n(X_F),B^n(Z_F))$. By Proposition \ref{geomod}, for any divisor 
$D$ on $X$ one has $${\rm Hom}_G(B^{n-1}(D_F),\pi B^n(Z_F))=
\pi B^{n-1}(Z\times_kk(D))=0,$$ so 
${\rm Hom}(X,M)=\pi{\rm Hom}_G(B^n(k(X)\otimes_kF),B^n(Z_F))$. 
By Lemma \ref{max-proj}, the space $B^n(k(X)\otimes_kF)$ coincides 
with the image of the projector $\Delta_{k(X)}$ acting on $B^n(X_F)$, 
so \begin{multline*}{\rm Hom}(X,M)=
\pi{\rm Hom}_G(\Delta_{k(X)}B^n(X_F),B^n(Z_F))\\ 
=\pi{\rm Hom}_G(B^n(X_F),B^n(Z_F))\Delta_{k(X)}.\end{multline*} 
Since the space ${\rm Hom}_G(B^n(X_F),B^n(Z_F))$ is a quotient of 
$B^n(Z\times_kX)$, we get that the space ${\rm Hom}(X,M)$ is a 
quotient of $\pi B^n(Z\times_kX)\Delta_{k(X)}={\rm Hom}
((X,\Delta_{k(X)}),M)$. On the other hand, the space ${\rm Hom}
((X,\Delta_{k(X)}),M)$ is a quotient of ${\rm Hom}(X,M)$. 
As both spaces are finite-dimensional, they coincide. 

The birational invariantness of $(X,\Delta_{k(X)})$ follows from the 
birational invariantness of $Y^{{\rm prim}}$ for any smooth projective 
variety $Y$ over $k$ explained in the beginning of \S\ref{functors-B} 
below.\footnote{One can show this directly as follows. By Corollary 
\ref{descr-proj}, for any primitive $n$-motive $M\cong(Y,\pi)$ with 
$\dim Y=n$ one has ${\rm Hom}(X,M):=\pi B^n(Y\times_kX)=\pi B^n(k(Y)
\times_kk(X))$, which is independent of the model of $k(X)$.} \qed 

\begin{corollary} \label{descr-proj} Let $X$ and $Y$ be smooth 
irreducible proper varieties over $k$, and $\dim X=\dim Y=n$. 
Then the unique submodule of the module 
$B^n(X\times_kY)$ over $\left(B^{\dim X}(X\times_kX)\otimes 
B^{\dim Y}(Y\times_kY)^{{\rm op}}\right)$ isomorphic to its quotient 
$B^n(k(X)\otimes_kk(Y))$ coincides with $\Delta_{k(X)}\cdot 
B^n(X\times_kY)=B^n(X\times_kY)\cdot\Delta_{k(Y)}=
\Delta_{k(X)}\cdot B^n(X\times_kY)\cdot\Delta_{k(Y)}$. \end{corollary} 
{\it Proof.} By Proposition \ref{max-prim}, $(X,\Delta_{k(X)})$ and 
$(Y,\Delta_{k(Y)})$ are the maximal primitive $n$-submotives in 
$(X,\Delta_X)$ and $(Y,\Delta_Y)$, so $${\rm Hom}(X,(Y,\Delta_{k(Y)}))
={\rm Hom}((X,\Delta_{k(X)}),(Y,\Delta_{k(Y)})).$$ By definition, 
\begin{multline*} {\rm Hom}(X,(Y,\Delta_{k(Y)})):=B^n(X\times_kY)\cdot
\Delta_{k(Y)}\\ \mbox{and}\quad{\rm Hom}((X,\Delta_{k(X)}),
(Y,\Delta_{k(Y)})):=\Delta_{k(X)}\cdot B^n(X\times_kY)\cdot
\Delta_{k(Y)},\end{multline*} and thus, 
$B^n(X\times_kY)\cdot\Delta_{k(Y)}=\Delta_{k(X)}\cdot 
B^n(X\times_kY)\cdot\Delta_{k(Y)}$. 

Similarly, $\Delta_{k(X)}\cdot B^n(X\times_kY)=
\Delta_{k(X)}\cdot B^n(X\times_kY)\cdot\Delta_{k(Y)}$. 

As in the proof of maximality of $(X,\Delta_{k(X)})$ in Proposition 
\ref{max-prim} we have $\Delta_{k(X)}\cdot B^n(X\times_kY)=
\Delta_{k(X)}\cdot B^n(X\times_kk(Y))$. By Lemma \ref{max-proj} 
the latter coincides with $B^n(k(X)\otimes_kk(Y))$. \qed 

\subsection{The functors ${\mathbb B}^{\bullet}$ and ${\mathfrak B}^q$} 
\label{functors-B} 
For a smooth projective variety $Y$ over $k$ let $Y^{{\rm prim}}$ 
be the motive defined by $Y^{{\rm prim}}:=\!\!\!\!\bigcap
\limits_{Y\stackrel{\varphi}{\longrightarrow}M\otimes{\mathbb L}}\!\!\!\!
\ker\varphi$, where $M$ runs over isomorphism classes of {\sl effective} 
motives, or equivalently, $$Y^{{\rm prim}}:={\rm coker}\left[
\bigoplus\limits_{M\otimes{\mathbb L}\stackrel{\varphi}{\longrightarrow}Y}
M\otimes{\mathbb L}\longrightarrow Y\right].$$ Clearly, 
$Y\longmapsto Y^{{\rm prim}}$ is a functor from the category of 
smooth projective varieties to the category of motives modulo 
numerical equivalence. Any birational map is a composition of 
blow-ups and blow-downs with smooth centers (\cite{akmw,wlo}). As 
a blow-up does not change $Y^{{\rm prim}}$ (cf. \cite{manin}), 
this implies that $Y^{{\rm prim}}$ is an invariant of the 
function field $k(Y)$. According to Hironaka, for any subfield 
$L$ of $F$ finitely generated over $k$ there exists a smooth 
projective variety $Y_{[L]}$ over $k$ with the function field 
$L$, and therefore, one gets a canonical projective system of 
motives $\{Y^{{\rm prim}}_{[L]}\}_L$ indexed by subfields $L$ 
of $F$ finitely generated over $k$. 

Now we define the functor ${\mathbb B}^{\bullet}=
\bigoplus{\mathbb B}^{[i]}$ of Theorem \ref{summ} from 
the category of motives modulo numerical equivalence to 
the category of graded ${\mathbb Q}$-spaces by setting 
${\mathbb B}^{[i]}=\lim\limits_{_L\longrightarrow}{\rm Hom}
\left(Y^{{\rm prim}}_{[L]}\otimes{\mathbb L}^{\otimes i},-\right)$ 
for its component of degree $i$. Let also ${\mathfrak B}^q$ 
denotes the restriction of ${\mathbb B}^{[0]}$ to 
the subcategory of the primitive $q$-motives. $G$ acts on the 
projective system $\{Y^{{\rm prim}}_{[L]}\}_L$ by $Y_{[L]}
\stackrel{\sigma}{\longrightarrow}Y_{[\sigma(L)]}$, 
$\sigma(L)\stackrel{\sigma^{-1}}{\longrightarrow}L$, so $G$ acts 
on the limits ${\mathfrak B}^q(M)$ and ${\mathbb B}^{\bullet}(M)$. 

\vspace{4mm}

\label{proof-of-prim-dec} 
{\sc Remark.} Any Grothendieck motive modulo numerical 
equivalence $M=(X,\pi)$ is isomorphic to $\bigoplus\limits
_{0\le i,j,i+j\le\dim X}M_{ij}\otimes{\mathbb L}^{\otimes i}$, 
where $M_{ij}$ is a primitive $j$-motive and 
${\mathbb L}=({\mathbb P}^1,{\mathbb P}^1\times\{0\})$, so 
${\mathbb B}^{[i]}(M)\cong\bigoplus_j{\mathfrak B}^j(M_{ij})$. One proves 
this by induction on dimension $d$ of $X$ as follows. Let $M_{0d}=
\bigcap_{\varphi}\ker(\varphi)$, where $\varphi$ runs over the set 
of morphisms from $M$ to motives of type $(Y\times{\mathbb P}^1,\Delta)$ 
for all $Y$ with $\dim Y<d$. (By Proposition \ref{max-prim}, 
$M_{0d}=(X,\pi\circ\Delta_{k(X)})$.) As the length of $M$ is 
$\le\dim_{{\mathbb Q}}{\rm End}(M)<\infty$, the motive $M/M_{0d}$ 
can be embedded into a finite direct sum of 
$(Y_j\times{\mathbb P}^1,\Delta)$ with $\dim Y_j<d$. As 
$(Y_j\times{\mathbb P}^1,\Delta)=(Y_j,\Delta)\oplus(Y_j,\Delta)
\otimes{\mathbb L}$, the induction is completed. In fact, the 
decomposition $M=\bigoplus\limits_{0\le i,j,i+j\le\dim X}
\widetilde{M}_{ij}$, where $\widetilde{M}_{ij}$ is isomorphic 
to $M_{ij}\otimes{\mathbb L}^{\otimes i}$, is canonical since 
$$\widetilde{M}_{ij}={\rm Im}\left[\bigoplus_{N:{\rm primitive}
~j\mbox{-}{\rm motives},N\otimes{\mathbb L}^{\otimes i}
\stackrel{\varphi}{\rightarrow}M}N\otimes{\mathbb L}^{\otimes i}
\stackrel{\sum\varphi}{\longrightarrow}M\right].$$ 

\begin{proposition} \label{b-explic} If $\dim X=q\le n$ 
and $M=(X,\pi)$ is a primitive $q$-motive then 
${\mathfrak B}^q(M)=\pi B^q(X_F)$. \end{proposition}
{\it Proof.} First we wish to show that ${\mathbb B}^{[0]}
(M)^{U_L}={\rm Hom}\left(Y_{[L]},M\right)$ for any $L\subset F$ 
of finite type over $k$, so (by Lemma \ref{sm-exh} below) 
${\mathbb B}^{[0]}(M)^{G_{F/K}}=\lim\limits_{_{K\supseteq L}\rightarrow}
{\rm Hom}\left(Y_{[L]},M\right)$. Any element $\alpha\in
{\mathbb B}^{[0]}(M)^{U_L}$ belongs to the image of ${\rm Hom}
\left(Y_{[L'']},M\right)$ for some $L''\supseteq L$. We may 
assume that $L''$ is a Galois extension of a purely transcendental 
extension $L'$ of $L$. Fix a finite affine open covering 
$\{U_{\gamma}\}$ of $Y_{[L']}$. 
Let $B_{\gamma}$ be the integral closure of ${\mathcal O}(U_{\gamma})$ 
in $L''$, $V_{\gamma}={\bf Spec}(B_{\gamma})$ and $Y''=
{\coprod\limits_{\gamma}}\phantom{}_{Y_{[L']}}
V_{\gamma}$. Then ${\rm Gal}(L''/L')$ acts on each 
$V_{\gamma}$, and therefore, ${\rm Gal}(L''/L')$ acts on $Y''$ 
with a smooth quotient $Y_{[L']}$, so $Y''$ is projective. By 
equivariant version of resolution of singularities, there is a smooth 
projective variety $Y_{[L'']}$ with a ${\rm Gal}(L''/L')$-action and 
a ${\rm Gal}(L''/L')$-equivariant birational morphism 
$Y_{[L'']}\longrightarrow Y''$, so $\alpha$ belongs to the image of 
${\rm Hom}\left(Y_{[L'']},M\right)^{{\rm Gal}(L''/L')}=\pi B^q\left(
X\times_kY_{[L'']}\right)^{{\rm Gal}(L''/L')}$. Let $\widetilde{Y}$ 
be a smooth projective variety admitting birational morphisms to 
$Y_{[L']}$ and to $Y_{[L'']}/{\rm Gal}(L''/L')$. Then 
\begin{multline*}{\rm Hom}\left(Y_{[L']}^{{\rm prim}},M\right)=
{\rm Hom}\left(\widetilde{Y}^{{\rm prim}},M\right)=
\pi B^q\left(X\times_k\widetilde{Y}\right)\\ \longrightarrow
\!\!\!\!\rightarrow\pi B^q\left(X\times_kY_{[L'']}/
{\rm Gal}(L''/L')\right)=\pi B^q\left(X\times_kY_{[L'']}
\right)^{{\rm Gal}(L''/L')}.\end{multline*} On the other hand, by the 
projection formula, $${\rm Hom}\left(Y_{[L']}^{{\rm prim}},M
\right)=\pi B^q\left(X\times_kY_{[L']}\right)\hookrightarrow
\pi B^q\left(X\times_kY_{[L'']}\right)^{{\rm Gal}(L''/L')},$$ so 
$$\pi B^q\left(X\times_kY_{[L'']}\right)^{{\rm Gal}(L''/L')}=
{\rm Hom}\left(Y_{[L']}^{{\rm prim}},M\right)={\rm Hom}\left(
Y_{[L]}^{{\rm prim}},M\right),$$ which means that $\alpha$ 
belongs to the image of ${\rm Hom}\left(Y_{[L]},M\right)$, and 
thus, ${\mathbb B}^{[0]}(M)^{U_L}={\rm Hom}\left(Y_{[L]},M\right)$. 

As ${\mathbb B}^{[0]}(M)^{U_L}={\rm Hom}\left(Y_{[L]},M\right)=\pi 
B^q(X\times_kY_{[L]})$ and $\pi B^q(X_F)^{U_L}=\pi B^q(X_L)$, 
it suffices to show that $\pi B^q(X\times_kY_{[L]})=\pi B^q(X_L)$ for 
$\pi=\Delta_{k(X)}$ and any sufficiently big (with $\dim Y_{[L]}\ge q$) 
subfield $L$ finitely generated over $k$. We may assume that $F$ is 
algebraic over $L$, so $n<\infty$. Then, by Lemma \ref{moving}, the 
natural map $Z^q(k(X)\otimes_kL)\longrightarrow B^q(X\times_kY_{[L]})$ 
is surjective, and thus, the composition $Z^q(k(X)\otimes_kF)=
\lim\limits_{_L\longrightarrow}Z^q(k(X)\otimes_kL)\longrightarrow
\lim\limits_{_L\longrightarrow}\pi B^q(X\times_kY_{[L]})=
\lim\limits_{_L\longrightarrow}{\rm Hom}\left(Y_{[L]},M\right)=
{\mathbb B}^{[0]}(M)$ is also surjective. 

This implies that there exists a natural embedding of the 
${\mathbb Q}$-algebra ${\rm End}_G{\mathbb B}^{[0]}(M)$ into the space 
${\rm Hom}_G\left(Z^q(k(X)\otimes_kF),{\mathbb B}^{[0]}(M)\right)$. 

Using Proposition \ref{max-prim} one gets 
\begin{multline*} {\rm Hom}_G\left(Z^q(k(X)\otimes_kF),
{\mathbb B}^{[0]}(M)\right)={\rm Hom}\left(Y_{[k(X)]},M\right) \\ 
={\rm Hom}\left((Y_{[k(X)]},\Delta_{k(X)}),M\right)=
\pi B^q(X\times_kX)\Delta_{k(X)}.\end{multline*} 
By Corollary \ref{descr-proj}, the latter coincides with 
$B^q(k(X)\otimes_kk(X))$, which is the same as 
${\rm End}_G\left(\pi B^q(X_F)\right)$. 

By \cite{jan} and Lemma \ref{loc-glo-ss}, the $G$-module 
${\mathbb B}^{[0]}(M)$ is semi-simple. As ${\mathbb B}^{[0]}(M)$ 
surjects onto $\pi B^q(X_F)$ and ${\rm End}_G{\mathbb B}^{[0]}(M)
\subseteq{\rm End}_G\pi B^q(X_F)$, one has ${\mathbb B}^{[0]}(M)=
\pi B^q(X_F)$. \qed 

\begin{corollary} \label{max-proj-gen} $\Delta_{k(X)}B^d
(X\times_kY)=B^d(k(X)\otimes_kk(Y))$ and the group $\Delta_{k(X)}
B^q(X\times_kY)$ vanishes for any $q<d:=\dim X$, any irreducible 
smooth proper $k$-variety $X$ and any irreducible smooth 
$k$-variety $Y$. \end{corollary} 
{\it Proof.} $\Delta_{k(X)}B^d(X\times_kY)={\rm Hom}
(Y,(X,\Delta_{k(X)}))$. By Proposition \ref{b-explic}, 
$${\rm Hom}(Y,(X,\Delta_{k(X)}))=\Delta_{k(X)}B^d\left(X_{k(Y)}
\right).$$ By Lemma \ref{max-proj}, $\Delta_{k(X)}B^d
\left(X_{k(Y)}\right)=B^d(k(X)\otimes_kk(Y))$. 

Suppose that $q<d$. Consider the projection 
$X\times_kY\times{\mathbb P}^{d-q}\stackrel{p}{\longrightarrow}
X\times_kY$. The pull-back induces an embedding $\Delta_{k(X)}
B^d(X\times_kY)\stackrel{p^{\ast}}{\longrightarrow}
\Delta_{k(X)}B^d(X\times_kY\times{\mathbb P}^{d-q})$, 
which is an isomorphism, since $\Delta_{k(X)}B^d(X\times_kY)$ 
coincides with $B^d(k(X)\otimes_kk(Y))$, and $\Delta_{k(X)}B^d
(X\times_kY\times{\mathbb P}^{d-q})=B^d(k(X)\otimes_k
k(Y\times{\mathbb P}^{d-q}))$ coincides with $B^d(k(X)\otimes_kk(Y))$. 

The push-forward induces a surjection $B^d(X\times_kY
\times{\mathbb P}^{d-q})\stackrel{p_{\ast}}{\longrightarrow}
B^q(X\times_kY)$. On the other hand, the composition $B^d
(X\times_kY)\stackrel{p_{\ast}p^{\ast}}{\longrightarrow}
B^q(X\times_kY)$ is zero, so $\Delta_{k(X)}B^q(X\times_kY)=0$. \qed 

\subsection{``Polarization'' on $B^n(k(X)\otimes_kF)$ 
and polarizable $G$-modules}
\begin{proposition} \label{polar} For any irreducible $k$-variety 
$X$ of dimension $n$ there is a symmetric $G$-equivariant 
non-degenerate pairing $$B^n(k(X)\otimes_kF)\otimes B^n(k(X)\otimes_kF)
\stackrel{\langle~,~\rangle}{\longrightarrow}{\mathbb Q}(\chi)$$ 
such that $\langle p^{\ast}(\cdot),~\cdot~\rangle=\langle~\cdot~,
p_{\ast}(\cdot)\rangle$ for any generically finite rational map $p$. 
In particular, $\langle~,~\rangle$ induces non-degenerate pairings 
between the submodules $W:=\pi B^n(k(X)\otimes_kF)$ and $~^tW:=~^t\pi 
B^n(k(X)\otimes_kF)$ for all projectors $\pi\in B^n(k(X)\otimes_kk(X))$. 

If for $(n-1)$-cycles on $2n$-dimensional complex varieties the 
numerical equivalence coincides with the homological one, then 
$\langle~,~\rangle$ is $(-1)^n$-definite. In particular, this 
holds for $n\le 2$. \end{proposition} 
{\it Proof.} We may suppose that $X$ is a smooth projective variety 
over $k$. For a pair $\alpha,\gamma\in B^n(k(X)\otimes_kF)$ fixed 
by a compact open subgroup $U\subset G$ we define $\langle\alpha,
\gamma\rangle\in{\mathbb Q}(\chi)$ by $\langle\widehat{\alpha}\cdot
\widehat{\gamma}\rangle\cdot[U]$, where $\widehat{\alpha},
\widehat{\gamma}$ are the images of $\alpha,\gamma\in 
B^n(k(X)\otimes_kk(Y_U))$ in $B^n(X\times_kY_U)$ in the sense of 
Proposition \ref{split}. Here $Y_U$ is a smooth proper variety over 
$k$ with the function field $k(Y_U)$ identified with $F^U$, and 
$\langle~\cdot~\rangle$ is the intersection form on 
$B^n(X\times_kY_U)$. By the projection formula, $\langle\alpha,
\gamma\rangle$ is independent of the choices, and $\langle 
p^{\ast}(\cdot),~\cdot~\rangle=\langle~\cdot~,p_{\ast}(\cdot)\rangle$. 

For a triplet of smooth proper varieties $X_1,X_2,X_3$, a triplet of 
integers $a,b,c\ge 0$ with $a+b+c=\dim(X_1\times X_2\times X_3)$ and 
a triplet $\alpha\in A^a(X_1\times X_2),\beta\in A^b(X_2\times X_3),
\gamma\in A^c(X_1\times X_3)$ one has $\langle\alpha\circ\beta\cdot
\gamma\rangle=\langle\beta\cdot\!\!~^t\alpha\circ\gamma\rangle=\langle
\alpha\cdot\gamma\circ\!\!~^t\beta\rangle$. For any 
$\alpha\in W-\{0\}$ fixed by $U$ there is $\beta\in B^n(X\times_kY_U)$ 
such that $\langle\widehat{\alpha}\cdot\beta\rangle\neq 0$. Then, as 
$\widehat{\alpha}=\pi\circ\widehat{\alpha}\circ\Delta_{k(Y_U)}$, 
one has 
$\langle\widehat{\alpha}\cdot\beta\rangle=\langle\widehat{\alpha}
\cdot\!~^t\pi\circ\beta\circ\Delta_{k(Y_U)}\rangle\neq 0$. As 
$~^t\pi\circ\beta\circ\Delta_{k(Y_U)}\in~^tW$, this shows that 
$\langle~,~\rangle$ induces a non-degenerate pairing between $W$ 
and $~^tW$ for an arbitrary projector $\pi$. 

By Lemma \ref{redvan}, ${\rm Hom}_{{\mathcal H}(U)}
(B^{n+1}(X\times_kY_U),B^n(k(X)\otimes_kk(Y_U)))=0$. By the 
semi-simplicity, this implies that the composition of the 
embedding of ${\mathcal H}(U)$-modules 
$B^n(k(X)\otimes_kk(Y_U))\hookrightarrow B^n(X\times_kY_U)$ 
with $B^n(X\times_kY_U)\stackrel{\cdot[L]}{\longrightarrow}
B^{n+1}(X\times_kY_U)$ is zero for any $L\in{\rm NS}(X)$. 
Interchanging the roles of $X$ and $Y_U$, we see that the image of 
$B^n(k(X)\otimes_kk(Y_U))$ in $B^n(X\times_kY_U)$ is annihilated 
by any $L\in{\rm NS}(X)\bigoplus{\rm NS}(Y_U)\subseteq
{\rm NS}(X\times_kY_U)$, i.e., it consists of primitive elements. 
Then, by the Hodge index theorem, if for $(n-1)$-cycles on 
$2n$-dimensional complex varieties the numerical equivalence coincides 
with the homological one, then the pairing $B^n(k(X)\otimes_kF)\otimes 
B^n(k(X)\otimes_kF)\stackrel{\langle~,~\rangle}{\longrightarrow}
{\mathbb Q}(\chi)$ is $(-1)^n$-definite. \qed 

\begin{proposition} \label{nt} Let $V$ be a finite-dimensional 
${\mathbb Q}$-vector space with a positive definite symmetric pairing 
$V\otimes V\stackrel{\langle~~,~~\rangle}{\longrightarrow}{\mathbb Q}$ 
and with an action $H\longrightarrow{\rm GL}(V)$, $p\longmapsto
\sigma_p$ of a subgroup $H\subseteq{\mathbb Q}^{\times}_+$, generated by 
almost all primes, such that $\langle\sigma_px,\sigma_py\rangle=
p\cdot\langle x,y\rangle$ for all $x,y\in V$ and for all $p\in H$. 
Then $V=0$. \end{proposition} 
{\it Proof.} Since a $\overline{{\mathbb Q}}$-multiple of $\sigma_p$ is 
orthogonal, $\sigma_p$ is diagonizable over $\overline{{\mathbb Q}}$ for 
all $p\in H$. As the group $H$ is abelian, there is a basis $\{e_i\}$ 
of $V\otimes\overline{{\mathbb Q}}$ and characters $\lambda_i:H
\longrightarrow\overline{{\mathbb Q}}^{\times}$ such that 
$\sigma_re_i=\lambda_i(r)\cdot e_i$. Note, that the elements $e_i$ 
belong to $V\otimes K$ for a finite extension $K$ of ${\mathbb Q}$, so 
the characters factor as $\lambda_i:H\longrightarrow K^{\times}$. 

For each embedding $\tau:K\hookrightarrow{\mathbb C}$ we define an 
hermitian form $\langle x,y\rangle_{\tau}:=\langle\tau x,
\overline{\tau y}\rangle$, where $\langle~~,~~\rangle$ is the 
bilinear form on $V\otimes{\mathbb C}$ induced by $\langle~~,~~\rangle$ 
and $\overline{\phantom{ta}}$ is the complex conjugation on 
$V\otimes{\mathbb C}$. One has $0\neq\langle\sigma_pe_i,\sigma_pe_i
\rangle_{\tau}=\tau(\lambda_i(p))\overline{\tau(\lambda_i(p))}
\cdot\langle e_i,e_i\rangle_{\tau}=p\cdot\langle e_i,e_i
\rangle_{\tau}$, so $|\tau(\lambda_i(p))|=p^{1/2}$ for any $\tau$. 

Let $L$ be the subfield in the normalization of $K$ over ${\mathbb Q}$ 
generated by all the conjugates over ${\mathbb Q}$ of the image of 
$\lambda_1$, so $L$ is a Galois extension of ${\mathbb Q}$. For any 
embedding $\tau:K\hookrightarrow{\mathbb C}$ the image of $L$ in 
${\mathbb C}$ is invariant under the complex conjugation, since 
$\overline{\tau\lambda_1(p)}=p\cdot(\tau\lambda_1(p))^{-1}$. 
As the latter does not depend on the embedding, the complex 
conjugation induces an element $c$ in the center of 
${\rm Gal}(L/{\mathbb Q})$. $L$ cannot be totally real, since then it 
would contain elements $p^{1/2}$ for almost all prime $p$, and thus 
it would be of infinite degree. As the field of invariants $R$ of 
$c$ is totally real, $L=R(\sqrt{-\alpha})$ for some totally positive 
integer $\alpha\in{\mathcal O}_R$. Let $d=[R:{\mathbb Q}]$. 

As the subgroup $\{1,c\}$ is normal in ${\rm Gal}(L/{\mathbb Q})$, 
the field $R$ is a Galois extension of ${\mathbb Q}$. By the Chebotar{\"e}v 
density theorem, there are infinitely many rational primes $p$ 
corresponding to the element (=the conjugacy class of) $c$. As 
restriction of $c$ to $R$ is trivial, such ideals $(p)$ split 
completely in the extension $R/{\mathbb Q}$, i.e., $(p)=\wp_1\cdots\wp_d$. 
On the other hand the ideals $\wp_1,\dots,\wp_d$ 
stay prime in the extension $L/R$ (one of them should stay 
prime, since $c$ is non-trivial, and the others lie in the 
${\rm Gal}(L/{\mathbb Q})$-orbit of such one). 

Now let integers $x,y\in{\mathcal O}_R$ be such that 
$x^2+\alpha y^2\in(p)\subseteq\wp_j$ for any $j$. Since $\wp_j$ 
remains prime in ${\mathcal O}_L$, one has $x+y\sqrt{-\alpha}\in\wp_j$ 
for a choice of $\sqrt{-\alpha}$. Since $\wp_j$ is 
$c$-invariant, $2x,2y\sqrt{-\alpha}\in\wp_j$, and thus, 
$x,y\in\wp_j$ for $p$ big enough with respect to $\alpha$. 
So we get $x,y\in\wp_1\cdots\wp_d=(p)$. 

If there is an element in $L$ with the modulus $p^{1/2}$ with respect to 
any complex embedding of $L$, then there are non-zero integers $x,y,z\in
{\mathcal O}_R$ of minimal possible norm of $xyz$ such that $x^2+\alpha y^2=
pz^2$. As this implies $z^2\in(p)$, one should also have $z\in(p)$ if $L$ 
is unramified over $p$, and therefore, the triplet $(x/p,y/p,z/p)$ also 
satisfies the above conditions and has a smaller ``norm''. This is 
contradiction. \qed 

\vspace{5mm}

We shall say that an admissible $G$-module $W$ is {\sl polarizable} 
(when $n<\infty$) if there is a symmetric positive definite 
$G$-equivariant pairing $W\otimes W\longrightarrow{\mathbb Q}(\chi)$. 
\begin{corollary} \label{hyperspe} Let $W$ be a polarizable 
$G$-module. Let $L$ be an extension of $k$ in $F$. Then $W^{U_{L(x)}}=0$ 
for any $x\in F$ transcendental over $L$. \end{corollary} 
{\it Proof.} Since $W$ is smooth, there is a finitely generated 
extension $L_1$ of $k$ such that the stabilizer of an element 
$w\in W^{U_{L(x)}}$ contains the open subgroup $\langle U_{L_1},
U_{L(x)}\rangle$. By Lemma \ref{2.0}, the latter contains 
$U_{L_2(x)\cap L_3(x)}=U_{L_4(x)}$, where $L_3$ is generated over $L$ 
by a transcendence basis of $F$ over $L(x)$, $L_2$ is generated over 
$L_1$ by a transcendence basis of $L_3$ over $k$, and $L_4:=L_2(x)
\bigcap L_3$ is finitely generated over $k$. We replace $L$ by 
$L_4$, thus assuming, that $L$ is finitely generated over $k$. 

As $U_{L(x)}\subset U_{L(x^p)}=\sigma U_{L(x)}\sigma^{-1}$, where 
$\sigma x=x^p$ and $\sigma|_L=id$, the element $\sigma$ induces 
an isomorphism $W^{U_{L(x)}}\stackrel{\sim}{\longrightarrow}
W^{U_{L(x^p)}}$, and $W^{U_{L(x^p)}}\subseteq W^{U_{L(x)}}$, the 
dimension argument shows that $\sigma$ induces an automorphism 
of $W^{U_{L(x)}}$. For any $w\in W$ one has $\langle\sigma w,
\sigma w\rangle=\chi(\sigma)\cdot\langle w,w\rangle$, so 
Proposition \ref{nt} says that $W^{U_{L(x)}}=0$. \qed

\begin{corollary} \label{no-fi-di} Any finite-dimensional 
polarizable $G$-module is zero. \end{corollary} 
{\it Proof.} This follows from Corollary \ref{liss-fin-inf} 
and Corollary \ref{hyperspe}. \qed

\section{Morphisms between certain $G$-modules and matrix coefficients} 
\subsection{Two remarks on the $G$-modules $F/k$ and 
$F^{\times}/k^{\times}$} 
\label{muladd}
\begin{proposition} \label{add-mult} For any $1\le n\le\infty$ 
the $G^{\circ}$-modules $F/k$ and $F^{\times}/k^{\times}$ are 
irreducible. \end{proposition} 
{\it Proof.} Let $A$ be the additive subgroup of $F$ generated by 
the $G^{\circ}$-orbit of some $x\in F-k$. For any $y\in A-k$ one has 
$\frac{2}{y^2-1}=\frac{1}{y-1}-\frac{1}{y+1}$. As $\frac{1}{y-1}$ 
and $\frac{1}{y+1}$ are in the $G^{\circ}$-orbit of $y$, this 
implies that $y^2\in A$. As for any $y,z\in A$ one has 
$yz=\frac{1}{4}((y+z)^2-(y-z)^2)$, the group $A$ is a subring of $F$. 

Let $M$ be the multiplicative subgroup of $F^{\times}$ generated by 
the $G^{\circ}$-orbit of some $x\in F-k$. Then for any $y,z\in M$ 
one has $y+z=z(y/z+1)$, so if $y/z\not\in k$ then $y+z\in M$, and 
thus, $M\bigcup\{0\}$ is a $G^{\circ}$-invariant subring of $F$. 

Since the $G^{\circ}$-orbit of an element $x\in F-k$ contains all 
elements of $F-\overline{k(x)}$, if $n\ge 2$ then each element 
of $F$ is the sum of a pair of elements in the orbit. Any 
$G^{\circ}$-invariant subring in $F$, but not in $k$, is a 
$k$-subalgebra, so if $n=1$ then ${\rm Gal}(F/{\mathbb Q}(G^{\circ}x))
\subset G^{\circ}$ is a compact subgroup normalized by 
$G^{\circ}$. Then by Theorem \ref{exactness} we have 
${\rm Gal}(F/{\mathbb Q}(G^{\circ}x))=\{1\}$. As any element of 
${\mathbb Q}(G^{\circ}x)$ is the fraction of a pair of elements in 
${\mathbb Z}[G^{\circ}x]$ and for any $y\in F-k$ the element $1/y$ 
belongs to the $G^{\circ}$-orbit of $y$, one has 
${\mathbb Z}[G^{\circ}x]=F$. \qed 

\begin{proposition} \label{ann-nul} For any $1\le n\le\infty$ 
the annihilator of $F/k$ in ${\bf D}_k$ and the annihilator of 
$F^{\times}/k^{\times}$ in ${\bf D}_{{\mathbb Q}}$ are trivial. 
\end{proposition} 
{\it Proof.} Let $\sum_{j=1}^Na_j\cdot\sigma_j$ be the image of some 
$\alpha\in{\bf D}_E$ in $E[G/U_L]$ for an open subgroup $U_L$, an 
integer $N\ge 1$, some $a_j\in E$ and some pairwise distinct 
embeddings $\sigma_j:L\stackrel{/k}{\hookrightarrow}F$. Here $E$ is 
either $k$, or ${\mathbb Q}$. Let $\tau\in G$ be such an element that the 
embeddings $\sigma_1,\dots,\sigma_N,\tau\sigma_1,\dots,\tau\sigma_N$ 
are pairwise distinct. 

Suppose $\alpha$ annihilates $F/k$. For any $x\in L$ one has 
$\sum_{j=1}^Na_j\cdot\sigma_jx\in k$, and thus, 
$\sum_{j=1}^Na_j\cdot\sigma_jx+ 
\sum_{j=1}^N(-a_j)\cdot\tau\sigma_jx=0$ for all $x\in L$. Then, by 
Artin's theorem on independence of characters, $a_1=\dots=a_N=0$. 

Similarly, suppose that $\alpha$ annihilates $F^{\times}/k^{\times}$. 
We may suppose that $a_j\in{\mathbb Z}$. 
For any $x\in L^{\times}$ one has $\prod_{j=1}^N(\sigma_jx)^{a_j}\in 
k^{\times}$, and thus, $\prod_{j=1}^{2N}
(\sigma_jx)^{a_j}=1$ for all $x\in L^{\times}$, where 
$\sigma_j=\tau\sigma_{j-N}$ and $a_j=-a_{j-N}$ for $N<j\le 2N$. 

Let $\tau_1,\dots,\tau_M$ be a collection of pairwise distinct 
embeddings of $L$ into $F$ and $b_1,\dots,b_M$ be such a collection of 
non-zero integers that $\sum_{j=1}^M|b_j|$ is minimal among those for 
which $\prod_{j=1}^M(\tau_jx)^{b_j}=1$ for all $x\in L^{\times}$. 
Then $\prod_{j=1}^M(\tau_jx+1)^{b_j}-\prod_{j=1}^M(\tau_jx)^{b_j}=0$, 
which is equivalent to $$\prod_{j:b_j>0}(\tau_jx+1)^{b_j}
\prod_{j:b_j<0}(\tau_jx)^{-b_j}-\prod_{j:b_j>0}(\tau_jx)^{b_j}
\prod_{j:b_j<0}(\tau_jx+1)^{-b_j}=0.$$ We rewrite this as 
$\sum_{c_1,\dots,c_M}A_{c_1,\dots,c_M}\prod_{j=1}^M(\tau_jx)^{c_j}=0$, 
where $0\le c_j\le|b_j|$ and $\sum_{j=1}^Mc_j<\sum_{j=1}^M|b_j|$. Then 
$A_{c_1,\dots,c_M}$ are integers which are all non-zero if either 
$c_j=b_j$ for all $j$ with $b_j>0$, or $c_j=-b_j$ for all $j$ with 
$b_j<0$. By Artin's theorem on independence of characters, 
$\prod_{j=1}^M(\tau_jx)^{c_j}=\prod_{j=1}^M(\tau_jx)^{c'_j}$ for 
a pair of distinct collections $(c_j)$ and $(c'_j)$ as before, 
with $|c_j-c'_j|<|b_j|$ for some $j$, and for all $x\in L$. 
But then $\sum_{j=1}^M|c_j-c'_j|<\sum_{j=1}^M|b_j|$, 
contradicting our assumptions, so $a_1=\dots=a_N=0$. \qed 

\subsection{Morphisms between certain $G$- and $G^{\circ}$-modules} 
\begin{proposition} \label{g-g-knought} Let ${\mathcal A}$ 
and ${\mathcal B}$ be simple commutative group schemes over $k$. 
Then for any $1\le n\le\infty$ 
$${\rm Hom}_{{\rm group~schemes}/k}({\mathcal A},{\mathcal B}
)_{{\mathbb Q}}\stackrel{\sim}{\longrightarrow}{\rm Hom}_G(W_{{\mathcal A}},
W_{{\mathcal B}})\stackrel{\sim}{\longrightarrow}{\rm Hom}_{G^{\circ}}
(W_{{\mathcal A}},W_{{\mathcal B}}).$$ \end{proposition}
{\it Proof.} \begin{enumerate} \item First, consider the cases 
${\mathcal A},{\mathcal B}\in\{{\mathbb G}_a,{\mathbb G}_m\}$. 
Let $\varphi\in{\rm Hom}_{G^{\circ}}(W_{{\mathcal A}},W_{{\mathcal B}})$ 
and $x\in W_{{\mathcal A}}-\{0\}$. Then $\varphi(x)$ is fixed by 
the stabilizer ${\rm St}_x$ of $x$ in $G^{\circ}$. The group ${\rm St}_x$ 
fits into an exact sequence $1\longrightarrow U_{k(x)}\bigcap G^{\circ}
\longrightarrow{\rm St}_x\longrightarrow{\mathcal A}(k)\longrightarrow 1$. 
As $W_{{\mathcal B}}^{U_{k(x)}\bigcap G^{\circ}}=({\mathcal B}(k(x))/
{\mathcal B}(k))_{{\mathbb Q}}$ and ${\mathcal A}(k)$ acts on $k(x)$ by 
the affine linear substitutions of $x$, one has 
$W_{{\mathcal B}}^{{\rm St}_x}=
({\mathcal B}(k(x))/{\mathcal B}(k))_{{\mathbb Q}}^{{\mathcal A}(k)}$. 
\begin{itemize} \item If ${\mathcal A}={\mathbb G}_m$ and 
${\mathcal B}={\mathbb G}_a$ 
then $W_{{\mathcal B}}^{{\rm St}_x}=(k(x)/k)^{k^{\times}}=0$, so $\varphi=0$. 
\item If ${\mathcal A}={\mathbb G}_a$ and ${\mathcal B}={\mathbb G}_m$ then 
$W_{{\mathcal B}}^{{\rm St}_x}=(k(x)^{\times}/k^{\times})^k_{{\mathbb Q}}=0$, 
so $\varphi=0$. 
\item If ${\mathcal A}={\mathcal B}={\mathbb G}_m$ then 
$W_{{\mathcal B}}^{{\rm St}_x}
=(k(x)^{\times}/k^{\times})^{k^{\times}}_{{\mathbb Q}}=
\{x^{\lambda}~|~\lambda\in{\mathbb Q}\}$, so $\varphi(x)=x^{\lambda}$ for 
some $\lambda\in{\mathbb Q}$. This implies that 
$\varphi(\sigma x)=(\sigma x)^{\lambda}$ for all $\sigma\in G^{\circ}$. 
As $F^{\times}/k^{\times}$ is an irreducible $G^{\circ}$-module, one has 
$\varphi y=y^{\lambda}$ for all $y\in F^{\times}/k^{\times}$.  
\item If ${\mathcal A}={\mathcal B}={\mathbb G}_a$ then 
$W_{{\mathcal B}}^{{\rm St}_x}
=(k(x)/k)^k=\{\lambda\cdot x~|~\lambda\in k\}$, so 
$\varphi(x)=\lambda\cdot x$ for some $\lambda\in k$. This implies that 
$\varphi(\sigma x)=\lambda\cdot\sigma x$ for all $\sigma\in G^{\circ}$. 
As $F/k$ is an irreducible $G^{\circ}$-module, one has 
$\varphi y=\lambda\cdot y$ for all $y\in F/k$. 
\end{itemize}
\item Fix a smooth irreducible curve $Z\subseteq{\mathcal A}$ over 
$k$. For any generic point $x:k(Z)\stackrel{/k}{\hookrightarrow}F$ 
and any $G^{\circ}$-homomorphism $W_{{\mathcal A}}
\stackrel{\varphi}{\longrightarrow}W_{{\mathcal B}}$ the element 
$\varphi(x)$ belongs to the space 
$W_{{\mathcal B}}^{U_{x(k(Z))}\cap G^{\circ}}
=({\mathcal B}(x(k(Z)))/{\mathcal B}(k))_{{\mathbb Q}}$, and therefore, 
there is an integer $N_x\ge 1$ and a rational map $Z
\stackrel{h_x}{--\rightarrow}{\mathcal B}$ defined over $k$ such that 
$\varphi(x)=\frac{1}{N_x}h_x(x)$. This implies that if $\varphi\neq 0$
and either ${\mathcal A}={\mathbb G}_m$ or ${\mathcal A}=
{\mathbb G}_a$ then 
${\mathcal B}={\mathbb G}_m$ or ${\mathcal B}={\mathbb G}_a$. 
\item Now suppose that ${\mathcal A}$ and ${\mathcal B}$ are simple 
abelian varieties over $k$. Then the maps $h_x$ are regular and factor as 
$Z\longrightarrow J_Z\stackrel{\varphi_x}{\longrightarrow}{\mathcal B}$ 
for some homomorphisms $\varphi_x$ from the Jacobian $J_Z$. 
\item If ${\mathcal A}$ and ${\mathcal B}$ are not isogeneous, there 
is such a curve $Z$ that any homomorphism of the group schemes $J_Z
\longrightarrow{\mathcal B}$ is trivial, and thus $\varphi(x)=0$, 
so we may suppose that ${\mathcal A}={\mathcal B}$. 
\item There is such a curve $Z$ that $${\rm End}_{{\mathcal A}}:=
{\rm End}_{{\rm group~scheme}/k}({\mathcal A})_{{\mathbb Q}}
\stackrel{\sim}{\longrightarrow}{\rm Hom}_{{\rm group~schemes}/k}
(J_Z,{\mathcal A})_{{\mathbb Q}}.$$ Then one can consider $\varphi_x$ 
as an element of ${\rm End}_{{\mathcal A}}$, i.e., 
$\varphi(x)=\frac{1}{N_x}\varphi_x(x)$. 
\item For any pair of generic points $x_1,x_2\in Z(F)$ there are 
generic points $$t,y_1,\dots,y_M,z_1,\dots,z_{M'}$$ of $Z(F)$ such that 
the elements of both collections $(x_1,t,y_1,\dots,y_M)$ and $(x_2,
t,z_1,\dots,z_{M'})$ are linearly independent in $W_{{\mathcal A}}$ 
over the algebra ${\rm End}_{{\mathcal A}}$, and $u:=x_1+t+\sum_jy_j$ 
and $v:=x_2+t+\sum_jz_j$ are generic points of $Z(F)$. 

Then, by definition, $\varphi(u)$ coincides with $$\frac{1}{N_u}
\varphi_u(u)=\frac{1}{N_u}\left(\varphi_u(x_1)+\varphi_u(t)+\sum_j
\varphi_u(y_j)\right),$$ and, as $\varphi$ is a homomorphism, 
$\varphi=\frac{1}{N_{x_1}}\varphi_{x_1}(x_1)+ 
\frac{1}{N_t}\varphi_t(t)+\frac{1}{N_{y_1}}\varphi_{y_1}(y_1)+
\dots+\frac{1}{N_{y_M}}\varphi_{y_M}(y_M)$. We can assume that 
$N_u=N_v=N_{x_1}=N_{x_2}=N_t$. Then one has 
$\varphi_u=\varphi_{x_1}=\varphi_t$. Similarly, one has 
$\varphi_v=\varphi_{x_2}=\varphi_t$, and therefore, the 
restriction of $\varphi$ to the set of generic points of 
$Z$ coincides with the restriction of $\psi$ to the set of 
generic points of $Z$, for some $\psi\in{\rm End}_{{\mathcal A}}$. 

As generic points of $Z$ generate ${\mathcal A}(F)$ as an abstract 
group, this implies that $\varphi=\psi$, i.e., that 
$\varphi\in{\rm End}_{{\mathcal A}}$. 
\item There remain the cases where ${\mathcal A}$ is an abelian variety 
and ${\mathcal B}$ is either ${\mathbb G}_a$ or ${\mathbb G}_m$. 
As $W_{{\mathcal B}}$ is an irreducible $G^{\circ}$-module, any non-zero 
$G^{\circ}$-homomorphism to $W_{{\mathcal B}}$ is surjective. Any 
surjection of smooth representations of $G^{\circ}$ induces a surjection 
of their subspaces fixed by a compact subgroup $K$ of $G$. Taking 
$K=U_{L'}$, where $L'$ is a purely transcendental extension 
of $k$, we get contradiction showing that 
${\rm Hom}_{G^{\circ}}(W_{{\mathcal A}},W_{{\mathcal B}})=0$. \qed 
\end{enumerate}

\begin{corollary} \label{for-fai} For any pair of pure $1$-motives 
$M_1,M_2$ one has $${\rm Hom}_G({\mathfrak B}^1(M_1),
{\mathfrak B}^1(M_2))={\rm Hom}_{G^{\circ}}({\mathfrak B}^1(M_1),
{\mathfrak B}^1(M_2)).\qquad\qed$$ \end{corollary}

\begin{corollary} \label{0-1} For any $G$-module $W$ there is 
at most $1$ character $\psi$ such that $W(\psi)\cong{\mathfrak B}^1(M)$ 
for a pure $1$-motive $M$. \end{corollary} 
{\it Proof.} Suppose that $W$ is irreducible, $W(\psi_1)\cong
{\mathfrak B}^1(M_1)$, $W(\psi_2)\cong{\mathfrak B}^1(M_2)$ for some pure 
$1$-motives $M_1$ and $M_2$ and some characters $\psi_1\neq\psi_2$. 
By the fully faithfulness and the previous corollary, 
$${\rm Hom}(M_1,M_2)={\rm Hom}_G({\mathfrak B}^1(M_1),
{\mathfrak B}^1(M_2))={\rm Hom}_{G^{\circ}}({\mathfrak B}^1(M_1),
{\mathfrak B}^1(M_2)),$$ so 
${\rm Hom}(M_1,M_2)={\rm Hom}_{G^{\circ}}(W(\psi_1),W(\psi_2))$, 
which non-zero, since the $G^{\circ}$-modules $W(\psi_1)$ and 
$W(\psi_2)$ are isomorphic, and thus, $M_1$ and $M_2$ are isomorphic, 
which implies that $W(\psi_1)\cong W(\psi_2)$ as $G$-modules. We may 
assume that $\psi_2=1$. Set $\psi=\psi_1\neq 1$. 

The $G$-module $W(\psi)$ coincides with $W$ as a vector space, 
but $G$ acts by $(\sigma,w)\longmapsto\psi(\sigma)\cdot\sigma w$. 
Suppose that there is an isomorphism $W\longrightarrow W(\psi)$, 
i.e., an automorphism $W\stackrel{\lambda}{\longrightarrow}W$ 
such that $\lambda(\sigma w)=\psi(\sigma)\cdot\sigma\lambda(w)$. 
Then $\lambda$ can be considered as an element of 
${\rm End}_{\ker\psi}(W)$. As $\ker\psi$ contains $G^{\circ}$, 
the automorphism $\lambda$ can also be considered as an element 
of ${\rm End}_{G^{\circ}}(W)$. 

This implies that ${\rm End}_{G^{\circ}}(W)\neq{\rm End}_G(W)$, 
contradicting the previous corollary. 

If $W$ is not irreducible, it should be semi-simple anyway, so 
its irreducible summands are motivic, implying the Corollary. \qed 

\subsection{Non-compactness of supports of matrix coefficients}
\label{geomor}
A {\sl matrix coefficient} of a smooth representation $W$ 
of a topological group is a function on the group of type 
$\langle\sigma w,\widetilde{w}\rangle$ for a vector $w\in W$ 
and a vector $\widetilde{w}$ in the dual representation with 
open stabilizer. 
\begin{proposition} Suppose that $n<\infty$, a subgroup $H$ of $G$ 
contains $G^{\circ}$ and the supports of the matrix coefficients of 
a representation of $H$ are compact. Then this representation is 
zero. \end{proposition} 
{\it Proof.} Let ${\mathbb C}$ be an algebraically closed extension 
of $k$ of cardinality strictly greater than the cardinality of $k$. 
Let $\omega\neq 0$ be an irreducible representation of $H$ over 
${\mathbb C}$ with compact supports of the matrix coefficients. 
We may replace $H$ by $G^{\circ}$, and replace $\omega$ by an 
irreducible subquotient of $\omega|_{G^{\circ}}$. Let $N$ be a 
compact open subgroup in $G^{\circ}$ such that $\omega(h_N)\neq 0$. 

Along the same lines as, e.g., in Claim 2.11 of \cite{bz}, 
one proves the Schur's lemma: the endomorphisms of a smooth 
irreducible ${\mathbb C}$-representation of $G^{\circ}$ are scalar. 
This allows one to modify Theorem 2.42 a) of \cite{bz} as follows. 
{\it For each open compact subgroup $N$ in $G^{\circ}$ 
there is an element $\varepsilon^{\omega}_N\in{\bf D}_{{\mathbb C}}$ 
such that $\omega(\varepsilon^{\omega}_N)=\omega(h_N)$ and 
$\pi(\varepsilon^{\omega}_N)=0$ for any smooth irreducible 
${\mathbb C}$-representation $\pi$ of $G^{\circ}$ distinct from $\omega$.} 

As the $G^{\circ}$-module $F^{\times}/k^{\times}$ is irreducible, 
and ${\rm End}_{G^{\circ}}(F^{\times}/k^{\times})={\mathbb Q}$, the 
${\mathbb C}$-representation $\pi=(F^{\times}/k^{\times})\otimes
{\mathbb C}$ of $G^{\circ}$ is irreducible.\footnote{Variant: $F/k$ is an 
irreducible $G^{\circ}$-module, ${\bf D}_{{\mathbb Q}}\hookrightarrow
{\rm End}_k(F/k)$, ${\rm End}_{G^{\circ}}(F/k)=k$, and therefore, $(F/k)
\otimes_k{\mathbb C}$ is an irreducible ${\mathbb C}[G^{\circ}]$-module.} 
On the other hand, the support of the 
matrix coefficient $\langle\sigma x,\widetilde{w}\rangle$ is not 
compact for any $x\in F^{\times}/k^{\times}-\{1\}$ and any vector 
$\widetilde{w}\neq 0$ in the representation dual to $\pi$, 
since the stabilizer of $x$ is not compact. This implies that 
there is an element $\varepsilon\in{\bf D}_{{\mathbb C}}$ such that 
$\omega(\varepsilon)=\omega(h_N)$ and $\pi(\varepsilon)=0$, 
contradicting Proposition \ref{ann-nul}. \qed 

\section{Some examples of (co-)homological calculations} 
\subsection{Examples of ${\rm Ext}$-calculation and torsors} 
\label{ext-examp} Let ${\mathcal S}m_G$ be the category of smooth 
$G$-modules. It is a full abelian subcategory in the category 
of $G$-modules. 
\begin{proposition} \label{Hilbert-90} Let $n=\infty$ and 
${\mathcal A}$ be an irreducible commutative algebraic group over $k$. 

Then ${\rm Ext}^1_{{\mathcal S}m_G}({\mathcal A}(F)_{{\mathbb Q}},
{\mathbb Q})={\rm Ext}^1_{G,{\rm cont}}
({\mathcal A}(F)_{{\mathbb Q}},{\mathbb Q})=0$. \end{proposition} 
{\it Proof.} Let $0\longrightarrow{\mathbb Q}\longrightarrow E
\longrightarrow{\mathcal A}(F)_{{\mathbb Q}}\longrightarrow 0$ be a 
continuous extension, i.e., with closed stabilizers. A choice of 
a linear section ${\mathcal A}(F)_{{\mathbb Q}}\stackrel{s}{\hookrightarrow}E$ 
defines a splitting $E\cong{\mathbb Q}\bigoplus
{\mathcal A}(F)_{{\mathbb Q}}$ as a ${\mathbb Q}$-vector space. The 
$G$-action is given by $\sigma(b,x)=(b+a_{\sigma}(x),\sigma x)$, 
where $a_{\sigma}(x+y)=a_{\sigma}(x)+a_{\sigma}(y)$ and $a_{\sigma}
(\tau x)+a_{\tau}(x)=a_{\sigma\tau}(x)$. $E$ is continuous if and only 
if the subgroup $\{\sigma\in{\rm Stab}_x~|~a_{\sigma}(x)=0\}$ is 
closed for any $x\in{\mathcal A}(F)_{{\mathbb Q}}$. In particular, the map 
${\rm Stab}_x\stackrel{a(x)}{\longrightarrow}{\mathbb Q}$ given by 
$\sigma\longmapsto a_{\sigma}(x)$ is a homomorphism with closed kernel. 
As $a(x)$ is continuous, the image of any compact subgroup of 
$U_{k(x)}$ is a compact subgroup in ${\mathbb Q}$, i.e., 0. 
By Lemma \ref{inf-cl}, the subgroup generated by compact subgroups 
is dense in $U_{k(x)}$, so $U_{k(x)}$ is in the kernel of $a(x)$, 
and thus, $a(x)$ factors through ${\rm Stab}_x/U_{k(x)}\subseteq
{\mathcal A}(k)_{{\rm tors}}$. Since any homomorphism from any torsion 
group to ${\mathbb Q}$ is zero, we get $a_{\sigma}(x)=0$ for any 
$\sigma\in{\rm Stab}_x$, and therefore, 
${\rm Stab}_{(b,y)}={\rm Stab}_y$ for any $(b,y)\in E$. This implies 
also that for any $y\in{\mathcal A}(F)_{{\mathbb Q}}$, any $\tau\in G$ and 
any $\sigma\in{\rm Stab}_y$ one has $a_{\tau\sigma}(y)=a_{\tau}(y)$. 
In particular, if $x$ is a generic point of ${\mathcal A}$ and 
$H=\{\tau\in G~|~\tau x=\mu\cdot x~{\rm for}~{\rm some}
~\mu\in{\mathbb Q}^{\times}\}$, then $a(x):H\longrightarrow{\mathbb Q}$ 
factors through ${\mathbb Q}^{\times}\longrightarrow{\mathbb Q}$ and 
$p\cdot a_{\tau}(x)+a_{\sigma}(x)=a_{\tau\sigma}(x)$ for any 
$\sigma,\tau\in G$ such that $\sigma x=p\cdot x$ and 
$\tau x=q\cdot x$ for some $p,q\in{\mathbb Q}^{\times}$. Then 
$a_{\sigma}(x)=(p-1)\cdot c(x)$ for some $c(x)\in{\mathbb Q}$. 
Clearly, $c(m\cdot x)=m\cdot c(x)$ for any $m\in{\mathbb Q}^{\times}$. 

Note, that $c(\phantom{x})$ is linear on the set of generic points of 
${\mathcal A}$. Indeed, if $x$ and $y$ (considered as embeddings $x,y:
k({\mathcal A})\stackrel{/k}{\hookrightarrow}F$) are algebraically 
independent over $k$ (i.e., $k(x,y):=x(k({\mathcal A}))y(k({\mathcal A}))
\subset F$ is of transcendence degree $2\dim{\mathcal A}$ over $k$) then 
there is an element $\sigma\in G$ such that $\sigma x=2\cdot x$ and 
$\sigma y=2\cdot y$, so $a_{\sigma}(x)=c(x)$, $a_{\sigma}(y)=c(y)$, 
and $a_{\sigma}(x+y)=c(x+y)$, so, by additivity of $a_{\sigma}$, one 
has $c(x+y)=c(x)+c(y)$. In general, for any collection of generic 
points $x_1,\dots,x_N$ there is some $z\in{\mathcal A}(F)$ such that the 
subfield $k(z,x_1,\dots,x_N)$ of $F$ is of transcendence degree 
$\dim{\mathcal A}$ over $k(x_1,\dots,x_N)$. By induction on $N$, one has 
$\sum_{j=1}^Nm_j\cdot c(x_j)=c(z+\sum_{j=1}^Nm_j\cdot x_j)-c(z)$. 
In particular, if $\sum_{j=1}^Nm_j\cdot x_j=0$ this means that 
$\sum_{j=1}^Nm_j\cdot c(x_j)=0$. This implies that $c(\phantom{x})$ 
extends to a linear functional on ${\mathcal A}(F)_{{\mathbb Q}}$. 

Subtracting the coboundary of $c(\phantom{x})$, we may assume that 
$a_{\sigma}(y)=0$ for any generic point $y\in{\mathcal A}(F)$ and any 
$\sigma\in G$ such that $\sigma y=\mu\cdot y$ for some 
$\mu\in{\mathbb Q}^{\times}$. 

Fix a generic point $y\in{\mathcal A}(F)$. Any $G$-equivariant section 
over the subset of generic points of ${\mathcal A}$ is of type $\sigma y
\stackrel{\alpha_{y,b}}{\longmapsto}\sigma(b,y)$ for some 
$b\in{\mathbb Q}$. As $G$ acts transitively on the set of generic 
points of ${\mathcal A}$, and the stabilizers of vectors in $E$ coincide 
with the stabilizers of their projections to ${\mathcal A}(F)_{{\mathbb Q}}$, 
this section is well-defined. For any $\mu\in{\mathbb Q}^{\times}$ let 
$\tau_{\mu}\in G$ be such an element that $\tau_{\mu}y=\mu\cdot y$. 
Then $\alpha_{y,b}(\mu\cdot\sigma y)=\alpha_{y,b}(\sigma\tau_{\mu}y)=
\sigma\tau_{\mu}\alpha_{y,b}(y)=\sigma\tau_{\mu}(b,y)=
(b+a_{\sigma\tau_{\mu}}(y),\sigma\tau_{\mu}y)=(b+a_{\sigma}(\tau_{\mu}y),
\sigma\tau_{\mu}y)=\mu\cdot(b/\mu+a_{\sigma}(y),\sigma y)$, since 
$a_{\tau_{\mu}}(y)=0$, so $\alpha_{y,b}(\mu\cdot x)=
\mu\cdot\alpha_{y,b}(x)$ for any generic point $x\in{\mathcal A}(F)$ 
and any $\mu\in{\mathbb Q}^{\times}$ if and only if $b=0$, i.e., 
$\alpha_{y,0}$ is the unique $G$-equivariant homogeneous (but, 
a priori, non-linear) section over the subset of generic points of 
${\mathcal A}$. As $\alpha_{\sigma y,0}$ is also a $G$-equivariant 
homogeneous section over the subset of generic points of ${\mathcal A}$ 
for arbitrary $\sigma\in G$, one has $\alpha_{y,0}(\sigma y)=
(a_{\sigma}(y),\sigma y)=\alpha_{\sigma y,0}(\sigma y)=(0,\sigma y)$, 
so $a_{\sigma}(y)=0$ for any $\sigma\in G$ and any generic point $y$. 
Since any element of ${\mathcal A}(F)_{{\mathbb Q}}$ is a sum of generic 
points of ${\mathcal A}$, we get $a_{\sigma}(z)=0$ for any $\sigma\in G$ 
and any $z\in{\mathcal A}(F)_{{\mathbb Q}}$. \qed

\begin{corollary} \label{Ext-1} Let $n=\infty$ and 
${\mathcal A}$ be an irreducible commutative algebraic group over $k$. 
Then one has 
$${\rm Ext}^1_{{\mathcal S}m_G}(W_{{\mathcal A}},{\mathbb Q})=
{\rm Ext}^1_{G,{\rm cont}}(W_{{\mathcal A}},{\mathbb Q})
={\rm Hom}({\mathcal A}(k),{\mathbb Q}).$$ \end{corollary} 
{\it Proof.} The functor ${\rm RHom}(-,{\mathbb Q})$ applied to 
$0\longrightarrow{\mathcal A}(k)_{{\mathbb Q}}\longrightarrow
{\mathcal A}(F)_{{\mathbb Q}}\longrightarrow W_{{\mathcal A}}
\longrightarrow 0$ gives an exact sequence 
${\rm Hom}({\mathcal A}(F)_{{\mathbb Q}},{\mathbb Q})
\longrightarrow{\rm Hom}({\mathcal A}(k)_{{\mathbb Q}},{\mathbb Q})
\longrightarrow{\rm Ext}^1(W_{\mathcal A},{\mathbb Q})\longrightarrow
{\rm Ext}^1({\mathcal A}(F)_{{\mathbb Q}},{\mathbb Q})$, where the exterior 
groups are zero by Proposition \ref{Hilbert-90}. \qed 

\begin{lemma} Let $n=\infty$, and either ${\mathcal A}={\mathbb G}_m$, 
or ${\mathcal A}={\mathbb G}_a$. 

Then ${\rm Ext}^1_{{\mathcal S}m_G}
({\mathbb Q},W_{{\mathcal A}})=0$. \end{lemma}
{\it Proof.} Let $0\longrightarrow W_{{\mathcal A}}\longrightarrow E
\longrightarrow{\mathbb Q}\longrightarrow 0$ be an exact sequence in 
${\mathcal S}m_G$. The stabilizer of an element of $E$ projecting to 
$1\in{\mathbb Q}$ contains an open subgroup $U_L$, so $E$ is a sum of 
$W_{{\mathcal A}}$ and a quotient of ${\mathbb Q}[G/U_L]$ by a submodule 
in the group ${\mathbb Q}[G/U_L]^{\circ}$ of degree-zero 0-cycles. 

When restricted to ${\mathbb Q}[G/U_L]^{\circ}$, the projection 
${\mathbb Q}[G/U_L]\longrightarrow E$ factors through $W_{{\mathcal A}}$. 
Denote it by $\alpha$. Fix elements $\sigma,\tau\in G$ such that 
$L$, $\sigma(L)$ and $\tau\sigma(L)$ are in general position, 
$\sigma|_L=\tau|_L$ and $\sigma^2|_L=id$. Then the generator 
$[1]-[\sigma]$ of ${\mathbb Q}[G/U_L]^{\circ}$ is fixed by 
$U_{L\sigma(L)}$, and therefore, \begin{multline*}
\alpha([1]-[\sigma])=f(x,y)\in\left({\mathcal A}(L\sigma(L))/
{\mathcal A}(k)\right)_{{\mathbb Q}}, \\ 
\alpha([\sigma]-[\tau\sigma])=f(y,z)\in\left({\mathcal A}
(\sigma(L)\tau\sigma(L))/{\mathcal A}(k)\right)_{{\mathbb Q}}, \\
\alpha([1]-[\tau\sigma])=f(x,z)\in\left({\mathcal A}(L\tau\sigma(L))/
{\mathcal A}(k)\right)_{{\mathbb Q}},\end{multline*} 
where $f(-,-)$ is a rational function and $x,y,z$ 
denote collections of elements in 
$L,\sigma(L),\tau\sigma(L)$, respectively. Then 
$f(x,y)+f(y,z)=f(x,z)$. Taking 
$\frac{\partial^2}{\partial x\partial z}$ in the case 
${\mathcal A}={\mathbb G}_a$, or 
$\frac{\partial}{\partial x}\log\cdot\frac{\partial}{\partial z}\log$ 
in the case ${\mathcal A}={\mathbb G}_m$, we get $f(x,z)=f(x)+g(z)$. 
As $\sigma([1]-[\sigma])=-([1]-[\sigma])$, one has $f(x,y)=-f(y,x)$, 
and thus, $f(x,y)=f(x)-f(y)$ for some $f(x)\in L$. 
Let $\beta\in G$ be such an element that $\beta f(x)=2f(x)$ (in the 
sense of the group law of ${\mathcal A}$). Then the image of 
$2\cdot[1]-[\beta]$ in $E$ is fixed by $G$ and projects to 
$1\in{\mathbb Q}$, so $E\cong{\mathbb Q}\oplus W_{{\mathcal A}}$. \qed 

\vspace{5mm}

It will follow from Propositions \ref{adm-i} and \ref{thick} 
that ${\rm Ext}^1_{{\mathcal S}m_G}({\mathbb Q},W)=0$ for any 
admissible representation $W$ of $G$ in the case $n=\infty$. (More 
generally, if $n=\infty$, $W_1,W_2\in{\mathcal I}^q_G$ and $W_1$ is 
projective in ${\mathcal I}^q_G$ then, using Lemma \ref{exa-i}, we get 
${\rm Ext}^1_{{\mathcal S}m_G}(W_1,W_2)={\rm Ext}^1_{{\mathcal S}m_G}
(W_1,W_2)={\rm Ext}^1_{{\mathcal I}^q_G}(W_1,W_2)=0$. Now take $q=0$, 
$W_1={\mathbb Q}$.) 

\vspace{5mm}

In the next example, we wish to show that ${\rm Ext}^1_{{\mathcal S}m_G} 
(W_{{\mathcal A}},W_{{\mathbb G}_m})\neq 0$ for any abelian variety 
${\mathcal A}$ over $k$. The $G$-module ${\rm Div}^{\circ}_{{\mathbb Q}}=
\lim\limits_{_U\longrightarrow}{\rm Div}_{{\rm alg}}(Y_U)_{{\mathbb Q}}$, 
introduced before Proposition \ref{pic-knot}, fits 
into the exact sequence $0\longrightarrow F^{\times}/k^{\times}
\longrightarrow{\rm Div}^{\circ}_{{\mathbb Q}}\longrightarrow
{\rm Pic}^{\circ}_{{\mathbb Q}}\longrightarrow 0$. By 
Proposition \ref{pic-knot}, any non-zero element of 
${\mathcal A}^{\vee}(k)_{{\mathbb Q}}$ determines an embedding of 
$W_{{\mathcal A}}$ into ${\rm Pic}^{\circ}_{{\mathbb Q}}$, 
thus inducing an extension of $W_{{\mathcal A}}$ by 
$W_{{\mathbb G}_m}$ inside ${\rm Div}^{\circ}_{{\mathbb Q}}$. 
This extension is non-split, since any generic $F$-point 
$x$ of ${\mathcal A}$, considered as an element of 
${\mathcal A}(F)_{{\mathbb Q}}$, identifies the space ${\rm Hom}_G(
{\mathcal A}(F)_{{\mathbb Q}},{\rm Div}^{\circ}_{{\mathbb Q}})$ 
with a subspace in 
$({\rm Div}^{\circ}_{{\mathbb Q}})^{{\rm Stab}_x}$ which is 
the same as $${\rm Div}_{{\rm alg}}({\mathcal A})_{{\mathbb Q}}^{
\langle\mbox{translations by torsion elements in 
${\mathcal A}(k)$}\rangle}=0.$$ 

\vspace{5mm}

For a smooth $G$-group ${\mathcal A}$ denote by $H^1_{{\mathcal S}m}
(G,{\mathcal A})$ the set of isomorphism classes of smooth 
$G$-torsors under ${\mathcal A}$, i.e., the set of classes 
of those $G$-torsors in $H^1(G,{\mathcal A})$ that become 
trivial on an open subgroup of $G$. 

In particular, $H^1(G,{\rm GL}_rF)$ is the set of isomorphism classes 
of semi-linear $r$-dimensional representations of $G$ over $F$, and 
$H^1_{{\mathcal S}m}(G,{\rm GL}_rF)$ is its subset of smooth ones. 

\begin{proposition} \label{torsors} If $n=\infty$ then 
$H^1_{{\mathcal S}m}(G,{\mathcal A}(F))
=\{\ast\}$ for any algebraic group ${\mathcal A}$ over $k$. \end{proposition} 
{\it Proof.} For any 1-cocycle $(a_{\sigma})$ presenting a class in 
$H^1_{{\mathcal S}m}(G,{\mathcal A}(F))$ there is an extension $L$ of $k$ in 
$F$ of finite type such that $a_{\xi}=1$ for any $\xi\in U_L$, so, as 
$a_{\sigma\tau}=a_{\sigma}\cdot\sigma a_{\tau}$, $a$ is a function on 
$G/U_L$, and $a_{\sigma}\in{\mathcal A}(L\sigma(L))$ for any $\sigma\in 
G$. 

For any $\tau:L\stackrel{/k}{\hookrightarrow}F$ in general 
position with respect to $L$ there is some $\sigma:L
\stackrel{/k}{\hookrightarrow}F$ such that $\sigma(L)$ and 
$\sigma\tau(L)$ are in general position with respect to $L$. 
Set $F'=\sigma(L)\tau(L)\sigma\tau(L)\subset F$. There is an 
$F'$-subalgebra $R$ of finite type in $LF'\subset F$ with the fraction 
field $LF'$ such that the fraction field of $R\bigcap L$ coincides 
with $L$ and $a_{\sigma},a_{\tau},a_{\sigma\tau}\in{\mathcal A}(R)$. There 
is a ring homomorphism $R\stackrel{s}{\longrightarrow}F'$ identical on 
$F'$ and inducing a homomorphism $R\bigcap L\longrightarrow k$, so 
$s(a_{\xi})\in{\mathcal A}(\xi(L))$ if $\xi\in\{\sigma,\tau,\sigma\tau\}$. 

Clearly, $\sigma a_{\tau}\in{\mathcal A}(F')$, so $s(a_{\sigma\tau})=
s(a_{\sigma}\cdot\sigma a_{\tau})=s(a_{\sigma})\cdot\sigma a_{\tau}$, 
and thus, $a_{\tau}=\sigma^{-1}s(a_{\sigma})^{-1}\cdot\sigma^{-1}
s(a_{\sigma\tau})$. This implies that $a_{\tau}=f_{\sigma}^{-1}\cdot
\tau f_{\sigma\tau}$, where $f_{\xi}:=\xi^{-1}s(a_{\xi})\in
{\mathcal A}(L)$ and $\xi\in\{\sigma,\sigma\tau\}$. 

In other words, for any $\tau:L\stackrel{/k}{\hookrightarrow}F$ in 
general position with respect to $L$ there exist $g_{\tau},h_{\tau}\in
{\mathcal A}(L)$ such that $a_{\tau}=g_{\tau}\cdot\tau h_{\tau}$. Any 
$\sigma:L\stackrel{/k}{\hookrightarrow}F$ in general position with 
respect to $L$ can be extended to $L\tau(L)
\stackrel{/k}{\hookrightarrow}F$ in such a way that $L,\sigma(L)$ and 
$\sigma\tau(L)$ will be in general position, so $g_{\sigma}^{-1}\cdot
a_{\sigma\tau}\cdot\sigma\tau h_{\tau}^{-1}=\sigma h_{\sigma}\cdot
\sigma g_{\tau}$ is an element of ${\mathcal A}(L\sigma\tau(L))$ and 
simultaneously of ${\mathcal A}(\sigma(L))$, i.e., $h_{\sigma}\cdot
g_{\tau}\in{\mathcal A}(k)$ for any $\sigma,\tau:L
\stackrel{/k}{\hookrightarrow}F$ in general position with respect to 
$L$. This means that $h_{\sigma}=b_{\sigma}\cdot f$ and $g_{\tau}=
f^{-1}\cdot c_{\tau}$, where $b_{\sigma},c_{\tau}\in{\mathcal A}(k)$ and 
$f\in{\mathcal A}(L)$, so $a_{\tau}=f^{-1}\cdot 
c_{\tau}b_{\tau}\cdot\tau f$. 

For any $\xi\in G$ there exist $\sigma$ and $\tau$ in general position 
with respect to $L$ such that $\xi=\sigma\tau$, and therefore, 
$a_{\xi}=(f^{-1}\cdot 
c_{\sigma}b_{\sigma}\cdot\sigma f)\cdot\sigma(f^{-1}\cdot c_{\tau}
b_{\tau}\cdot\tau f)=f^{-1}\cdot(c_{\sigma}b_{\sigma}c_{\tau}b_{\tau})
\cdot\sigma\tau f=f^{-1}\cdot d_{\xi}\cdot\xi f$ for some $d_{\xi}\in
{\mathcal A}(k)$, which means that the 1-cocycle $(a_{\xi})$ is 
cohomological to the image in $H^1_{{\mathcal S}m}(G,{\mathcal A}(F))$ of an 
element of $H^1_{{\mathcal S}m}(G,{\mathcal A}(k))={\rm Hom}_{{\mathcal S}m}
(G,{\mathcal A}(k))$, trivial by Theorem \ref{exactness}. \qed 

\vspace{4mm}

{\sc Remarks.} 1. If $n<\infty$ then there exist uncountably many 
semi-linear smooth representations of $G$ over $F$ of any finite 
dimension. Namely, any extension of coefficients from ${\mathbb Q}$ 
to $F$ of a non-trivial finite-dimensional smooth 
${\mathbb Q}$-representations of $G$ gives such semi-linear 
representations of $G$. (In fact, the natural map 
$H^1(G,{\mathcal A}(k))\longrightarrow H^1(G,{\mathcal A}(F))$ is 
injective for any algebraic group ${\mathcal A}$ over $k$. Otherwise, 
for a pair of 1-cocycles $(a_{\sigma})$ and $(a'_{\sigma})$ on $G$ 
with values in ${\mathcal A}(k)$ there would exist such 
$B\in{\mathcal A}(F)$ that 
$\sigma B=a_{\sigma}\cdot B\cdot{a'_{\sigma}}^{-1}$. If $f(B)\not\in 
k$ for some $f\in k({\mathcal A})$ defined at $B$ then there is $\sigma
\in G$ sending $f(B)$ out of the field of definition of $B$, so 
$\sigma B\neq A\cdot B\cdot A'$ for any $A,A'\in{\mathcal A}(k)$. \qed) 

2. If $n<\infty$ then there exist non-trivial finite-dimensional 
semi-linear smooth $F$-representations of $G^{\circ}$, e.g, 
$\Omega^q_{F/k}$ for any $1\le q\le n$. (Moreover, the class of 
$\det_F\Omega^q_{F/k}=\left(\Omega^n_{F/k}\right)^{\otimes_F
\binom{n-1}{q-1}}$ in $H^1(G^{\circ},\Omega^1_{F/k})$ under the 
$d\log$ map is non-trivial. Otherwise, for a generator 
$\omega=dx_1\wedge\dots\wedge dx_n$ of $\Omega^n_{F/k}$ there would 
exist some $\psi\in\Omega^1_{F/k}$ such that 
$\frac{d(\sigma\omega/\omega)}{\sigma\omega/\omega}=\sigma\psi-\psi$ 
for any $\sigma\in G^{\circ}$. For any $\sigma\in G^{\circ}$ with 
$\sigma x=Ax+B$ for an invertible $(n\times n)$-matrix $A$ over $k$ 
and a $k$-vector $B$ one has 
$\frac{d(\sigma\omega/\omega)}{\sigma\omega/\omega}=0$, and thus, 
$\sigma\psi-\psi=0$. But it is clear, that there are no non-zero 
1-forms invariant under all such $\sigma$'s, so $\psi=0$, and thus, 
$\tau\omega/\omega\in k$ for any $\tau\in G^{\circ}$. However, it is 
not the case if $\tau x_j=-x_j^{-1}$ for all $1\le j\le n$. \qed) 

\subsection{An example of $H_0$-calculation} 
\begin{lemma} \label{contr-van} Let $n=\infty$, $X$ a smooth 
projective variety over $k$, ${\mathcal A}$ an irreducible commutative 
algebraic group over $k$, and either $W=CH_0(X_F)^0_{{\mathbb Q}}$, 
or $W={\mathcal A}(F)_{{\mathbb Q}}$. Then for any open subgroup 
$U_L$ in $G$ one has $H_0(U_L,W)=0$. \end{lemma}
{\it Proof.} Let $K$ be the function field of $X$, or of 
${\mathcal A}$. The embeddings $\sigma:K\hookrightarrow F$ over 
$k$, in general position with respect to $L$ (i.e., with 
${\rm tr.deg}(\sigma(K)L/L)={\rm tr.deg}(K/k)$), form a 
single $U_L$-orbit. By Corollary \ref{cor-mov}, for any 
generic point $w:K\stackrel{/k}{\hookrightarrow}F$ in general 
position with respect to $L$ the ${\mathbb Q}$-space $W$ is generated 
by $\tau w-w$ for all $\tau\in U_L$, so $H_0(U_L,W)=0$. \qed

\section{The category ${\mathcal I}_G$} \label{I-G}
\begin{lemma} \label{sm-exh} Let $W$ be a smooth representation 
of $G$, and $L$ be an extension of $k$ in $F$. Then 
$W^{U_L}=\bigcup_{L_0\subseteq L}W^{U_{L_0}}$, where $L_0$ 
runs over extensions of $k$ of finite type. \end{lemma} 
{\it Proof.} For any $w\in W^{U_L}$ there is an extension $L_1$ 
of $k$ of finite type such that $w\in W^{U_{L_1}}$, so $w\in 
W^H$, where $H=\langle U_L,U_{L_1}\rangle$. 

Consider first the case ${\rm tr.deg}(L/k)<\infty$. Let $L_2$ be 
generated over $L_1$ by a transcendence basis of $L$ over $k$. 
Then $H\supseteq\langle U_L,U_{L_2}\rangle$ and 
$\overline{L_2\bigcap L}=\overline{L}$. One has the following 
evident inclusions $U_{\overline{L}}\subseteq
\langle U_L,U_{L_2}\rangle=:H'\subseteq U_{L_2\bigcap L}$. 
Consider the quotients $H'/U_{\overline{L}}=\langle 
{\rm Gal}(\overline{L}/L),U_{L_2}/U_{L_2\overline{L}}\rangle$ 
and $U_{L_2\bigcap L}/U_{\overline{L}}=
{\rm Gal}(\overline{L}/L_2\bigcap L)$. 

By the standard Galois theory (e.g., S.Lang, Algebra, Chapter VIII, 
\S1, Theorem 4), $U_{L_2}/U_{L_2\overline{L}}={\rm Gal}(L_2
\overline{L}/L_2)\cong{\rm Gal}(\overline{L}/L_2\bigcap\overline{L})$. 

According to Lemma \ref{2.0}, $\langle{\rm Gal}(\overline{L}/
L_2\bigcap\overline{L}),{\rm Gal}(\overline{L}/L)\rangle=
{\rm Gal}(\overline{L}/L_2\bigcap L)$, so $H'/U_{\overline{L}}
=U_{L_2\bigcap L}/U_{\overline{L}}$, and therefore, 
$H'=U_{L_2\bigcap L}$. Finally, $H\supseteq U_{L_0}$, where 
$L_0=L_2\bigcap L$. 

Now consider the case ${\rm tr.deg}(L/k)=\infty$. The group $H$ 
contains the subgroup $\langle U_{\overline{L}},U_{\overline{L_1}}
\rangle$, which coincides, by Proposition \ref{2.14}, with 
$U_{\overline{L}\bigcap\overline{L_1}}$. Let $L_2$ be generated over 
$L_1$ by a transcendence basis $S$ of 
$\overline{L_1}\bigcap\overline{L}$ over $k$ (in particular, 
$L_2\subset\overline{L_1}$). Then one has embeddings 
$L_1\subseteq L_2\subseteq\overline{L_1}$, and therefore, 
$\overline{L_1}=\overline{L_2}$, and thus, 
$\overline{L_1}\bigcap\overline{L}=\overline{L_2}\bigcap\overline{L}$. 
Similarly, $k(S)\subseteq L_2\bigcap\overline{L}\bigcap
\overline{L_1}\subseteq\overline{L}\bigcap\overline{L_1}$, 
and therefore, $\overline{L_2\bigcap(\overline{L}\bigcap
\overline{L_1})}=\overline{L}\bigcap\overline{L_1}=\overline{L_2}
\bigcap(\overline{L}\bigcap\overline{L_1})$. By Proposition 
\ref{2.14}, this implies that $\langle U_{L_2},
U_{\overline{L}\bigcap\overline{L_1}}\rangle=
U_{L_2\bigcap\overline{L}\bigcap\overline{L_1}}$, so $H$ contains 
$U_{L_3}$, where $L_3=L_2\bigcap\overline{L}$. 
Let $L_4$ be the minimal ${\rm Gal}(\overline{L}/L)$-invariant 
extension of $L_3$. As $L_3$ is a subfield of $\overline{L}$ 
finitely generated over $k$, $L_4$ is an extension of finite type 
of $L_3$. Then $\overline{L_4\bigcap L}=\overline{L_4}=
\overline{L_4}\bigcap\overline{L}$, so by Proposition \ref{2.14}, 
$\langle U_L,U_{L_4}\rangle=U_{L_0}$, where $L_0=L_4\bigcap L=
L_4^{{\rm Gal}(\overline{L}/L)}$. 

Finally, $H\supseteq U_{L_0}$, where $L_0$ is 
a subfield of $L$ of finite type over $k$. \qed 

\begin{corollary} \label{h-inv-fi-ge} 
Let $W$ be a smooth representation of $G$ such that 
$W^{U_{L_1}}=W^{U_{L_1(t)}}$ for any extension 
$L_1$ of $k$ of finite type and any $t\in F-\overline{L_1}$. 
Then $W^{U_L}=W^{U_{L'}}$ for any extension $L$ of $k$ and any 
purely transcendental extension $L'$ of $L$. \end{corollary}
{\it Proof.} By Lemma \ref{sm-exh}, $W^{U_{L'}}=
\bigcup_{L_0\subseteq L'}W^{U_{L_0}}$, where $L_0$ runs over 
extensions of $k$ of finite type. Let $L'=L(x_1,x_2,x_3,\dots)$ 
for some $x_1,x_2,x_3,\dots$ algebraically independent over $L$. 
Each $L_0$ of this type is a subfield in $L_1(x_1,\dots,x_N)$ for 
some $L_1\subseteq L$ of finite type over $k$ and some integer 
$N\ge 0$. Then $W^{U_{L'}}=\bigcup_{L_0\subseteq L'}W^{U_{L_0}}
\subseteq\bigcup_{L_1\subseteq L}W^{U_{L_1}}=W^{U_L}$, so 
$W^{U_{L'}}=W^{U_L}$. \qed 

\begin{lemma} \label{h-inv-alg-cl} Let $F'\subsetneqq F$ 
be an algebraically closed extension of $k$. Suppose that 
$W^{U_{F'}}=W^{U_{F'(x)}}$ for a smooth representation $W$ 
of $G$ and some $x\in F-F'$. Then $W^{U_L}=W^{U_{L(x)}}$ 
for any extension $L$ of $k$ inside $F'$. \end{lemma}
{\it Proof.} Fix a transcendence basis $S$ of $F'$ 
over $L$ and set $L_1=L(S)$. 

As the group $U_{L_1}$ is an extension of $G_{F'/L_1}$ by 
$G_{F/F'}$, and $U_{L_1(x)}$ is an extension of the group 
$G_{F'(x)/L_1(x)}=G_{F'/L_1}$ by $G_{F/F'(x)}$, one has 
$$W^{U_{L_1}}=\left(W^{G_{F/F'}}\right)^{G_{F'/L_1}}=\left(
W^{G_{F/F'(x)}}\right)^{G_{F'(x)/L_1(x)}}=W^{U_{L_1(x)}}.$$ 

For any $y\in S\cup\{x\}$ there is $\sigma\in U_L$ inducing a 
permutation of $S\cup\{x\}$ that transforms $(S\cup\{x\})-\{y\}$ to 
$S$. Such $\sigma$ induces an automorphism of $W^{U_{L_1(x)}}$ 
transforming $W^{U_{L((S\cup\{x\})-\{y\})}}$ to $W^{U_{L_1}}$, 
and therefore, the latter two spaces coincide. 

This implies that $W^{U_{L_1}}$ is fixed by the subgroup of $G$ 
generated by $U_{L((S\cup\{x\})-\{y\})}$ for all $y\in S\cup\{x\}$. 
By Lemma \ref{pur-tra}, this subgroup is dense in $U_L$, 
and therefore, $W^{U_{L_1}}=W^{U_L}$. \qed 

\vspace{4mm}

Let ${\mathcal I}_G$ be the full subcategory in ${\mathcal S}m_G$ consisting 
of those representations $W$ of $G$ for which $W^{G_{F/L}}=W^{G_{F/L'}}$ 
for any extension $L$ of $k$ in $F$ and any purely transcendental extension 
$L'$ of $L$ in $F$. For each integer $q\ge 0$ let ${\mathcal I}^q_G$ be 
the full subcategory in ${\mathcal I}_G$ consisting of those representations 
$W$ of $G$ for which $W^{G_{F/F'}}=0$ for any algebraically closed $F'$ with 
${\rm tr.deg}(F'/k)=q-1$. 

\begin{proposition} \label{adm-i} Any admissible representation of 
$G$ is an object in ${\mathcal I}_G$ if $n=\infty$. \end{proposition} 
{\it Proof.} Let $W$ be an admissible representation of $G$, $L$ an 
extension of $k$ in $F$ of finite type and $x,y\in F$ are algebraically 
independent over $L$. Then the finite-dimensional space $W^{U_L}$ is included 
into the finite-dimensional spaces $W^{U_{L(x)}}$ and $W^{U_{L(y)}}$; 
and the latter ones are included into the finite-dimensional space 
$W^{U_{L(x,y)}}$. As the group $U_{L(x+y,xy)}$ is an extension of 
the group $\{1,\alpha\}={\rm Gal}(L(x,y)/L(x+y,xy))$ (so $\alpha x=y$ 
and $\alpha y=x$) by $U_{L(x,y)}$, one has 
$W^{U_{L(x+y,xy)}}=\left(W^{U_{L(x,y)}}\right)^{\langle\alpha\rangle}$. 
As the subgroups $U_{L(x+y,xy)}$ and $U_{L(x,y)}$ are conjugated in $G$, 
the spaces $W^{U_{L(x+y,xy)}}$ and $W^{U_{L(x,y)}}$ are of the same 
dimension. This implies that $W^{U_{L(x+y,xy)}}=W^{U_{L(x,y)}}$, and 
thus, $\alpha$ acts trivially on $W^{U_{L(x,y)}}$. 

Notice, however, that $\alpha$ permutes $W^{U_{L(x)}}$ and 
$W^{U_{L(y)}}$, so $W^{U_{L(x)}}=W^{U_{L(y)}}$. By Lemma \ref{pur-tra}, 
the group generated by $U_{L(x)}$ and $U_{L(y)}$ is dense in 
$U_L$, and therefore, $W^{U_L}=W^{U_{L(x)}}$. \qed 

\begin{corollary} \label{adm-ab} The category of 
admissible representations of $G$ over $E$ is abelian. It is 
closed under extensions in ${\mathcal S}m_G(E)$. \end{corollary} 
{\it Proof.} It suffices to check that for any short exact sequence 
$$0\longrightarrow W_1\longrightarrow W_2\longrightarrow 
W_3\longrightarrow 0$$ of representations of $G$ 
with admissible $W_2$ the representations $W_1$ and $W_3$ 
of $G$ are also admissible. For any subextension $L$ 
of finite type over $k$ and a transcendence basis $t_1,t_2,t_3,\dots$ 
of $F$ over $L$ set $L'=L(t_1,t_2,t_3,\dots)$. Then the sequence 
$0\longrightarrow W_1^{U_{L'}}\longrightarrow W_2^{U_{L'}}
\longrightarrow W_3^{U_{L'}}\longrightarrow 0$ is exact. If $n<\infty$, 
we can assume that $F=\overline{L}$. As $W_2\in{\mathcal I}_G(E)$ in 
the case $n=\infty$, the middle term coincides with $W_2^{U_L}$, 
which is a finite-dimensional space, and therefore, so are the 
terms $W_1^{U_{L'}}$ and $W_3^{U_{L'}}$, containing $W_1^{U_L}$ and 
$W_3^{U_L}$, respectively. This implies that $W_1$ and $W_3$ are 
admissible. \qed 

\vspace{4mm}

{\sc Remark.} \label{no-proj} If $n=\infty$, the category of smooth 
representations of $G$ has such disadvantage that it has no non-zero 
projective objects. 

{\it Proof.} Let $W$ be a projective object in the category of 
smooth $E$-representations of $G$. Choose a system of generators 
$\{e_j\}_{j\in J}$ of $W$. This determines a surjection 
$$\bigoplus_{j\in J}E[G/U_{L_j}]\stackrel{\pi}{\longrightarrow}W,$$ 
where $U_{L_j}\subseteq{\rm Stab}_{e_j}$. Fix an element $i_0\in J$ 
and for each $j\in J$ fix an extension $L'_j$ of $L_j$ such that 
${\rm tr.deg}(L'_j/k)>{\rm tr.deg}(L_{i_0}/k)$. As $W$ is projective, 
the composition of $\pi$ with the surjection $\bigoplus_{j\in J}
E[G/U_{L'_j}]\longrightarrow\bigoplus_{j\in J}E[G/U_{L_j}]$ splits, 
and therefore, there is an element in $\bigoplus_{j\in J}E[G/U_{L'_j}]$ 
with the same stabilizer as $e_{i_0}$. However, as 
$E[G/U_{L'_j}]^{U_{L_{i_0}}}=0$, this implies that $e_{i_0}=0$, 
and thus, $W=0$. \qed 

\begin{lemma} \label{exa-i} The functor $H^0(G_{F/L},-):
{\mathcal I}_G\longrightarrow{\mathcal V}ect_{{\mathbb Q}}$ 
is exact for any extension $L$ of $k$ in $F$. ${\mathcal I}_G$ 
is closed under taking subquotients in ${\rm Sm}_G$. 

$\{{\mathcal I}^q_G\}_{q\ge 0}$ is a decreasing filtration of 
the category ${\mathcal I}_G$ by Serre subcategories.\footnote{A 
full subcategory of an abelian category ${\mathcal A}$ is called 
a Serre subcategory if it is stable under taking subquotients 
and extensions in ${\mathcal A}$.} \end{lemma}
{\it Proof.} By definition, the functors $H^0(G_{F/L},-)$ and 
$H^0(G_{F/L'},-)$ coincide on ${\mathcal I}_G$ for any purely
transcendental extension $L'$ of $L$. The group $G_{F/L'}$ 
is compact, if $F$ is algebraic over $L'$, so $H^0(G_{F/L'},-)$ 
is exact on ${\mathcal S}m_G$, and thus, its restriction to 
${\mathcal I}_G$ is also exact. 

Let $W\in{\mathcal I}^q_G$ and $W_1\subseteq W$ a 
subrepresentation of $G$. Then for any extension $L$ of $k$ in $F$ 
and any purely transcendental extension $L'$ of $L$ in $F$ one has 
$W_1^{G_{F/L}}=W_1\bigcap W^{G_{F/L}}=W_1\bigcap W^{G_{F/L'}}=
W_1^{G_{F/L'}}$, so $W_1\in{\mathcal I}_G$. Now it is clear that 
$W_1\in{\mathcal I}^q_G$. 

As $V^{G_{F/L}}\subseteq V^{G_{F/L'}}$ for any representation $V$ of 
$G$, to show that $V:=W/W_1\in{\mathcal I}_G$ we may suppose that $F$ 
is algebraic over $L'$. Then $V^{G_{F/L}}\subseteq V^{G_{F/L'}}
=W^{G_{F/L'}}/W_1^{G_{F/L'}}=W^{G_{F/L}}/W_1^{G_{F/L}}\subseteq
V^{G_{F/L}}$, so $V^{G_{F/L}}=V^{G_{F/L'}}$. It follows from the 
exactness of $H^0(G_{F/F'},-)$ that $V\in{\mathcal I}^q_G$. 

Let $0\longrightarrow W_1\longrightarrow E
\stackrel{\beta}{\longrightarrow}W_2\longrightarrow 0$ 
be a short exact sequence in ${\mathcal I}_G$ with 
$W_1,W_2\in{\mathcal I}^q_G$. Then for any algebraically closed $F'$ with 
${\rm tr.deg}(F'/k)=q-1$ the restriction of $\beta$ to $E^{G_{F/F'}}$ 
factors through $W_2^{G_{F/F'}}=0$, so $E^{G_{F/F'}}\subseteq 
W_1\bigcap E^{G_{F/F'}}=0$, and therefore, $E\in{\mathcal I}^q_G$. \qed 

\begin{lemma} \label{h-0-equi-cat} If $F'$ is a subfield in $F$ with 
${\rm tr.deg}(F'/k)=\infty$ then the functor 
$H^0(G_{F/\overline{F'}},-)$ from ${\mathcal S}m_G$ to 
${\mathcal S}m_{G_{\overline{F'}/k}}$ is an equivalence 
of categories (inducing an equivalence of ${\mathcal I}_G$ and 
${\mathcal I}_{G_{\overline{F'}/k}}$). The functor $H^0(G_{F/K},-)$ 
from ${\mathcal S}m_G$ to ${\mathcal V}ect_{{\mathbb Q}}$ is exact if 
and only if ${\rm tr.deg}(K/k)={\rm tr.deg}(F/k)(\le\infty)$. \end{lemma}
{\it Proof.} There exists a field isomorphism $\varphi:F
\stackrel{\sim}{\longrightarrow}\overline{F'}$ identical on $k$. 
Then $\varphi$ induces an isomorphism of topological 
groups $G_{\overline{F'}/k}\stackrel{\sim}{\longrightarrow}G$ by 
$\tau\longmapsto\varphi^{-1}\tau\varphi$ and an equivalence of 
the categories of representations of $G$ and representations of 
$G_{\overline{F'}/k}$ by $\pi\longmapsto\varphi^{\ast}\pi$, where 
$\varphi^{\ast}\pi(\tau)=\pi(\varphi^{-1}\tau\varphi)$. 

For any subfield $L\subset F$ finitely generated over $k$ there 
exists an element $\sigma\in G$ such that $\sigma|_L=\varphi|_L$. 
Let $W$ be a smooth representation of $G$. Then $\varphi$ and 
$\sigma$ induce the same isomorphism 
$W^{U_L}\stackrel{\sim}{\longrightarrow}W^{U_{\sigma(L)}}=
W^{U_{\varphi(L)}}$. Passing to the direct limit 
with respect to $L$, we get an isomorphism 
$W\stackrel{\sim}{\longrightarrow}W^{G_{F/\overline{F'}}}$. 
For any $\tau\in G_{\overline{F'}/k}$ and any $w\in W$ one has 
$\varphi\pi(\varphi^{-1}\tau\varphi)w=\tau\varphi w$ 
(since $\varphi$ is a limit of elements of $G$), i.e., 
$\varphi^{\ast}\pi\cong W^{G_{F/\overline{F'}}}$. 

Now, if ${\rm tr.deg}(K/k)=\infty$, then $H^0(G_{F/K},-)$ is the 
composition of exact functors $H^0(G_{F/\overline{K}},-)$ and 
$H^0({\rm Gal}(\overline{K}/K),-)$. Otherwise, if 
${\rm tr.deg}(K/k)<{\rm tr.deg}(L/k)$, then the only 
$G_{F/\overline{K}}$-invariant element in ${\mathbb Q}
[G/U_L]$ is zero, so $H^0(G_{F/K},-)$ transforms 
the surjection ${\mathbb Q}[G/U_L]\longrightarrow{\mathbb Q}$ to 
$0\longrightarrow{\mathbb Q}$. \qed 

\subsection{The functor ${\mathcal I}$} \label{func-I}
For a representation $M$ of $G$ define $N_jM$ as the subspace generated 
by the invariants $M^{G_{F/F_j}}$ for all subfields $F_j\subseteq F$ 
with ${\rm tr.deg}(F_j/k)=j$. Clearly, $N_jM\subseteq N_{j+1}M$, 
$M=\bigcup_{j\ge 0}N_jM$ if $M$ is smooth, 
\begin{itemize} \item 
$N_jM$ is the subrepresentation $M$ of $G$ in $M$ generated by 
$M^{G_{F/F_j}}$ for some algebraically closed $F_j$, 
\item restriction to $N_jM$ of each $G$-homomorphism 
$M\longrightarrow M'$ factors through $N_jM'$; 
\item $N_{i+j}(M_1\otimes M_2)\supseteq N_iM_1\otimes N_jM_2$. 
\end{itemize} 

\begin{proposition} \label{def-I} For any object $W\in{\mathcal S}m_G$ 
and any integer $q\ge 0$ there is its quotient ${\mathcal I}^qW\in
{\mathcal I}^q_G$ such that any $G$-homomorphism from $W$ to an object 
of ${\mathcal I}^q_G$ factors through ${\mathcal I}^qW$. The 
functor\footnote{cf. \cite{gm}, Chapter II, \S3.23} 
${\mathcal S}m_G\stackrel{{\mathcal I}^q}{\longrightarrow}{\mathcal I}^q_G$ 
given by $W\longmapsto{\mathcal I}^qW$ is right exact and 
${\mathcal I}^qW={\mathcal I}W/N_{q-1}{\mathcal I}W$. \end{proposition} 
{\it Proof.} Let $W'\in{\mathcal I}^q_G$. Any $G$-homomorphism  $W
\stackrel{\alpha}{\longrightarrow}W'$ factors through $\alpha(W)$, 
which is an object in ${\mathcal I}^q_G$, so we may assume that $\alpha$ 
is surjective. Let $L$ be an extension of $k$ in $F$ and $L'$ a purely 
transcendental extension of $L$ in $F$ over which $F$ is algebraic. 
As the functor $H^0(U_{L'},-)$ is exact on ${\mathcal S}m_G$, the morphism 
$\alpha$ induces a surjection $W^{U_{L'}}\longrightarrow(W')^{U_{L'}}$. 
As $(W')^{U_L}=(W')^{U_{L'}}$, the subgroup $U_L$ acts trivially 
on $(W')^{U_{L'}}$, and therefore, the subrepresentation $W_L=\langle
\sigma w-w~|~\sigma\in U_L,~w\in W^{U_{L'}}\rangle_G$ of $G$ is in the 
kernel of $\alpha$ (it is independent of $L'$ as all possible $L'$ 
form a single $U_L$-orbit and $\sigma\tau w-\tau w=(\sigma\tau)w-w-
(\tau w-w)$ for any $\tau\in U_L$ and any $w\in W^{U_{L'}}$; 
moreover, $W_L$ depends only on the $G$-orbit of $L$). This
implies that $\alpha$ factors through ${\mathcal I}W:=W/\sum_LW_L$. 

The representation ${\mathcal I}W$ of $G$ is smooth, so the map 
$W^{U_{L'}}\longrightarrow({\mathcal I}W)^{U_{L'}}$ induced by 
the projection is surjective, and therefore, one can lift any 
element $\overline{w}\in({\mathcal I}W)^{U_{L'}}$ to an element 
$w\in W^{U_{L'}}$. Then $\sigma\overline{w}-\overline{w}$ coincides 
with the projection of $\sigma w-w$ for any $\sigma\in U_L$. Note, 
that $\sigma w-w\in W_L$, so its projection is zero, and therefore, 
$\sigma\overline{w}=\overline{w}$ for any $\sigma\in U_L$. As 
$({\mathcal I}W)^{U_L}\subseteq({\mathcal I}W)^{U_{L'}}$, this means 
that $({\mathcal I}W)^{U_L}=({\mathcal I}W)^{U_{L'}}$, and thus, 
${\mathcal I}W\in{\mathcal I}_G$. 

We may further suppose that $W\in{\mathcal I}_G$. As $N_{q-1}W'=0$ and 
$N_{q-1}$ is functorial, $N_{q-1}W$ is in the kernel of $W
\stackrel{\alpha}{\longrightarrow}W'$, so $\alpha$ factors through 
${\mathcal I}^qW:=W/N_{q-1}W$. 

As the functor $H^0(G_{F/F'},-)$ is exact on ${\mathcal I}_G$, one 
has a short exact sequence $$0\longrightarrow(N_{q-1}W)^{G_{F/F'}}
\longrightarrow W^{G_{F/F'}}\longrightarrow({\mathcal I}^qW)^{G_{F/F'}}
\longrightarrow 0,$$ where $(N_{q-1}W)^{G_{F/F'}}=W^{G_{F/F'}}$ for 
any algebraically closed $F'$ of transcendence degree $q-1$ 
over $k$, so $({\mathcal I}^qW)^{G_{F/F'}}=0$, i.e., 
${\mathcal I}^qW\in{\mathcal I}^q_G$. 

As ${\rm Hom}_{{\mathcal I}^q_G}({\mathcal I}^qW,W')=
{\rm Hom}_{{\mathcal S}m_G}(W,W')$ for any $W\in{\mathcal S}m_G$ and 
$W'\in{\mathcal I}^q_G$, i.e., ${\mathcal I}^q$ is left adjoint to the 
identity functor ${\mathcal I}^q_G\hookrightarrow{\mathcal S}m_G$, it 
is right exact (cf., e.g., \cite{gm}, Chapter II, \S6.20). \qed 

\vspace{4mm}

{\sc Example.} \label{ell-ex} Let ${\mathcal A}$ be a one-dimensional
group scheme over $k$, $m\ge 1$ an integer, and $W=N_qW$ be a smooth 
representation of $G$, where $q\le n-1$. Then ${\mathcal I}
(S^m{\mathcal A}(F)\otimes W)=0$ if either $m$ is even, or 
if ${\mathcal A}={\mathbb G}_a$, or if 
${\mathcal A}={\mathbb G}_m$. In particular, the natural projection 
${\mathcal A}(F)^{\otimes N}_{{\mathbb Q}}\longrightarrow\bigwedge^N
{\mathcal A}(F)_{{\mathbb Q}}$ induces an isomorphism 
${\mathcal I}({\mathcal A}(F)^{\otimes N}_{{\mathbb Q}})\longrightarrow
{\mathcal I}(\bigwedge^N{\mathcal A}(F)_{{\mathbb Q}})$ if $n\ge N-1$. 

{\it Proof.} For any vector space $V$ and its proper subspace $V'$ 
the space $S^mV$ is spanned by the elements $x^{\otimes m}$ for all 
$x\in V-V'$, so $S^m{\mathcal A}(F)\otimes W$ is spanned by the elements 
$x^{\otimes m}\otimes w$ for all $w\in W$ with ${\rm Stab}_w\supseteq 
U_{k'}$, ${\rm tr.deg}(k'/k)<q+1$, and all $x\in{\mathcal A}(F)-
{\mathcal A}(k')$. If $m$ is even then the stabilizer of such $x^{\otimes m}
\otimes w$ contains the group $\{\sigma\in U_{k'}~|~\sigma x=\pm x\}$, 
which contains the subgroup $U_{k'(t)}$ for an element $t\in k(x):=
x(k({\mathcal A}))$ with quadratic extension $k(x)/k(t)$. If ${\mathcal A}$ 
is rational then the stabilizer of $x^{\otimes m}\otimes w$ contains the 
subgroup $U_{k'(x)}$. This implies that the image of $x^{\otimes m}
\otimes w$ in ${\mathcal I}(S^m{\mathcal A}(F)\otimes W)$ spans a trivial 
subrepresentation of $U_{k'}$. On the other hand, as $n\ge q+1>{\rm tr.deg}
(k'/k)$, there is $\sigma\in U_{k'}$ such that $\sigma x=2\cdot x$, 
and thus, $\sigma(x^{\otimes m}\otimes w)=2^m\cdot(x^{\otimes m}
\otimes w)$, which means that this trivial subrepresentation of $U_{k'}$ 
is zero, so ${\mathcal I}(S^m{\mathcal A}(F)\otimes W)=0$. \qed 

\vspace{4mm}

{\sc Remark.} ${\mathcal I}$ is not left exact, e.g., it transforms 
the injection $k\hookrightarrow F$ to $k\longrightarrow 0$. 

\begin{conjecture} \label{ss-conj} If $n=\infty$ then 
for any $j\ge 0$ and any object $W$ of ${\mathcal I}_G$ the 
$G$-module $gr^N_jW$ is semi-simple. \end{conjecture} 
This is evident when $j=0$, and deduced easily 
from Corollary \ref{str-lev-1} when $j=1$. 

\begin{corollary} \label{str-comp} If $n=\infty$ and Conjecture 
\ref{ss-conj} holds for any $0\le j\le q-1$ then the functor 
${\mathcal I}^q$ is exact on ${\mathcal I}_G$. This is equivalent to 
the strict compatibility of the filtration 
$N_0\subseteq\dots\subseteq N_{q-1}$ 
with morphisms in ${\mathcal I}_G$. \end{corollary} 
{\it Proof.} Let $0\longrightarrow
W_1\stackrel{\alpha}{\longrightarrow}W_2
\stackrel{\beta}{\longrightarrow}W_3\longrightarrow 0$ be a short 
exact sequence in ${\mathcal I}_G$. Then the first term of the induced 
exact sequence $$0\longrightarrow(W_1\bigcap N_{q-1}W_2)/N_{q-1}W_1
\longrightarrow{\mathcal I}^qW_1\longrightarrow{\mathcal I}^qW_2
\longrightarrow{\mathcal I}^qW_3\longrightarrow 0$$ measures deviation 
from the strict compatibility of $N_{q-1}$ with the morphism $\alpha$. 

To show that the filtration $N_{\bullet}$ on objects 
of ${\mathcal I}_G$ is strictly compatible with the morphisms, 
i.e., $\varphi(N_jW_1)=\varphi(W_1)\bigcap N_jW_2$ for 
any morphism $W_1\stackrel{\varphi}{\longrightarrow}W_2$ 
in ${\mathcal I}_G$, we proceed by induction on $j\ge 0$. 
This is clear when $j=0$. By Lemma \ref{exa-i}, 
$\varphi(W_1^{G_{F/F'}})=\varphi(W_1)^{G_{F/F'}}$ for any subfield 
$F'$ with ${\rm tr.deg}(F'/k)=j$, so $\varphi(N_jW_1)=
N_j\varphi(W_1)$, and thus, we may further suppose that $\varphi$ is 
an embedding: $W_1\subset W_2$. Suppose that $W_1\bigcap N_jW_2\neq 
N_jW_1$. Replace $W_1$ by $W'_1:=\frac{W_1\bigcap N_jW_2}{N_jW_1}$ and 
$W_2$ by $W'_2:=\frac{N_jW_2}{N_jW_1}$. Then $0\neq W'_1\subset 
W'_2=N_jW'_2$ and $N_jW'_1=0$. (By Proposition 
\ref{def-I}, $W'_1\subseteq{\mathcal I}^{j+1}W_1$, so by Lemma \ref{exa-i}, 
$W'_1\in{\mathcal I}^{j+1}_G$, which is equivalent to $N_jW'_1=0$.)

The induction assumption excludes the case $W'_1\subseteq
N_{j-1}W'_2$, so we may replace $W'_1$ by $W_1:=
\frac{W'_1}{N_{j-1}W'_2}$ and $W'_2$ by $W_2:=gr^N_jW'_2$. 
Then $0\neq W_1\subset W_2=N_jW_2$ and $N_jW_1=N_{j-1}W_2=0$. 

Assuming Conjecture \ref{ss-conj}, there is 
a morphism $W_2\stackrel{\xi}{\longrightarrow}W_1$ splitting the 
embedding $W_1\subset W_2$. However, $\xi(W_2)=\xi(N_jW_2)\subseteq 
N_jW_1=0$, leading to contradiction. \qed 

\vspace{4mm}

{\sc Remark.} The inclusion ${\mathbb Q}[G/U_L]^{\circ}\hookrightarrow
{\mathbb Q}[G/U_L]$ is an example of morphism of smooth $G$-modules which 
is {\sl not} strictly compatible with the filtration $N_{\bullet}$: 
$N_{{\rm tr.deg}(L/k)}{\mathbb Q}[G/U_L]^{\circ}$ coincides with 
$$\left\{\sum_{[\sigma]\in G/U_L}a_{\sigma}
[\sigma]~\bigg\vert~\begin{array}{l} \mbox{for any $F'$ with 
${\rm tr.deg}(F'/k)={\rm tr.deg}(L/k)$} \\ \mbox{one has 
$\sum_{\sigma(L)\subset\overline{F'}}a_{\sigma}=0$} \end{array} \right\}$$ 
which is different from ${\mathbb Q}[G/U_L]^{\circ}$, 
but ${\mathbb Q}[G/U_L]=N_{{\rm tr.deg}(L/k)}{\mathbb Q}[G/U_L]$.

\begin{corollary} \label{projectivity} Let $s\ge 0$ be an integer and 
$L$ be an extension of $k$ in $F$ of finite type with ${\rm tr.deg}(L/k)=q$. 
Then $C_L/N_{s-1}:={\mathcal I}^s{\mathbb Q}[G/U_L]$ is a projective object 
in ${\mathcal I}^s_G$ left orthogonal to ${\mathcal I}^{q+1}_G$. 
In particular, there are sufficiently many projective objects 
in ${\mathcal I}^s_G$ for any $s\ge 0$. \end{corollary} 
{\it Proof.} Any exact sequence $0\longrightarrow W_1
\longrightarrow W_2\longrightarrow W_3\longrightarrow 0$ 
in ${\mathcal I}^s_G$ gives an exact sequence $0\longrightarrow 
W_1^{U_L}\longrightarrow W_2^{U_L}\longrightarrow W_3^{U_L}
\longrightarrow 0$, where $W_j^{U_L}={\rm Hom}_G
({\mathbb Q}[G/U_L],W_j)={\rm Hom}_G(C_L/N_{s-1},W_j)$. \qed 

\begin{lemma} \label{no-0} For any $1\le n\le\infty$, any subfield 
$L_1\subset F$ of finite type over $k$, and any a unirational 
extension $L_2$ of $L_1$ in $F$ of finite type there is a natural 
isomorphism $C_{L_2}\stackrel{\sim}{\longrightarrow}C_{L_1}$. \end{lemma}
{\it Proof.} Let $x\in F-\overline{L_1}$ and let ${\mathbb Q}[G/U_{L_1(x)}]
\stackrel{\alpha}{\longrightarrow}{\mathcal I}{\mathbb Q}[G/U_{L_1(x)}]
=C_{L_1(x)}$ be the projection. For any $\sigma\in U_{L_1}$ one has 
$[\sigma U_{L_1(x)}]-[U_{L_1(x)}]=\sigma[U_{L_1(x)}]-[U_{L_1(x)}]\in 
W_{L_1}\subset{\mathbb Q}[G/U_{L_1(x)}]$ (see definition of ${\mathcal I}$ in 
the proof of Proposition \ref{def-I}), and thus, $\alpha$ factors through 
${\mathbb Q}[G/U_{L_1}]$, and therefore, through ${\mathcal I}
{\mathbb Q}[G/U_{L_1}]=C_{L_1}$, i.e., the surjection ${\mathbb Q}
[G/U_{L_1(x)}]\longrightarrow{\mathbb Q}[G/U_{L_1}]$ induces an 
isomorphism $C_{L_1(x)}\longrightarrow C_{L_1}$. 

One has $L_1\subseteq L_2\subseteq L_1(x_1,\dots,x_N)$ for 
some $x_1,\dots,x_N$ algebraically independent over $L_1$. Then 
the surjections $${\mathbb Q}[G/U_{L_1(x_1,\dots,x_N)}]\longrightarrow
{\mathbb Q}[G/U_{L_2}]\longrightarrow{\mathbb Q}[G/U_{L_1}]$$ 
induce surjections $C_{L_1(x_1,\dots,x_N)}\longrightarrow C_{L_2}
\longrightarrow C_{L_1}$, where the composition is an isomorphism. \qed 

\begin{lemma} Let $k\subseteq L\subset F'$ be subfields in $F$. 
Then $${\mathbb Q}[G/G_{F/F'}]\otimes_{{\mathbb Q}[G_{F'/k}]}
{\mathbb Q}[G_{F'/k}/G_{F'/L}]={\mathbb Q}[G/G_{\{F,F'\}/L}].$$ \end{lemma}
{\it Proof.} The module on the left coincides with 
${\mathbb Q}[G/G_{F/F'}]/\langle [\sigma\tau]
-[\sigma]~|~\sigma\in G, \tau\in G_{\{F,F'\}/L}\rangle$, 
which is the same as ${\mathbb Q}[G/G_{\{F,F'\}/L}]$. \qed

\begin{lemma} If $L$ is a finitely generated extension of $k$ 
and $F'$ is an algebraically closed subfield in $F$ with 
${\rm tr.deg}(F'L/L)={\rm tr.deg}(F'/F'\bigcap L)<\infty$ 
then there is a canonical isomorphism 
${\mathcal I}E[G/G_{\{F,F'\}/L}]\stackrel{\sim}{\longrightarrow}
C_L\otimes E$. \end{lemma} 
{\it Proof.} Let $t_1,\dots,t_N$ be a transcendence basis of $F'$ over 
$F'\bigcap L$, and thus, a transcendence basis of $F'L$ over $L$. 
Then the surjections $${\mathbb Q}[G/U_{L(t_1,\dots,t_N)}]\longrightarrow
{\mathbb Q}[G/G_{\{F,F'\}/L}]\longrightarrow{\mathbb Q}[G/U_L]$$ induce 
surjections $C_{L(t_1,\dots,t_N)}\longrightarrow{\mathcal I}{\mathbb Q}
[G/G_{\{F,F'\}/L}]\longrightarrow C_L$. By Lemma \ref{no-0}, their 
composition is an isomorphism, so both arrows are isomorphisms. \qed

\begin{proposition} \label{thick} ${\mathcal I}_G(E)$ is a Serre 
subcategory in ${\mathcal S}m_G(E)$ if $n=\infty$. \end{proposition} 
{\it Proof.} By Lemma \ref{exa-i}, the category ${\mathcal I}_G(E)$ 
is closed under taking subquotients in ${\mathcal S}m_G(E)$, so we have 
only to check that it is closed under extensions in ${\mathcal S}m_G(E)$. 

Let $0\longrightarrow W_1\longrightarrow W\longrightarrow W_2
\longrightarrow 0$ be an extension in ${\mathcal S}m_G(E)$ with 
$W_1,W_2\in{\mathcal I}_G(E)$. By Corollary \ref{h-inv-fi-ge} and Lemma 
\ref{h-inv-alg-cl}, it suffices to check that $W^{U_{\overline{L}}}=
W^{U_{L'}}$ for any extension $L$ of $k$ in $F$ of finite type 
and a purely transcendental extension $L'$ of $\overline{L}$ with 
$\overline{L'}=F$. As the functor $H^0(U_{L'},-)=
{\rm Hom}_{{\mathcal S}m_{U_{L'}}(E)}(E,-)$ is exact on 
${\mathcal S}m_G(E)$ and $W_1,W_2\in{\mathcal I}_G(E)$, this is 
equivalent to the vanishing of 
${\rm Ext}^1_{{\mathcal S}m_{U_{\overline{L}}}(E)}(E,-)$ on 
${\mathcal I}_G(E)$. As the forgetful functor ${\mathcal S}m_G(E)
\longrightarrow{\mathcal S}m_{U_{\overline{L}}}(E)$ induces 
${\mathcal I}_G(E)\longrightarrow{\mathcal I}_{U_{\overline{L}}}(E)$, 
one can replace $k$ with $\overline{L}$ and then it remains to show the 
vanishing of ${\rm Ext}^1_{{\mathcal S}m_G(E)}(E,-)$ on ${\mathcal I}_G(E)$. 

Let $0\longrightarrow W_1\longrightarrow W\longrightarrow E
\longrightarrow 0$ be an extension in ${\mathcal S}m_G(E)$ with 
$W_1\in{\mathcal I}_G(E)$. Choose some $e\in W$ projecting to 
$1\in E$. The stabilizer of $e$ contains an open subgroup 
$U_L$. Fix a maximal purely transcendental extension $L'$ of $k$ in 
$L$. Let $L''$ be a Galois extension of $L'$ containing $L$. Then 
$\frac{1}{[L'':L']}\sum_{\sigma\in{\rm Gal}(L''/L')}\sigma e$ is fixed 
by $U_{L'}$ and projects to $1\in E$, so $W$ is a sum of $W_1$ 
and a quotient of $E[G/U_{L'}]$. Restriction to $E[G/U_{L'}]^{\circ}$ 
of the projection $E[G/U_{L'}]\longrightarrow W$ factors through $W_1$ 
which is an object in ${\mathcal I}_G(E)$. 

Let $\sigma\in G$ be such an element that $L'$ and $\sigma(L')$ are 
in general position and $\sigma^2|_{L'}=id$. Then $[1]-[\sigma]$ is 
a generator of $E[G/U_{L'}]^{\circ}$, so there is a surjection 
$E[G/U_{L'\sigma(L')}]\longrightarrow E[G/U_{L'}]^{\circ}$ sending 
$[\tau]$ to $[\tau]-[\tau\sigma]$, and therefore, a surjection $E=
C_{L'\sigma(L')}\otimes E\longrightarrow{\mathcal I}E[G/U_{L'}]^{\circ}$. 
However, $\sigma$ changes sign of the generator $[1]-[\sigma]$ of 
$E[G/U_{L'}]^{\circ}$, and thus, ${\mathcal I}E[G/U_{L'}]^{\circ}=0$. 
This means that the projection $E[G/U_{L'}]\longrightarrow W$ factors 
through $E$, i.e., $W\cong W_1\oplus E$. \qed

\begin{conjecture} \label{C=CH} If $n=\infty$ then 
$CH_0(X_F)_{{\mathbb Q}}=C_{k(X)}$ for any smooth 
irreducible proper variety $X$ over $k$. \end{conjecture}
{\sc Remark.} This is true, by Corollary \ref{cl-pic} and 
Lemma \ref{no-0}, if $k(X)$ is unirational over a one-dimensional 
field. Another example is given by $k(X)=k(x_1,\dots,
x_m)^{\langle e_1e_2^2\cdots e_m^m\rangle}$, where 
$\sum_{j=1}^mx_j^d=1$, $m$ is odd and $d\in\{m+1,m+2\}$, 
$e_ix_j=\zeta^{\delta_{ij}}\cdot x_j$ for a primitive $d^{{\rm th}}$ root 
of unity $\zeta$. (The proof is the same as for $CH_0(X)={\mathbb Z}$.) 

\begin{proposition} \label{descr-inf-w} If $n=\infty$ then for any 
irreducible variety $X$ over $k$ the kernel of ${\mathbb Q}[\{k(X)
\stackrel{/k}{\hookrightarrow}F\}]\longrightarrow C_{k(X)}$ is the sum 
over all curves $y\in(k(X)\otimes_kF)_1$ of subspaces spanned by those 
linear combinations of generic (with respect to a field of definition 
of $y$) $F$-points of $\overline{\{y\}}$ that are linearly equivalent
to zero on any compactification of $\overline{\{y\}}$. \end{proposition} 
{\it Proof.} Let ${\mathfrak H}$ be the set of algebraically 
closed extensions $F_{\infty}$ of $k$ in $F$ such that 
${\rm tr.deg}(F/F_{\infty})=1$. Set $h_K:=h_{U_K}$. By Lemma 
\ref{h-inv-alg-cl} and Corollary \ref{h-inv-fi-ge}, the kernel 
$_{\infty}W$ of $W\longrightarrow{\mathcal I}W$ coincides with 
\begin{multline*} \sum_{F_{\infty}\in{\mathfrak H}}\langle\sigma w-w
~|~w\in W^{G_{F/F_{\infty}(t)}},\sigma\in G_{F/F_{\infty}},
t\in F-F_{\infty}\rangle_{{\mathbb Q}}\\
=\sum_{F_{\infty}\in{\mathfrak H}}\langle\sigma h_{F_{\infty}(t)}w-
h_{F_{\infty}(t)}w~|~w\in W,\sigma\in G_{F/F_{\infty}},t\in 
F-F_{\infty}\rangle_{{\mathbb Q}}. \end{multline*}

If $W={\mathbb Q}[G/U_L]$ this is the same as 
$$\sum_{F_{\infty}\in{\mathfrak H}}\langle\sigma h_{F_{\infty}(t)}\xi-
h_{F_{\infty}(t)}\xi~|~\xi:L\stackrel{/k}{\hookrightarrow}F,\sigma
\in G_{F/F_{\infty}},t\in F-F_{\infty}\rangle_{{\mathbb Q}}.$$ 
We may suppose that in this sum $\xi(L)\not\subset F_{\infty}$, as 
otherwise $\sigma h_{F_{\infty}(t)}\xi-h_{F_{\infty}(t)}\xi=0$. 
Then the pair $(\xi,F_{\infty})$ determines the $F_{\infty}$-curve 
$C^{F_{\infty},\xi}:={\bf Spec}((\xi\cdot id)(L\otimes_kF_{\infty}))$ 
on ${\bf Spec}(L\otimes_kF_{\infty})$. 

As for any  pair of elements $\sigma\in G_{F/F_{\infty}}$ 
and $\alpha\in G_{F/F_{\infty}(t)}$ the homomorphism 
$L\otimes_kF_{\infty}\stackrel{\sigma\alpha\xi\cdot id}%
{-\!\!\!-\!\!\!-\!\!\!\longrightarrow}F$ is the composition of 
$L\otimes_kF_{\infty}\stackrel{\xi\cdot id}{\longrightarrow}F$ with 
the automorphism $\sigma\alpha$, the kernels of $\sigma\alpha\xi\cdot id$ 
and of $\xi\cdot id$ coincide, and thus, $\sigma h_{F_{\infty}(t)}
\xi-h_{F_{\infty}(t)}\xi$ is an $F$-divisor on $C^{F_{\infty},\xi}$, 
which is a linear combination of generic $F$-points of 
$C^{F_{\infty},\xi}$. 

The triplet $(\xi,F_{\infty},t)$ determines the $F_{\infty}$-curve 
$C^{F_{\infty},\xi}_t:={\bf Spec}((\xi\cdot id)
(L\otimes_kF_{\infty}[t]))$ endowed with the projections 
${\mathbb A}^1_{F_{\infty}}\stackrel{T}{\longleftarrow} 
C^{F_{\infty},\xi}_t\stackrel{\beta}{\longrightarrow}C^{F_{\infty},\xi}$ 
induced by the inclusions $F_{\infty}[t]\subset(\xi\cdot id)(L\otimes_k
F_{\infty}[t])\supset(\xi\cdot id)(L\otimes_kF_{\infty})$. 
Now one can rewrite $\sigma h_{F_{\infty}(t)}\xi-h_{F_{\infty}(t)}\xi$ 
as $\beta_{\ast}T^{\ast}{\rm div}(\frac{T-\sigma t}{T-t})$, which is 
the divisor of a rational function on any compactification of 
$C^{F_{\infty},\xi}\times_{F_{\infty}}F={\bf Spec}((L\otimes_kF)/
{\mathfrak p})=:C_{{\mathfrak p}}$, where ${\mathfrak p}:=\ker(L\otimes_kF
\stackrel{\xi\otimes id}{\longrightarrow}F\otimes_{F_{\infty}}F)$. 

It remains to show that the divisor of a rational function $f$ on 
$C_{{\mathfrak p}}$, which is generic with respect to a field $F_0$ of 
definition of $C_{{\mathfrak p}}$, and independent of any 
compactification, belongs to $_{\infty}W$. As the map 
$Z^1(F_0(\xi(L))\otimes_{F_0}F)\longrightarrow{\mathcal I}W$ factors 
through ${\mathcal I}_{F_0}Z^1(F_0(\xi(L))\otimes_{F_0}F)$, where $\xi$ 
is a generic $F$-point of the model of $C_{{\mathfrak p}}$ over $F_0$, 
it follows from Proposition \ref{irat-1} that $f$ is sent to zero, 
i.e., ${\rm div}(f)\in~_{\infty}W$. \qed 

\subsection{Objects of ${\mathcal I}_G$ of level 1} \label{level-1} 
For a subfield $L$ in $F$ of finite type over $k$ denote by 
$Z^{{\rm rat}}_0(L\otimes_kF)$ the kernel of the natural 
projection ${\mathbb Q}[G/U_L]\longrightarrow CH_0(X_F)_{{\mathbb Q}}$ 
for a smooth proper model $X$ of $L$ over $k$. 

For any $W\in{\mathcal S}m_G$ one has a surjection 
$\bigoplus_{e\in W^{G_{F/F'}}}\langle e\rangle_G\longrightarrow 
N_1W$, where $F'$ is an algebraically closed extension of $k$ in $F$ 
with ${\rm tr.deg}(F'/k)=1$. This means that to describe the objects 
of ${\mathcal I}_G$ of level 1 it suffices to treat the case $W=\langle 
e\rangle_G$, where ${\rm Stab}_e\supseteq U_L$ with $L\cong k(X)$ 
for a smooth projective curve $X$ over $k$ of genus $g\ge 0$. Then 
$W$ is dominated by $C_L$. Let $J_X$ be the Jacobian of $X$. 
\begin{lemma} \label{cyclic-1} If $n=\infty$ then the  $G$-module 
$Z^{{\rm rat}}_0(L\otimes_kF)$ is generated by 
$w_N=\sum^N_{j=1}\sigma_j-\sum^N_{j=1}\tau_j$ for all $N>g$, where 
$(\sigma_1,\dots,\sigma_N;\tau_1,\dots,\tau_N)$ is a generic 
$F$-point of the fiber over $0$ of the map 
$X^N\times_kX^N\stackrel{p_N}{\longrightarrow}{\rm Pic}^{\circ}X$ 
sending $(x_1,\dots,x_N;y_1,\dots,y_N)$ to the class of 
$\sum^N_{j=1}x_j-\sum^N_{j=1}y_j$. \end{lemma} 
{\it Proof.} Let $\gamma_1,\dots,\gamma_s:L
\stackrel{/k}{\hookrightarrow}F$ and 
$\delta_1,\dots,\delta_s:L\stackrel{/k}{\hookrightarrow}F$ be generic 
points of $X$ such that $\sum^s_{j=1}\gamma_j-\sum^s_{j=1}\delta_j$ 
is the divisor of a rational function on $X_F$. 

We need to show that $\sum^s_{j=1}\gamma_j-\sum^s_{j=1}\delta_j$ 
belongs to the $G$-submodule in $Z^{{\rm rat}}_0(L\otimes_kF)$ 
generated by $w_N$'s. 

There is a collection $\alpha_1,\dots,\alpha_g:L
\stackrel{/k}{\hookrightarrow}F$ of generic points of $X$ such 
that the class of $\sum^s_{j=1}\gamma_j+\sum^g_{j=1}\alpha_j$ in 
${\rm Pic}^{s+g}X$ is a generic point. Then there is a collection 
$\xi_1,\dots,\xi_{s+g}:L\stackrel{/k}{\hookrightarrow}F$ 
of generic points of $X$ in general position such that 
$\sum^s_{j=1}\gamma_j+\sum^g_{j=1}\alpha_j-\sum^{s+g}_{j=1}\xi_j$ 
is divisor of a rational function on $X_F$ (so the same holds 
also for $\sum^s_{j=1}\delta_j+\sum^g_{j=1}\alpha_j-
\sum^{s+g}_{j=1}\xi_j$). We may, thus, suppose that 
$\delta_1,\dots,\delta_s$ are in general position. 

Fix a collection $\varkappa_{ij}$ of generic points of $X$ in general 
position, also with respect to $\gamma_1,\dots,\gamma_s$ and to 
$\delta_1,\dots,\delta_s$, for $1\le i\le g$ and $1\le j\le s$ 
such that the classes of $\gamma_1+\sum^g_{i=1}\varkappa_{i1},
\dots,\gamma_s+\sum^g_{i=1}\varkappa_{is}$ in ${\rm Pic}^{g+1}X$ 
are generic points in general position. Then one can choose 
a collection $\xi_{ij}$ of generic points of $X$ in general 
position for $0\le i\le g$ and $1\le j\le s$ such that 
$\gamma_j+\sum^g_{i=1}\varkappa_{ij}-\sum^g_{i=0}\xi_{ij}$ 
is divisor of a rational function on $X_F$ (so the same holds 
also for $\sum^s_{j=1}\sum^g_{i=0}\xi_{ij}-\left(\sum^s_{j=1}
\delta_j+\sum^s_{j=1}\sum^g_{i=1}\varkappa_{ij}\right)$). 
We may, thus, suppose that both $\gamma_1,\dots,\gamma_s$ 
and $\delta_1,\dots,\delta_s$ are in general position. 

Then there is a collection of generic points $\xi_1,\dots,\xi_s:
L\stackrel{/k}{\hookrightarrow}F$ such that the points 
$(\gamma_1,\dots,\gamma_s;\xi_1,\dots,\xi_s)$ and 
$(\delta_1,\dots,\delta_s;\xi_1,\dots,\xi_s)$ are generic 
on $p^{-1}_s(0)$. Then $\sum^s_{j=1}\gamma_j-\sum^s_{j=1}\xi_j$ 
and $\sum^s_{j=1}\delta_j-\sum^s_{j=1}\xi_j$ are divisors 
of rational functions on $X_F$. Clearly, such elements 
belong to the $G$-orbit of $w_s$. \qed 

\begin{lemma} \label{ifixed} The images of the generators 
$w_N=\sum^N_{j=1}\sigma_j-\sum^N_{j=1}\tau_j$ of 
$Z^{{\rm rat}}_0(L\otimes_kF)$ (from Lemma \ref{cyclic-1}) in 
$W:={\mathcal I}Z^{{\rm rat}}_0(L\otimes_kF)$ are fixed by $G$. 
\end{lemma} 
{\it Proof.} Set $M:=N+g$. Then for any generic point $\gamma$ of 
${\rm Pic}^M(X_F)$ in general position with respect to 
$\sigma_1,\dots,\sigma_N,\tau_1,\dots,\tau_N$ there is a collection 
$\alpha_1,\dots,\alpha_g$ of generic points of $X_F$, also in general 
position with respect to 
$\sigma_1,\dots,\sigma_N,\tau_1,\dots,\tau_N$, such that the class of 
$\sum_{j=1}^N\sigma_j+\sum_{j=1}^g\alpha_j$ coincides with $\gamma$. 

There exists an $M$-tuple $\xi=(\xi_1,\dots,\xi_M)$ of generic points 
such that both $2M$-tuples 
$(\sigma_1,\dots,\sigma_N,\alpha_1,\dots,\alpha_g;\xi)$ 
and $(\tau_1,\dots,\tau_N,\alpha_1,\dots,\alpha_g;\xi)$ 
are generic points of the irreducible variety $p_M^{-1}(0)$. 

The subfields $$L_{\xi}:=(\xi_1(L)\cdots\xi_M(L))^{{\mathfrak S}_M}, 
\quad L_{\sigma}:=(\sigma_1(L)\cdots\sigma_N(L)\alpha_1(L)\cdots
\alpha_g(L))^{{\mathfrak S}_M}$$ are isomorphic to the function field of 
the $M^{{\rm th}}$ symmetric power of $X$, and $L_{\sigma}L_{\xi}$ is 
isomorphic to the function field of the ${\mathfrak S}_M\backslash 
p_M^{-1}(0)/{\mathfrak S}_M$. As $M^{{\rm th}}$ symmetric power of $X$ is 
birational to the product of the Jacobian of $X$ and a projective 
space, the subfields $L_{\xi}$, $L_{\sigma}$, as well as 
$L_{\sigma}L_{\xi}$, are purely transcendental extensions of the 
subfield $\gamma(k({\rm Pic}^M(X)))$ in $F$. Clearly, the same is true 
for the subfields $L_{\tau}:=(\tau_1(L)\cdots\tau_N(L)\alpha_1(L)
\cdots\alpha_g(L))^{{\mathfrak S}_M}$ and $L_{\tau}L_{\xi}$. 

The elements $w_{\sigma}:=\sum_{j=1}^N\sigma_j+\sum_{j=1}^g\alpha_j-
\sum_{j=1}^M\xi_j$ and $w_{\tau}:=\sum_{j=1}^N\tau_j+
\sum_{j=1}^g\alpha_j-\sum_{j=1}^M\xi_j$ of $Z^{{\rm rat}}_0
(L\otimes_kF)$ are fixed by $U_{L_{\sigma}L_{\xi}}$ and by 
$U_{L_{\tau}L_{\xi}}$, respectively. Then the classes of $w_{\sigma}$ 
and of $w_{\tau}$ in $W$ are fixed by $U_{\gamma(k({\rm Pic}^M(X)))}$, 
and thus, so is their difference $w_N$. 

Fix a purely transcendental extension $L_0$ of $k$ in $F$ with 
${\rm tr.deg}(L_0/k)=g$ in general position with respect to 
$\sigma_1,\dots,\sigma_N,\tau_1,\dots,\tau_N$. Then $$L_0=
\bigcap_{\gamma\in G,~\gamma(k({\rm Pic}^M(X)))\supset L_0}
\gamma(k({\rm Pic}^M(X))),$$ so by Lemma \ref{2.0}, the subgroup in 
$G$ generated by $U_{\gamma(k({\rm Pic}^M(X)))}$ for such $\gamma$ 
contains $U_{L_0}$, and thus, the image of $w_N$ in $W$ is fixed by 
$U_{L_0}$. 

Finally, as $W\in{\mathcal I}_G$, one has $w_N\in W^{U_{L_0}}=W^G$. \qed

\begin{proposition} \label{irat-1} If $n=\infty$ then 
${\mathcal I}Z^{{\rm rat}}_0(L\otimes_kF)=0$. \end{proposition} 
{\it Proof.} By Lemma \ref{ifixed}, the images of 
the generators $w_N=\sum^N_{j=1}\sigma_j-\sum^N_{j=1}\tau_j$ of 
$Z^{{\rm rat}}_0(L\otimes_kF)$ (from Lemma \ref{cyclic-1}) in 
$W:={\mathcal I}Z^{{\rm rat}}_0(L\otimes_kF)$ are fixed by $G$. 
As $(\sigma_1,\dots,\sigma_N;\tau_1,\dots,\tau_N)$ and 
$(\tau_1,\dots,\tau_N;\sigma_1,\dots,\sigma_N)$ are both generic 
points of the irreducible variety $p_N^{-1}(0)$ and the generic 
points form a single $G$-orbit, there is an element $\beta\in G$ 
such that $\beta\sigma_j=\tau_j$ and $\beta\tau_j=\sigma_j$ 
for all $1\le j\le N$, so  $\beta w=-w$. As $W$ is generated by the 
images of $w_N$'s, which are fixed by $G$, we get $W=0$. \qed

\begin{corollary} \label{cl-pic} Let $X$ be a smooth projective 
curve over $k$, $L=k(X)$ be its function field and 
${\mathbb Q}[G/U_L]^{\circ}$ be the group of degree-zero 0-cycles. 
If $n=\infty$ then  ${\mathcal I}{\mathbb Q}[G/U_L]^{\circ}=
{\rm Pic}^{\circ}(X_F)_{{\mathbb Q}}$ 
and $C_L={\rm Pic}(X_F)_{{\mathbb Q}}$. \end{corollary} 
{\it Proof.} Using Proposition \ref{irat-1}, this can be done by 
applying the right-exact functor ${\mathcal I}$ to the short exact 
sequences $0\longrightarrow Z^{{\rm rat}}_0(L\otimes_kF)
\longrightarrow{\mathbb Q}[G/U_L]^{\circ}\longrightarrow
CH_0(X_F)^0_{{\mathbb Q}}\longrightarrow 0$ and 
$0\longrightarrow Z^{{\rm rat}}_0(L\otimes_kF)
\longrightarrow{\mathbb Q}[G/U_L]\longrightarrow
CH_0(X_F)_{{\mathbb Q}}\longrightarrow 0$. \qed 

\begin{corollary} \label{str-lev-1} If $n=\infty$ then any object of 
${\mathcal I}_G$ of level 1 is a direct sum of a trivial module and a 
quotient of direct sum of modules ${\mathcal A}(F)_{{\mathbb Q}}$ for 
some abelian varieties ${\mathcal A}$ over $k$ by a trivial 
submodule. \end{corollary} 
{\it Proof.} By Corollary \ref{cl-pic}, any object $W$ of 
${\mathcal I}_G$ of level 1 is a quotient of $\bigoplus_{X\in I}
{\rm Pic}(X_F)_{{\mathbb Q}}$ for a set $I$ of smooth projective 
curves over $k$. As there is a splitting 
${\rm Pic}(X_F)_{{\mathbb Q}}\cong{\mathbb Q}\oplus
{\rm Pic}^{\circ}(X_F)_{{\mathbb Q}}$, we get that $W$ is a quotient of 
$W^G\oplus\bigoplus_{{\mathcal A}\in J}{\mathcal A}(F)_{{\mathbb Q}}$ 
for a set $J$ of simple abelian varieties over $k$. 

In particular, 
$W/W^G$ is semi-simple, so there is a subset $J'\subseteq J$ such that 
the projection $\bigoplus_{{\mathcal A}\in J'}W_{{\mathcal A}}
\longrightarrow W/W^G$ is an isomorphism. 

Then $W^G\oplus
\bigoplus_{{\mathcal A}\in J'}{\mathcal A}(F)_{{\mathbb Q}}
\longrightarrow W$ is a surjection with trivial kernel. \qed

\begin{corollary} If $n=\infty$ then ${\mathcal A}(F)_{{\mathbb Q}}$ 
(resp., $W_{{\mathcal A}}$) is a projective object of ${\mathcal I}_G$ 
(resp., of ${\mathcal I}^1_G$) for any abelian $k$-variety 
${\mathcal A}$. \end{corollary} 
{\it Proof.} ${\mathcal A}(F)_{{\mathbb Q}}$ (resp., 
$W_{{\mathcal A}}$) is a direct summand of 
${\rm Pic}(X_F)_{{\mathbb Q}}$ (resp., $W_{J_X}$) for a smooth curve 
$X$ on ${\mathcal A}$, which is a projective object in 
${\mathcal I}_G$ (resp., in ${\mathcal I}^1_G$) by Corollary 
\ref{cl-pic} and Corollary \ref{projectivity}. \qed 

\vspace{4mm}

\label{filtra-2}
Define the following decreasing filtration on the objects of 
${\mathcal I}_G$: ${\mathcal F}^jW=\bigcap_{\varphi}\ker\varphi$, 
where $\varphi$ runs over the set of $G$-homomorphisms from $W$ to 
objects of ${\mathcal I}_G$ of level $j$. If $N_{\bullet}$ is strictly 
compatible with the morphisms of ${\mathcal I}_G$ then one can assume 
that $\varphi$'s are surjective. 

As $\ker(W\stackrel{id}{\longrightarrow}W)=0$, one has 
${\mathcal F}^{q+1}W=0$, if $W=N_qW$. 

\begin{corollary} \label{filt-2} If $n=\infty$ then 
$gr_{{\mathcal F}}^1C_{k(X)}={\rm Alb}X(F)_{{\mathbb Q}}$ 
and $C_{k(X)}\cong{\mathbb Q}\oplus
{\rm Alb}X(F)_{{\mathbb Q}}\oplus{\mathcal F}^2C_{k(X)}$ 
for any smooth proper $k$-variety $X$. \end{corollary} 
{\it Proof.} By Corollary \ref{str-lev-1}, ${\mathcal F}^1W=\bigcap
\ker(W\stackrel{\varphi}{\longrightarrow}{\mathbb Q}\oplus W')$, 
where $W'$ runs over all quotients of direct sum of modules 
${\mathcal A}(F)_{{\mathbb Q}}$ for some abelian varieties 
${\mathcal A}$ over $k$ by a trivial submodule. We may suppose that 
$\varphi$'s are surjective. As $C_{k(X)}$ is projective, any such 
$\varphi$ lifts to a homomorphism to the direct sum of 
${\mathcal A}(F)_{{\mathbb Q}}$. 

Let $C^{\circ}_{k(X)}=\ker(C_{k(X)}\stackrel{\deg}{\longrightarrow}
{\mathbb Q})$, so $C_{k(X)}\cong{\mathbb Q}\oplus  C^{\circ}_{k(X)}$. 
One has \begin{multline*} {\rm Hom}_G(C^{\circ}_{k(X)},{\mathcal A}
(F)_{{\mathbb Q}})={\rm Hom}_G(C_{k(X)},{\mathcal A}(F)_{{\mathbb Q}})/
{\rm Hom}_G({\mathbb Q},{\mathcal A}(F)_{{\mathbb Q}}) \\ 
=({\mathcal A}(k(X))/{\mathcal A}(k))_{{\mathbb Q}} 
=({\rm Mor}(X,{\mathcal A})/{\mathcal A}(k))_{{\mathbb Q}}\\ 
={\rm Hom}({\rm Alb}X,{\mathcal A})_{{\mathbb Q}}={\rm Hom}_G({\rm Alb}
X(F)_{{\mathbb Q}},{\mathcal A}(F)_{{\mathbb Q}}), \end{multline*} 
so any $G$-homomorphisms from $C^{\circ}_{k(X)}$ to any object 
of level $1$ factors through ${\rm Alb}X(F)_{{\mathbb Q}}$. 

The filtration splits, since ${\rm Alb}X(F)_{{\mathbb Q}}$ 
is projective. \qed 

\subsection{The inner ${\mathcal H}om$}
\begin{corollary} \label{red-gro} The inclusion ${\mathbb Q}[\{L(X)
\stackrel{/L}{\hookrightarrow}F\}]\subseteq{\mathbb Q}[\{k(X)
\stackrel{/k}{\hookrightarrow}F\}]$ induces a surjection of
$U_L$-modules ${\mathbb Q}[\{L(X)\stackrel{/L}{\hookrightarrow}F\}]
\longrightarrow C_{k(X)}$ for any extension $L$ of $k$ in $F$ with 
${\rm tr.deg}(F/L)=\infty$ and any irreducible $k$-variety 
$X$. \end{corollary} 
{\it Proof.} For any $\sigma:k(X)\stackrel{/k}{\hookrightarrow}F$ 
there is a generic curve $Y$ on $X_F$ defined over some $k'$ such that 
$\sigma$ is its $F$-point generic with respect to $k'$. Then the class 
of $\sigma$ in ${\mathcal I}_{k'}{\mathbb Q}[\{k'(Y)
\stackrel{/k'}{\hookrightarrow}F\}]$ can be presented by a linear
 combination of $F$-points of $Y$ generic with respect to 
$k'L$. As generic points of $Y$ are generic points of $X$, this means 
that ${\mathbb Q}[\{L(X)\stackrel{/L}{\hookrightarrow}F\}]
\longrightarrow C_{k(X)}$ is surjective. \qed

\vspace{5mm}

Propositions \ref{thick} and the following one suggest that 
the category ${\mathcal I}_G$ should be related to the category of 
{\sl effective} homological motives. 
\begin{proposition} \label{inner-hom} The inner ${\mathcal H}om$ 
functor\footnote{$(W_1,W_2)\longmapsto\lim\limits_{U\rightarrow}
{\rm Hom}_U(W_1,W_2)$, where $U$ runs over the set of open subgroups 
of $G$. Clearly, ${\rm Hom}_G(W_1,{\mathcal H}om(W_2,W_3))=
{\rm Hom}_G(W_1\otimes W_2,W_3)$ for any smooth $G$-modules 
$W_1,W_2,W_3$.} on ${\mathcal S}m_G$ induces the inner ${\mathcal H}om$ 
functor on ${\mathcal I}_G$ if $n=\infty$. Level of ${\mathcal H}om(W_1,W_2)$ 
does not exceed $q$ if $W_1,W_2=N_qW_2\in{\mathcal I}_G$ and $q\le 1$. 
\end{proposition} 
{\it Proof.} For any $W_1\in{\mathcal I}_G$ there is a collection $I$ of
irreducible varieties over $k$ and a surjection of $\bigoplus_{X\in I}
C_{k(X)}$ onto $W_1$, and thus, an inclusion of $G$-modules
$${\mathcal H}om(W_1,W_2)\hookrightarrow\prod_{X\in I}{\mathcal H}om(C_{k(X)},
W_2).$$ Clearly, for any collection $\{M_{\alpha}\}_{\alpha\in I}$ of 
objects of ${\mathcal I}_G$ any {\sl smooth} submodule in 
$\prod_{\alpha\in I}M_{\alpha}$ is also an object of 
${\mathcal I}_G$, and thus, to show that ${\mathcal H}om(W_1,W_2)\in
{\mathcal I}_G$ for any $W_1,W_2\in{\mathcal I}_G$ it suffices to check 
that ${\mathcal H}om(W_1,W_2)\in{\mathcal I}_G$ 
for any $W_2\in{\mathcal I}_G$ and any $W_1$ of type $C_{k(X)}$. 

By its definition, ${\mathcal H}om(W_1,W_2)$ is the union of 
${\rm Hom}_{U_L}(W_1,W_2)$ over all open subgroups $U_L$ in $G$. 

Fix an algebraically closed subfield $F'$ of $F$ such that 
${\rm tr.deg}(F/F')={\rm tr.deg}(F'/k)=\infty$ and 
an embedding $\sigma:k(X)\stackrel{/k}{\hookrightarrow}F$ 
in general position with respect to $F'$. 

As, by Corollary \ref{red-gro}, for any $L\subset F'$ with 
${\rm tr.deg}(F/L)=\infty$ there is a surjection of $U_L$-modules 
$W:={\mathbb Q}[\{L(X)\stackrel{/L}{\hookrightarrow}F\}]\longrightarrow 
W_1$, one has ${\rm Hom}_{U_L}(W_1,W_2)\subseteq{\rm Hom}_{U_L}(W,W_2)
\cong W_2^{U_{L\sigma(k(X))}}$, and thus, the $G_{F'/k}$-module 
${\mathcal H}om(W_1,W_2)^{G_{F/F'}}$ embeds into the $G_{F'/k}$-module 
$W_2^{U_{F'\sigma(k(X))}}=:W_0$. Here the group $G_{F'/k}=
G_{F'\sigma(k(X))/\sigma(k(X))}$ acts on $W_0$ as the quotient of 
$G_{F/\sigma(k(X))}$ by $G_{F/F'\sigma(k(X))}$. Clearly, $W_0\in
{\mathcal I}_{G_{F'/k}}$, so the smooth $G_{F'/k}$module 
${\mathcal H}om(W_1,W_2)^{G_{F/F'}}$ belongs to 
${\mathcal I}_{G_{F'/k}}$. By Lemma \ref{h-0-equi-cat}, the 
functor ${\mathcal S}m_G\stackrel{H^0(G_{F/F'},-)}%
{-\!\!\!-\!\!\!-\!\!\!-\!\!\!-\!\!\!-\!\!\!-\!\!\!\longrightarrow}
{\mathcal S}m_{G_{F'/k}}$ is an equivalence of categories, inducing 
an equivalence of ${\mathcal I}_G$ and ${\mathcal I}_{G_{F'/k}}$, so 
the object ${\mathcal H}om(W_1,W_2)$ belongs to ${\mathcal I}_G$. 

If $W_2$ is a trivial $G$-module then ${\mathcal H}om(C_{k(X)},W_2)=W_2$, 
so the representation ${\mathcal H}om(W_1,W_2)$ is a submodule of a 
trivial $G$-module, and thus, ${\mathcal H}om(W_1,W_2)$ is trivial itself. 

To show that ${\mathcal H}om(W_1,W_2)$ is of level 1, if $W_2=N_1W_1$, 
it is sufficient to treat the case $W_1=\bigoplus_{X\in I}
C^{\circ}_{k(X)}$ and $W_2=(\bigoplus_{j\in J}{\mathcal A}_j
(F)_{{\mathbb Q}})/\Lambda$ for some $\Lambda\subseteq
\bigoplus_{j\in J}{\mathcal A}_j(k)_{{\mathbb Q}}$. Then, using 
Corollaries \ref{filt-2} and \ref{red-gro}, we get that 
$${\rm Hom}_{U_L}({\mathbb Q}[\{L(X)\stackrel{/L}{\hookrightarrow}F
\}]^{\circ},W_2)={\rm Hom}_{U_L}({\rm Alb}X(F)_{{\mathbb Q}},W_2).$$ 
As this is independent of $L$, the $G$-module 
${\mathcal H}om(W_1,W_2)$ is trivial. \qed 

\vspace{5mm}

{\sc Remark.} The $G$-equivariant pairing ${\mathbb Q}[G/U]\otimes
{\mathbb Q}[G/U]\longrightarrow{\mathbb Q}$ given by 
$[\sigma]\otimes[\tau]\longmapsto 0$ if $[\sigma]\neq[\tau]$ and 
$[\sigma]\otimes[\sigma]\longmapsto 1$ defines an embedding of 
${\mathbb Q}[G/U]$ into its contragredient, so, unlike the objects of 
${\mathcal I}_G$ in the case $n=\infty$, for any $1\le n\le\infty$ 
there exist many smooth $G$-modules with non-trivial contragredients. 

\subsection{A tensor structure on ${\mathcal I}_G$} As Example after
Proposition \ref{def-I} on p.\pageref{ell-ex} shows, ${\mathcal I}_G$ 
is not closed under tensor products in ${\mathcal S}m_G$. 
We define $W_1\otimes_{{\mathcal I}}W_2$ by ${\mathcal I}(W_1\otimes W_2)$. 

This operation is not associative on ${\mathcal S}m_G$ as one can see 
from the following example. Let $W_j={\mathbb Q}[G/U_j]$ for some open 
subgroups $U_j=U_{L_j}$ in $G$, $1\le j\le N$, $N\ge 2$. 
Then one has $W_1\otimes\dots\otimes W_N={\mathbb Q}[\prod_{j=1}^NG/U_j]=
\bigoplus_{\tau\in G\backslash(\prod_jG/U_j)}{\mathbb Q}
[G\cdot(\tau_1,\tau_2,\dots,\tau_N)]$. 

Clearly, the representation 
${\mathbb Q}[G\cdot(\tau_1,\tau_2,\dots,\tau_N)]$ is isomorphic 
to the representation ${\mathbb Q}[G/(\bigcap_{j=1}^N\tau_jU_j\tau_j^{-1})]$, so 
$${\mathcal I}(W_1\otimes\dots\otimes W_N)\cong\bigoplus_{(\tau_j)\in 
G\backslash(\prod_jG/U_j)}C_{L_1\tau_2(L_2)\dots\tau_N(L_N)}.\footnote{More 
symmetrically, $W_1\otimes\dots\otimes W_N\cong\!\!\!\!\!
\bigoplus\limits_{x\in{\bf Spec}(L_1\otimes_k\dots\otimes_kL_N)}
{\mathbb Q}[\{k(x)\stackrel{/k}{\hookrightarrow}F\}]$, so 
${\mathcal I}(W_1\otimes\dots\otimes W_N)\cong\!\!\!\!\!
\bigoplus\limits_{x\in{\bf Spec}(L_1\otimes_k\dots\otimes_kL_N)}C_{k(x)}$.}$$ 
If $U_1=U_2=U_{k(x)}$, one has ${\mathcal I}W_1={\mathcal I}W_2={\mathbb Q}$, 
and therefore, $W_1\otimes_{{\mathcal I}}(W_2\otimes_{{\mathcal I}}
{\mathbb Q})=
W_1\otimes_{{\mathcal I}}{\mathcal I}W_2={\mathcal I}W_1={\mathbb Q}$. 

On the other hand, by Noether normalization, $(W_1
\otimes_{{\mathcal I}}W_2)\otimes_{{\mathcal I}}{\mathbb Q}={\mathcal I}
(W_1\otimes W_2)$ contains submodules isomorphic to $C_L$ 
for any field $L$ finitely generated over $k$ and with 
${\rm tr.deg}(L/k)=1$.

\begin{lemma} \label{kuenneth} If $n=\infty$ then for any finite 
collection of smooth irreducible proper $k$-varieties $X_1,\dots,X_N$ 
there is a canonical surjection of $G$-modules $C_{k(X_1\times_k
\cdots\times_kX_N)}\stackrel{{\mathcal I}(\alpha)}{\longrightarrow}
{\mathcal I}\left(C_{k(X_1)}\otimes\cdots\otimes C_{k(X_N)}\right)$. 

If $C_{k(X_1\times_k\cdots\times_kX_N)}=CH_0(X_1\times_k\cdots\times_k 
X_N)_{{\mathbb Q}}$ then ${\mathcal I}(\alpha)$ is an isomorphism. \end{lemma} 
{\it Proof.} It suffices to check that the canonical $G$-homomorphism 
${\mathbb Q}[\{k(X_1)\otimes_k\cdots\otimes_kk(X_N)
\stackrel{/k}{\hookrightarrow}F\}]\stackrel{\alpha}{\longrightarrow}
C_{k(X_1)}\otimes\cdots\otimes C_{k(X_N)}$, given by $\tau\longmapsto
\tau|_{k(X_1)}\otimes\cdots\otimes\tau|_{k(X_N)}$, is surjective, and 
its kernel is contained in the kernel of ${\mathbb Q}[\{k(X_1)\otimes_k
\cdots\otimes_kk(X_N)\stackrel{/k}{\hookrightarrow}F\}]\longrightarrow 
C_{k(X_1\times_k\cdots\times_kX_N)}$ (so that a canonical 
$G$-homomorphism $C_{k(X_1)}\otimes\cdots\otimes C_{k(X_N)}
\longrightarrow C_{k(X_1\times_k\cdots\times_kX_N)}$ is defined, and 
the composition $$C_{k(X_1\times_k\cdots\times_kX_N)}
\stackrel{{\mathcal I}(\alpha)}{\longrightarrow}
{\mathcal I}\left(C_{k(X_1)}\otimes\cdots\otimes C_{k(X_N)}\right)
\longrightarrow C_{k(X_1\times_k\cdots\times_kX_N)}$$ is the identity). 

We have to check that for any collection of generic points 
$\sigma_j\in X_j$, $1\le j\le N$, the class $\sigma$ of 
$\sigma_1\otimes\cdots\otimes\sigma_N$ in 
$C_{k(X_1)}\otimes\cdots\otimes C_{k(X_N)}$ is a linear combination of 
the images of generic points of $X_1\times_k\cdots\times_kX_N$. 

By induction on $j$ we show that $\sigma$ is a linear combination of 
elements of type $\sigma'_1\otimes\cdots\otimes\sigma'_N$, where 
$\sigma'_1,\dots,\sigma'_j$ are in general position. For $j=1$ there 
is nothing to prove. If $j>1$ there is a curve $Y$ on $X_j$ defined 
over a subfield $k'\subset F$ with ${\rm tr.deg}(k'/k)<\infty$ such 
that $\sigma'_j\in Y(F)-Y(\overline{k'})$. Clearly, for any $G$-module 
$W$ the canonical $G$-homomorphism $W\longrightarrow{\mathcal I}W$ factors 
through the $G_{F/\overline{k'}}$-homomorphism $W\longrightarrow
{\mathcal I}_{\overline{k'}}W$. Here ${\mathcal I}_{\overline{k'}}:
{\mathcal S}m_{G_{F/\overline{k'}}}\longrightarrow
{\mathcal I}_{G_{F/\overline{k'}}}$ denotes the same functor as 
${\mathcal I}$, but in the context of $G_{F/\overline{k'}}$-modules. 
The embedding $Y\hookrightarrow(X_j)_{k'}$ induces the 
$G_{F/\overline{k'}}$-homomorphism ${\mathbb Q}[Y(F)]\longrightarrow
{\mathbb Q}[X_j(F)]$, and therefore, using Corollary \ref{cl-pic}, the 
$G_{F/\overline{k'}}$-homomorphism $$Z_0(Y_{\overline{k'}})\oplus
{\rm Pic}(Y_F)_{{\mathbb Q}}\longrightarrow
\bigoplus_{x\in(X_j)_{\overline{k'}}}C_{\overline{k'}(x)}.$$ This 
implies that the image of $\sigma'_j\in Y(F)$ in ${\mathcal I}
{\mathbb Q}[X_j(F)]=\bigoplus_{x\in X_j}C_{k(x)}$, which is equal to 
$[\sigma'_j]\in C_{k(X_j)}$, coincides with the image of a linear 
combination of some 
generic points of $Y$ that are in general position with respect to 
$\sigma'_1,\dots,\sigma'_{j-1}$, i.e., of some points of 
$Y(F)-Y(k'')$, where $k''=\overline{k'\sigma'_1(k(X_1))\cdots 
\sigma'_{j-1}(k(X_{j-1}))}$. This completes the induction, so we 
may suppose that $\sigma_1,\dots,\sigma_j$ are in general position. 
Then there is some $\tau:k(X_1)\otimes_k\cdots\otimes_kk(X_N)
\stackrel{/k}{\hookrightarrow}F$ such that $\tau|_{k(X_j)}=\sigma_j$ 
for all $1\le j\le N$, thus implying that $\alpha$ is surjective. 

Besides, the representation $C_{k(X_1)}\otimes\cdots\otimes C_{k(X_N)}$ 
surjects onto the representation 
$CH_0((X_1)_F)\otimes\cdots\otimes CH_0((X_N)_F)_{{\mathbb Q}}$, and 
the latter one surjects onto $CH_0((X_1\times_k\dots\times_kX_N
)_F)_{{\mathbb Q}}=C_{k(X_1\times_k\dots\times_kX_N)}$. \qed 

\begin{corollary} \label{associativ} If $n=\infty$ and Lemma 
\ref{kuenneth} is true then $\otimes_{{\mathcal I}}$ is associative, 
the class of projective objects in ${\mathcal I}_G$ is stable under 
$\otimes_{{\mathcal I}}$, and $W_1\otimes_{{\mathcal I}}\cdots
\otimes_{{\mathcal I}}W_N={\mathcal I}\left(W_1\otimes\cdots\otimes 
W_N\right)$. \end{corollary} 
{\it Proof.} By Lemma \ref{kuenneth}, the class of $G$-modules of type 
$C_L$ is stable under $\otimes_{{\mathcal I}}$, and $\otimes_{{\mathcal I}}$ 
is associative on this class. As any projective object is a direct 
summand of a direct sum of $G$-modules of type $C_L$, the same holds 
for the class of projective objects in ${\mathcal I}_G$, and also $W_1
\otimes_{{\mathcal I}}\cdots\otimes_{{\mathcal I}}W_N={\mathcal I}\left(W_1
\otimes\cdots\otimes W_N\right)$ for projective $W_1,\dots,W_N$. 

Any object $W_j\in{\mathcal I}_G$ is the cokernel of a map ${\mathcal Q}_j
\stackrel{\alpha_j}{\longrightarrow}{\mathcal P}_j$, where ${\mathcal P}_j$ 
and ${\mathcal Q}_j$ are direct sums of $G$-modules of type $C_L$. This 
implies that $W_i\otimes_{{\mathcal I}}W_j$ is the cokernel of 
${\mathcal P}_i\otimes_{{\mathcal I}}{\mathcal Q}_j\oplus{\mathcal Q}_i\otimes
_{{\mathcal I}}{\mathcal P}_j\stackrel{id\otimes\alpha_j+\alpha_i\otimes id}
{-\!\!\!-\!\!\!-\!\!\!-\!\!\!-\!\!\!-\!\!\!\longrightarrow}{\mathcal P}_i
\otimes_{{\mathcal I}}{\mathcal P}_j$, and in general, $(\dots(W_1
\otimes_{{\mathcal I}}W_2)\otimes_{{\mathcal I}}\cdots)
\otimes_{{\mathcal I}}W_N$ is the cokernel of $$\bigoplus_{j=1}^N
{\mathcal I}({\mathcal P}_1\otimes\cdots\otimes{\mathcal Q}_j\otimes\cdots
\otimes{\mathcal P}_N)\stackrel{\sum_jid\otimes
\cdots\otimes\alpha_j\otimes\cdots\otimes id}{-\!\!\!-\!\!\!-\!\!\!-
\!\!\!-\!\!\!-\!\!\!-\!\!\!-\!\!\!-\!\!\!-\!\!\!-\!\!\!\longrightarrow}
{\mathcal I}({\mathcal P}_1\otimes\cdots\otimes{\mathcal P}_N).$$ 
Clearly, this is independent of rearrangements of the brackets. \qed 

\vspace{4mm}

{\sc Remarks.} 1. As it follows from Example on p.\pageref{ell-ex}, 
The form $W_1\otimes W_2\longrightarrow W_1\otimes_{{\mathcal I}}W_2$ 
can be degenerate. (If $W_1=E(F)_{{\mathbb Q}}$ for an elliptic 
curve $E$ over $k$, and $W_2=W_1/\Lambda$ for some subspace 
$0\neq\Lambda\subseteq E(k)_{{\mathbb Q}}$ then $W_1\otimes_{{\mathcal I}}W_1
={\mathcal I}(\bigwedge^2W_1)$ surjects onto $W_1\otimes_{{\mathcal I}}W_2$ 
with kernel dominated by $W_1\otimes\Lambda$. As the form is 
skew-symmetric, when lifted to $W_1\otimes W_1$, its left kernel 
contains $\Lambda$. \qed) 

2. The functor $E(F)\otimes_{{\mathcal I}}$ is not exact. (Applying it to 
$\Lambda\hookrightarrow W_1$, we get $W_1\otimes\Lambda\longrightarrow
{\mathcal I}(W_1^{\otimes 2})={\mathcal I}(\bigwedge^2W_1)$, with the kernel 
containing $S^2\Lambda$, since $S^2\Lambda$ is in the kernel of the 
composition $W_1\otimes\Lambda\hookrightarrow W_1^{\otimes 2}
\longrightarrow\bigwedge^2W_1$. $\Box$) 

In particular, if we denote by ${\rm Tor}^{{\mathcal I}}_{\bullet}(W,-)$ 
the left derivatives of the functor $W\otimes_{{\mathcal I}}$ then 
${\rm Tor}^{{\mathcal I}}_{\bullet}(W_1,W_2)\not\cong
{\rm Tor}^{{\mathcal I}}_{\bullet}(W_2,W_1)$. 
(As the functor $W\otimes_{{\mathcal I}}$ is right exact, any 
projective object of ${\mathcal I}_G$ is acyclic, cf., e.g., 
\cite{gm}, Ch.III, \S6.12., so if ${\rm Tor}^{{\mathcal I}}_{\bullet}
(W_1,W_2)\cong{\rm Tor}^{{\mathcal I}}_{\bullet}(W_2,W_1)$ 
then $${\rm Tor}^{{\mathcal I}}_j(E(F)_{{\mathbb Q}},W)\cong
{\rm Tor}^{{\mathcal I}}_j(W,E(F)_{{\mathbb Q}})=0$$ if $j>0$, i.e., 
the functor $E(F)\otimes_{{\mathcal I}}$ should be exact. \qed)

\section{Representations induced from the compact open subgroups} 
\label{ind-comp} In this section we give an example (Corollary 
\ref{example}) of a pair of essentially different open compact 
subgroups $U$ and $U'$ in $G$ with embeddings of $E$-representations 
$E[G/U]\hookrightarrow E[G/U']$ and $E[G/U']\hookrightarrow E[G/U]$ 
of $G$. This implies that $E[G/U]$ and $E[G/U']$ have the same irreducible 
subquotients. Proposition \ref{example-2} contains one more example 
of this phenomenon. However, it seems crucial for these examples that 
the primitive motives of maximal level of models of $F^U$ and $F^{U'}$ 
coincide (and trivial). 

But first, two general remarks. 

{\sc Remarks.} 1. {\bf Representation of $G/G^{\circ}$.} 
Let $U$ be a compact open subgroup in $G$. Then there is a 
surjection of the $E$-representation $E[G/U]$ of $G$ onto any 
irreducible $E$-representation of $G$ factorizing through 
the quotient $G/G^{\circ}$. 

2. {\bf Twists by 1-dimensional representations.} 
Let $\varphi$ be a homomorphism from $G/G^{\circ}$ to 
${\mathbb Q}^{\times}$. We consider $E[G/U](\varphi)$ as the same 
$E$-vector space as $E[G/U]$, but with the $G$-action 
$[\sigma]\stackrel{\tau}{\longmapsto}\varphi(\tau)\cdot[\tau\sigma]$. 
Then $\lambda_{\varphi}([\sigma]):=\varphi(\sigma)\cdot[\sigma]$ 
defines an isomorphism of representations $E[G/U]
\stackrel{\lambda_{\varphi}}{\longrightarrow}E[G/U](\varphi)$ of $G$. 

This implies that for any irreducible $E$-representation $W$ of $G$ 
the multiplicities of $W$ and $W(\varphi)$ in $E[G/U]$ coincide. 

\subsection{Purely transcendental extensions of quadratic extensions} 
\label{indomp} 
\begin{lemma} \label{1.2} Let $U$ and $U'$ be open compact 
subgroups in $G$ such that $U\bigcap U'$ is of index 2 in $U$: 
$U=(U\bigcap U')\bigcup\sigma(U\bigcap U')$, and 
$U'\bigcap\sigma U'\sigma\subseteq U$. Then $U$ is the 
only right $U$-coset in $UU'$. Equivalently, for any 
$\sigma_1,\sigma_2\in G$ if $\sigma_1UU'=\sigma_2UU'$ 
then $\sigma_1U=\sigma_2U$. \end{lemma}
{\it Proof.} Equivalence. If $U$ is the only right $U$-coset 
in $UU'$ and $\sigma_1UU'=\sigma_2UU'$ then 
$\sigma_2^{-1}\sigma_1U\subseteq UU'$, so 
$\sigma_2^{-1}\sigma_1U=U$, i.e., $\sigma_2^{-1}\sigma_1\in U$. 
Conversely, suppose that $\sigma_1UU'=\sigma_2UU'$ implies 
$\sigma_1U=\sigma_2U$. Then if $\sigma_1U\subseteq UU'$ one also 
has $\sigma_1UU'\subseteq UU'$ and, by the measure argument, 
$\sigma_1UU'=UU'$, and thus, $\sigma_1U=U$. 

Now suppose that $\sigma_1UU'=\sigma_2UU'$. As
$UU'=U'\bigcup\sigma U'$, one has either $\sigma_1U'=\sigma_2U'$ 
and $\sigma_1\sigma U'=\sigma_2\sigma U'$, or 
$\sigma_1U'=\sigma_2\sigma U'$ and $\sigma_1\sigma 
U'=\sigma_2U'$. The second case can be reduced to the first one 
by replacing $\sigma_2$ with $\sigma_2\sigma$ (as this does not 
change $\sigma_2U$). Now one has $\sigma^{-1}_1\sigma_2\in 
U'\bigcap\sigma U'\sigma\subseteq U$, and thus, 
$\sigma_1U=\sigma_2U$. \qed

{\sc Remark.} One obviously has $U\bigcap U'=
\sigma(U\bigcap U')\sigma\subseteq U'\bigcap(\sigma U'\sigma)$. 
Under assumptions of Lemma \ref{1.2}, this means 
$U'\bigcap(\sigma U'\sigma)=U\bigcap U'$. 

\begin{lemma} \label{1.3} Let $U$ and $U'$ be some open compact 
subgroups in 
$G$ such that $U\bigcap U'$ is of index 2 in $U$: $U=(U\bigcap U')
\bigcup\sigma(U\bigcap U')$. Suppose that for any integer $N\ge 
1$ and any collection $\tau_1,\dots,\tau_N\in U'\sigma-U$ one has 
$\tau_1\cdots\tau_N\neq 1$. Then the morphism of $E$-representations 
$E[G/U]\stackrel{[\xi]\mapsto[\xi\sigma]+[\xi]}%
{-\!\!\!-\!\!\!-\!\!\!-\!\!\!-\!\!\!-\!\!\!\longrightarrow}
E[G/U']$ of $G$ is injective. \end{lemma}
{\it Proof.} First, we check that $U'\bigcap(\sigma U'\sigma)
\subseteq U$. If $\tau\in U'\bigcap(\sigma U'\sigma)-U$ then 
$\tau^{-1}\in U'\bigcap(\sigma U'\sigma)-U$, so 
$1=\tau\tau^{-1}\in(U'-U)((\sigma U'\sigma)-U)=
(U'\sigma-U)(U'\sigma-U)$, which contradicts our assumption 
when $N=2$. 

Suppose that $\sum_{j=1}^Mb_j[\sigma_j]$ is in the kernel, i.e., 
$\sum_{j=1}^Mb_j([\sigma_j]+[\sigma_j\sigma])=0$, where $b_j\neq 0$, 
$\sigma_j$ are pairwise distinct as elements of $G/U$ and $M\ge 2$. 
Then, by Lemma \ref{1.2}, $\sigma_iUU'\neq\sigma_jUU'$ for $i\neq j$. 

One considers the graph whose vertices are the right $U'$-cosets 
in the union $\bigcup_{j=1}^M\sigma_jUU'$, and whose edges are 
the sets $\sigma_jUU'$ for all $1\le j\le M$ which join the 
vertices $\sigma_jU'$ and $\sigma_j\sigma U'$. There are at least 
2 edges entering to a given vertex, since otherwise this ``vertex'' 
is contained in the support of $\sum_{j=1}^Mb_j([\sigma_j]+
[\sigma_j\sigma])$, so there exists a simple cycle in the graph, 
say, formed by edges $\sigma_1UU',\dots,\sigma_sUU'$ for some 
$s\ge 3$, i.e., the intersection of the subsets $\sigma_iUU'$ and 
$\sigma_jUU'$ in $G$ is non-empty if and only if $|i-j|\in\{0,1,s-1\}$. 

We may suppose that for any $1\le j<s$ one has $\sigma_jUU'
\bigcap\sigma_{j+1}UU'=\sigma_jU'$, and $\sigma_1UU'
\bigcap\sigma_sUU'=\sigma_sU'$, and therefore, $\sigma_jU'=
\sigma_{j+1}\sigma U'$ for any $1\le j<s$, and $\sigma_sU'=
\sigma_1\sigma U'$. 

Then $\sigma_j^{-1}\sigma_{j+1}\in U'\sigma-U$ for any $1\le j<s$, 
and $\sigma_s^{-1}\sigma_1\in U'\sigma-U$. 

As $(\sigma_1^{-1}\sigma_2)\cdots
(\sigma_j^{-1}\sigma_{j+1})\cdots(\sigma_{s-1}^{-1}\sigma_s)
(\sigma_s^{-1}\sigma_1)=1$, we get contradiction. \qed

\begin{corollary} \label{example} Let $2\le n<\infty$ and $L''\subset F$ 
be a subfield finitely generated over $k$ with ${\rm tr.deg}
(F/L'')=1$. For some $u\in\sqrt{(L'')^{\times}}-(L'')^{\times}$ and 
some $t\in F$ transcendental over $L''$ set $L=L''(u,T)$, where 
$T=(2t-u)^2$, and $L'=L''(t)$. Then for $U=U_L$ and $U'=U_{L'}$ 
there exist embeddings $E[G/U']\hookrightarrow E[G/U]$ and 
$E[G/U]\hookrightarrow E[G/U']$. \end{corollary}
{\it Proof.}  One has $U\bigcap U'=U_{L''(t,u)}$, 
$U=(U\bigcap U')\bigcup(U\bigcap U')\sigma$, where 
$\sigma t=u-t$ and $\sigma|_{L''(u)}=id$, and 
$U'=(U\bigcap U')\bigcup(U\bigcap U')\tau$, where 
$\tau u=-u$ and $\tau|_{L''(t)}=id$. This implies that 
$U'\sigma-U=(U\bigcap U')\tau\sigma$ and $(\tau\sigma)^2u=u$, 
$(\tau\sigma)^2t=t+2u$, so for any $N\ge 1$ and any
$\tau_1,\dots,\tau_N\in U'\sigma-U$ one has 
$\tau_1\cdots\tau_N\neq 1$. Similarly, $U\tau-U'=
(U\bigcap U')\sigma\tau=(U'\sigma-U)^{-1}$, so for any 
$N\ge 1$ and any $\tau_1,\dots,\tau_N\in U\tau-U'$ one has 
$\tau_1\cdots\tau_N\neq 1$. It follows from Lemma \ref{1.3} that 
there exist embeddings $E[G/U']\hookrightarrow E[G/U]$ and 
$E[G/U]\hookrightarrow E[G/U']$. \qed 

\begin{proposition} \label{example-2} Fix an odd integer $m\ge 1$, and 
let $m-1\le n\le\infty$. Fix a collection $x_1,\dots,x_m$ of elements
of $F$ with the only relation $\sum_{j=1}^mx_j^d=1$ over $k$, where 
$d\in\{m+1,m+2\}$. Set $L''=k(x_1,\dots,x_m)$ and 
$L=(L'')^{\langle e_1e_2^2\cdots e_m^m\rangle}$, where 
$e_ix_j=\zeta^{\delta_{ij}}\cdot x_j$ for a primitive $d^{{\rm th}}$ root 
of unity $\zeta$. Let $L'$ be a maximal purely transcendental extension 
of $k$ in $L$. Then for $U=U_L$ and $U'=U_{L'}$ the $E$-representations 
$E[G/U]$ and $E[G/U']$ of $G$ have the same irreducible 
subquotients. \end{proposition}
{\it Proof.} As $E[G/U']$ embeds naturally into $E[G/U]$, any subquotient 
of the module $E[G/U']$ is a subquotient of $E[G/U]$. Let $W$ be the 
quotient of $E[G/U_{L''}]$ by the sum of the images of $E[G/U']$ 
under all possible $E[G]$-homomorphisms to $E[G/U_{L''}]$. 
As in the proof of $CH_0(Y_{[L]})={\mathbb Z}$,\footnote{Let $A$ be 
the image of ${\mathbb Q}[e_1,e_2,\dots,e_m]$ in 
${\rm End}_GW$. It is a semisimple algebra, so we want to show 
that $e_1e_2^2\cdots e_m^m\not\equiv 1$ modulo any maximal ideal in 
$A$. For this we note that $(L'')^{\langle e_{i_1}\cdots e_{i_l}
\rangle}$ is rational for any $1\le i_1<\dots<i_l\le m$, so 
$\sum_{j=0}^{d-1}(e_{i_1}\cdots e_{i_l})^j=0$. The assumptions on $m$ 
and $d$ imply that modulo any maximal ideal in $A$ the element 
$e_1e_2^2\cdots e_m^m$ is a non-trivial root of unity.} one checks 
that $W^{\langle e_1e_2^2\cdots e_m^m\rangle}=0$, and thus, 
$E[G/U]$ coincides with the sum of the images of 
$E[G/U']$ under all $E[G]$-homomorphisms to $E[G/U]$. \qed 

\vspace{4mm}

{\sc Remark.} Let $L$ be an extension of $k$ of finite type 
and of transcendence degree $q$ in $F$. Then, at least assuming some 
conjectures, any motivic $G$-module of level $<q$ is a subquotient of 
${\mathbb Q}[G/U_L]$ with infinite multiplicity. To see this, fix a 
transcendence basis $x_1,\dots,x_q$ of $L$ over $k$. Then there is 
a surjection ${\mathbb Q}[G/U_L]\longrightarrow\Omega^s_{F/k}$, given 
by $[1]\longmapsto x_{s+1}dx_1\wedge\dots\wedge dx_s$ for any $s<q$. 
Any motivic $G$-module of level $s$ is a submodule of 
$\Omega^s_{F/k}$ with infinite multiplicity. 

\appendix
\section{The centers of the Hecke algebras} \label{Hecke-center} 
\begin{lemma} \label{utv1} Let $K$ be a compact open subgroup 
in $G$. Let $\nu\in{\mathcal H}_E(K)$ be an element which is 
not a $E$-multiple of $h_K$. Then there exist elements 
$x_1,\dots,x_n\in F^K$ algebraically independent over $k$ 
such that $\nu h_U\not\in E\cdot h_U$, where 
$U=U_{k(x_1,\dots,x_n)}\supseteq K$. \end{lemma}
{\it Proof.} Let $\nu=\sum a_j\sigma_j h_K$, where the classes 
of $\sigma_j$ in $G/K$ are pairwise distinct. After subtracting 
a multiple of $h_K$, if necessary, we may suppose 
that $\sigma_j\not\in K$ for any $j$. Then the sets $\{x\in F^K~|~
\sigma_ix=\sigma_jx\}$ for $i\neq j$ and $\{x\in F^K~|~
\sigma_jx=x\}$ for any $j$ are proper $k$-subspaces in $F^K$, 
and therefore, there exist elements $x_1,\dots,x_n\in F^K$ 
algebraically independent over $k$ with $x_1$ outside their union. 
These conditions on $x_1,\dots,x_n$ imply that 
$\sigma_i|_{k(x_1,\dots,x_n)}\neq\sigma_j|_{k(x_1,\dots,x_n)}$ 
and $\sigma_j|_{k(x_1,\dots,x_n)}\neq id$ for any $i\neq j$. 

Set $U=U_{k(x_1,\dots,x_n)}$. Then the support of the element 
$\nu\ast h_U$ coincides with $\bigcup_j\sigma_j U$, which is 
not a subset in $U$, so  $\nu\ast h_U$ is not a multiple of $h_U$. 
\qed

\begin{lemma} \label{utv2} Let $U=U_{k(x_1,\dots,x_n)}$ 
for some $x_1,\dots,x_n$ algebraically independent over $k$, 
and let $\nu$ be a central element either in the Hecke algebra 
${\mathcal H}_E(U)$, or in the Hecke algebra ${\mathcal H}^{\circ}_E(U)$. 
Then $\nu\in E\cdot h_U$. \end{lemma} 
{\it Proof.} For any $\tau$ in the normalizer of $U$ one has 
$\nu(h_U\tau h_U)=\nu h_U\tau=\nu\tau\neq 0$ if $\nu\neq 
0$, and $(h_U\tau h_U)\nu=\tau h_U\nu=\tau\nu$. We may
suppose that the support of $\nu$ does not contain 1, i.e., 
${\rm Supp}(\nu)=\coprod_{\sigma\in S}U\sigma U$ for 
a finite subset $S$ in $G-U$. 

Let $H=\{\tau\in G~|~\tau|_{k(x_j)}\in{\rm Aut}(k(x_j)/k)
~\mbox{for all $1\le j\le n$}\}$, and for each $1\le j\le n$ let 
the subfield $L_j$ be generated over $k$ by 
$x_1,\dots,\widehat{x_j},\dots,x_n$. As $\nu$ is a central element 
in ${\mathcal H}(U)$, one has $\tau\nu\tau^{-1}=\nu$ for all 
$\tau\in H$. In particular, 
${\rm Supp}(\tau\nu\tau^{-1})={\rm Supp}(\nu)$, so each $\tau\in H$ 
induces a permutation of the set $S$ of double $U$-classes. 
The subgroup $U\subset H$ acts trivially on $S$, so the action of $H$ 
on $S$ factors through the quotient $H/U\cong
\left({\rm PGL}_2k\right)^n$. Any homomorphism from 
$\left({\rm PGL}_2k\right)^n$ to the permutation group of the set $S$, 
is trivial, since any element of ${\rm PGL}_2k$ is $(\# S)!$-th power 
of another element of ${\rm PGL}_2k$, and therefore, 
$U\tau\sigma\tau^{-1}U=U\sigma U$ for any $\sigma\in S$. In particular, 
$\tau\sigma\tau^{-1}x_j$ is in the finite set $U\sigma x_j$ for all 
$\tau\in H$; or, even more particularly, the set of fields 
$k(\tau\sigma x_j)$ for all $\tau\in H$ is finite. 

Fix some $j$. Suppose that $\sigma x_j\not\in\overline{k(x_j)}$ 
(this implies that $n>1$). Then there is $1\le s\le n$ different 
from $j$ such that $F$ is algebraic over $L_s(\sigma x_j)$. Set 
$H_j=\{\tau\in U_{L_s}~|~\tau|_{k(x_s)}\in{\rm Aut}(k(x_s)/k)\}$. 
Then for any $\tau\in H_j$ one has $\tau\sigma\tau^{-1}x_j=
\tau\sigma x_j$, so the $H_j$-orbit of $\sigma x_j$ should be finite, 
and thus, a subgroup of finite index in $H_j$ should be compact, 
so the group $H_j$ should be compact itself, which is false. 

As $U\tau\sigma\tau^{-1}U=U\sigma U$ is equivalent to 
$U\tau\sigma^{-1}\tau^{-1}U=U\sigma^{-1}U$, we get 
$\sigma^{\pm 1}x_j\in\overline{k(x_j)}$. If $\sigma^{\pm 1}x_j
\not\in k(x_j)$ then $k(\sigma^{\pm 1}x_j,x_j)/k(x_j)$ has a non-empty 
branch locus. The ${\rm PGL}_2k$-orbit of any point on ${\mathbb P}^1_k$ 
is infinite, so the ${\rm PGL}_2k$-orbit of the branch locus is also 
infinite, which means that the set of fields $k(\tau\sigma^{\pm 1}x_j)$ 
is infinite, unless $k(\sigma^{\pm 1}x_j)$ is a subfield in 
$k(x_j)$. Then $k(x_j)=\sigma k(\sigma^{-1}x_j)\subseteq
\sigma k(x_j)=k(\sigma x_j)\subseteq k(x_j)$. As the 
center of ${\rm PGL}_2k$ is trivial, this shows that 
$\sigma|_{k(x_j)}=id$. When varying $j$, we get $\sigma\in U$, 
contradicting our assumptions. \qed 

\begin{lemma} \label{utv4} Let $K$ be a compact subgroup in $G$. 
If $n<\infty$ and $\nu\in{\mathcal H}_E(K)-E\cdot h_K$ then there 
exists a compact open subgroup $U$ containing $K$ such that 
$\nu\ast h_U\not\in E\cdot h_U$. \end{lemma} 
{\it Proof.} There is some $\sigma$ in the support of $\nu$ 
outside of $K$, i.e., if $U'$ is an open compact subgroup in $G$ 
not containing $\sigma$ then there is an open subgroup $U\subseteq 
U'$ such that $\nu(\sigma U)\neq 0$, and therefore, the support 
of $\nu h_U$ contains $\sigma$, so it is non-empty and it 
does not coincide with $U$. \qed 

\vspace{4mm}

As a corollary of these statements we get 
\begin{theorem} \label{final} Let $K$ be a compact 
subgroup in $G$. Then the centers of the Hecke algebras 
${\mathcal H}_E(K)$ and ${\mathcal H}^{\circ}_E(K)$ coincide 
with $E\cdot h_K$ if $n<\infty$. \end{theorem} 
{\it Proof.} Clearly, for any pair of compact subgroups
$K\subseteq U$ the multiplication by $h_U$: 
$\nu\longmapsto\nu h_U$ gives homomorphisms of the centers 
$Z({\mathcal H}_E(K))\stackrel{h_U\ast}{\longrightarrow}
Z({\mathcal H}_E(U))$ and $Z({\mathcal H}^{\circ}_E(K))
\stackrel{h_U\ast}{\longrightarrow}Z({\mathcal H}^{\circ}_E(U))$. 
Then by Lemma \ref{utv4}, we may suppose that $K$ is open. 
By Lemma \ref{utv1}, we may further suppose that $K=U_{k(x_1,\dots,x_n)}$. 
Then, by Lemma \ref{utv2}, the centers of ${\mathcal H}_E(K)$ and 
${\mathcal H}^{\circ}_E(K)$ coincide with $E\cdot h_K$. \qed

\section{The case of positive characteristic} \label{subgr-p} 
In this appendix we show that all results of \S\ref{subgroups} 
and \S\ref{I-G} remain valid in the case of ${\rm char}(k)=p>0$. 

The topology on $G$ is the same as described in Introduction. The 
group $G$ is also Hausdorff, locally compact if $n<\infty$, and totally 
disconnected; the subgroups $G_{\{F,(F_{\alpha})_{\alpha\in I}\}/k}$ 
are closed in $G$, the fibers of the morphism of unitary semigroups 
$$\{\mbox{{\rm subfields in} $F$ {\rm over} $k$}\}
\longrightarrow\{\mbox{{\rm closed subgroups in} $G$}\}$$ 
given by $K\longmapsto{\rm Aut}(F/K)$ 
consist of subfields of $F$ with the same sets of perfect 
subfields containing them (with the same ``perfectization''), its 
image is stable under passages to sub-/sup-groups with compact 
quotients, and it induces bijections 
\begin{itemize} \item $\left\{\begin{array}{c} 
\mbox{{\rm perfect subfields} $K\subset F$} \\ 
\mbox{{\rm over} $k$ {\rm with} $F=\overline{K}$}\end{array}\right\}
\leftrightarrow\{\mbox{{\rm compact subgroups of} $G$}\}$; 
\item $\left\{\begin{array}{c} \mbox{{\rm perfect
subfields} $K$ {\rm of} $F$ {\rm minimal over}} \\ 
\mbox{{\rm subfields of finite type over} $k$ {\rm with} 
$F=\overline{K}$} \end{array} \right\}\leftrightarrow
\left\{\begin{array}{c} \mbox{{\rm compact open}} \\ 
\mbox{{\rm subgroups of} $G$}\end{array} \right\}$. 

The inverse correspondences are given by $G\supset H\longmapsto F^H$. 
\end{itemize} 

{\it Proof of Lemma \ref{2.0}.} We may suppose that $L$ and all 
$L_{\alpha}$ are perfect. Then, after replacing the reference to 
the Galois correspondence of \S\ref{subgroups} with the reference 
to the Galois correspondence of this appendix, the proof in 
\S\ref{subgroups} goes through. \qed 

{\it Proof of Lemma \ref{norm-int-phi-1}.} First, replace $L$ with 
its ``perfectization''. Let $\ell\in\{2,3\}-\{p\}$. Then 
any element $\tau$ in the common normalizer in $G$ of all closed 
subgroups of index $\le\ell$ in $U_L$ satisfies $\tau(L(f^{1/\ell}))
=L(f^{1/\ell})$ for all $f\in L^{\times}$. If $\tau\not\in U_L$ then 
there is an element $x\in L^{\times}$ such that $\tau x/x\neq 1$. 
Then $\tau x/x=y^{\ell}$ for some $y\in F^{\times}-\mu_{\ell}$. Set 
$f=x+\lambda$ for a variable $\lambda\in k$. By Kummer theory, $\tau 
f/f\in L^{\times\ell}$, and therefore, $L$ contains $L_0:=k\left(
y\frac{(x+\lambda y^{-\ell})^{1/\ell}}{(x+\lambda)^{1/\ell}}~|~
\lambda\in k\right)\subset\overline{k(x,y)}$. 

Now we come back to our original $L$ and replace $L_0$ with the 
subfield generated by appropriate $p$-primary powers of $y
\frac{(x+\lambda y^{-\ell})^{1/\ell}}{(x+\lambda)^{1/\ell}}$, where 
$p={\rm char}(k)$. As ${\rm tr.deg}(\overline{k(x,y)}/k)\le 2$, 
by our assumption on $L$, the subfield $L_0$ of $L$ should be finitely 
generated over $k$. But this is possible only if $y^{\ell}=1$, i.e., 
if $\tau\in U_L$. (To see this, one can choose 
a smooth model of the extension $L_0(x)/k(x,y)$ 
over $k$ and look at its branch locus.) \qed 

{\it Proof of Lemma \ref{uutuun}.} Let $\sigma\in H\bigcap U-\{1\}$ 
and $k'$ be the algebraic closure in $F$ of any subfield in 
$F^{\langle\sigma\rangle}$ with ${\rm tr.deg}(F/k')=1$. As the 
extension $F/F^{\langle\sigma\rangle}$ is abelian there is an element 
$x\in F-k'$ and an integer $N\ge 2$ such that $\sigma x\neq x$ 
and either $\sigma x^N=x^N$, or $\sigma(x^p-x)=x^p-x$. 
Then one has $\sigma(k')=k'$ and $\sigma(k'(x))=k'(x)$. 

The rest of the proof is the same as the last two 
paragraphs of the proof in \S\ref{subgroups}. \qed 

{\it Proof of Lemma \ref{opnor}} in \S\ref{subgroups} remains valid, 
after we replace $L'$ with its perfect closure in $LL'$, but we do 
not claim that the inclusion ${\rm PGL}_2k\hookrightarrow 
N_{G^{\circ}}U_{L'}/U_{L'}$ is an isomorphism. \qed 

{\it Proof of Lemma \ref{choice}.} We proceed by induction on $m$, 
the case $m=0$ being trivial. We wish to find $w_m\in F$ such that 
$w$ and $\xi w_m$ are algebraically independent over $k'$ generated 
over $k$ by $w_1,\dots,w_{m-1},\xi w_1,\dots,\xi w_{m-1}$. Suppose that 
there is no such $w_m$. Then for any $u\in F-\overline{k'}$ and 
any $v\in F-\overline{k'(\xi u)}$ one has the following vanishings 
in $\Omega^2_{k'(u,v,\xi u,\xi v)/k'}$: $du\wedge d\xi u=
dv\wedge d\xi v=0$, $d(u+v)\wedge d\xi(u+v)=0$, and $d(u+v^{\ell})
\wedge d\xi(u+v^{\ell})=0$ for any integer $\ell\ge 2$ prime to 
$p:={\rm char}(k)$. Applying the first two to the third, we get 
$\ell(v^{\ell-1}-\xi v^{\ell-1})dv\wedge d\xi u=0$, which means that 
either $\xi v^{\ell-1}=v^{\ell-1}$ for any $v\in 
F-\overline{k'(\xi u)}$, or $dv=0\in\Omega^1_{k'(v,\xi v)/k'}$ for all 
$v\in F-\overline{k'}$, or $d\xi u=0\in\Omega^1_{k'(u,\xi u)/k'}$ 
for all $u\in F-\overline{k'}$. In the first case $\xi v=v$ for 
any $v\in F$, i.e., $\xi=1$. 

If $\xi(\overline{k'})\neq\overline{k'}$ then there exists 
$u\in\overline{k'}$ such that $\xi u\in F-\overline{k'}$. Fix some 
$v\in F-(\overline{k'(\xi u)}\cup\xi^{-1}(\overline{k'}))$. Even if 
$\xi v\in\overline{k'(v)}$, the element $\xi(uv)$ does not belong to 
$\overline{k'(v)}=\overline{k'(uv)}$, i.e., $\xi(uv)$ and $uv$ are 
algebraically independent over $k'$. 

We may, thus, suppose that $\xi(\overline{k'})=\overline{k'}$. 
Replacing $\xi$ with $\xi^{-1}$, we reduce the case 
$d\xi u=0\in\Omega^1_{k'(u,\xi u)/k'}$ to the case 
$du=0\in\Omega^1_{k'(u,\xi u)/k'}$ for all $u\in F-\overline{k'}$. 
Let $P(X,Y^{p^s})$ be the minimal polynomial of $\xi u$ over 
$\overline{k'}[u]$ with maximal possible integer $s\ge 1$. 

Replacing $\xi$ by ${\rm Fr}^s\xi$, where ${\rm Fr}$ is the Frobenius 
automorphism, we get that $du\neq0\in\Omega^1_{k'(u,\xi u)/k'}$. 
As this implies $\xi=1$, we get contradiction, since no non-zero 
power of the Frobenius automorphism is identical on $k$. 
This shows that there exists desired $w_m\in F$. \qed 

{\it Proof of Proposition \ref{2.14}.} We may suppose that $L_1$ 
and $L_2$ are perfect. Then the proof in \S\ref{subgroups} goes 
through. \qed 

Proofs of Lemmas \ref{reduction}, \ref{inf-cl}, \ref{lll}, \ref{G-0-dec} 
and \ref{pur-tra}, Theorem \ref{exactness}, and Corollaries \ref{tri-cent}, 
\ref{liss-fin-inf}, \ref{ex-fin-inf} and \ref{ccc} given in 
\S\ref{subgroups} remain valid without any changes. 

\begin{acknowledgement}I would like to thank Uwe Jannsen 
for his encouraging interest in this work and many suggestions, and 
Maxim Kontsevich for a helpful discussion on Proposition \ref{nt}. 
The present form of the proof of Proposition \ref{adm-i} 
is a response to a question posed by Dmitry Kaledin. 
Several improvements were suggested by the referee. 
I am grateful to the I.H.E.S. for its hospitality and to 
the European Post-Doctoral Institute for its support during 
the period when the prime part of this work was done. 
\end{acknowledgement}

\end{document}